\documentclass[11pt,letterpaper]{amsart}
\usepackage[latin1]{inputenc}
\usepackage[francais,english]{babel}    % Pour du fran\c{c}ais
\usepackage{amssymb}
\usepackage{amsmath}
\usepackage{amsthm}
\usepackage{multicol}
\usepackage{graphicx}
\usepackage{color}
\usepackage[T1]{fontenc}
\usepackage{yfonts}
\usepackage{placeins}
\usepackage{array}
\usepackage{dsfont}
\usepackage{stmaryrd}
\usepackage{verbatim}
\usepackage{caption}
\usepackage{transparent}
\usepackage{bbold}
\usepackage{hyperref}
\usepackage{enumitem}
\usepackage{tikz}
\usepackage{float}
\usetikzlibrary{patterns}
\usetikzlibrary{arrows.meta}
\usetikzlibrary{decorations.markings}
\usetikzlibrary{matrix,arrows,calc,fit,cd,positioning,intersections,arrows.meta,braids}
\usepackage{pgfplots}
\usepackage{tikz-3dplot}

\newtheorem{thm}{Theorem}[section]
\newcommand{\bthm}{\begin{thm}}
\newcommand{\ethm}{\end{thm}}

\newtheorem{thmi}{Theorem}
\newcommand{\bthmi}{\begin{thmi}}
\newcommand{\ethmi}{\end{thmi}}

\newtheorem{cori}[thmi]{Corollary}
\newcommand{\bcori}{\begin{cori}}
\newcommand{\ecori}{\end{cori}}

\newtheorem{mthm}{Theorem}
\newcommand{\bmthm}{\begin{mthm}}
\newcommand{\emthm}{\end{mthm}}

\newtheorem{mcor}[mthm]{Corollary}
\newcommand{\bmcor}{\begin{mcor}}
\newcommand{\emcor}{\end{mcor}}
\newtheorem{mconj}[mthm]{Conjecture}
\newcommand{\bmconj}{\begin{mconj}}
\newcommand{\emconj}{\end{mconj}}
\newtheorem{mpro}[mthm]{Proposition}
\newcommand{\bmpro}{\begin{mpro}}
\newcommand{\empro}{\end{mpro}}

\newtheorem*{conj}{Conjecture}
\newcommand{\bconj}{\begin{conj}}
\newcommand{\econj}{\end{conj}}

\newtheorem*{question}{Question}
\newcommand{\bq}{\begin{question}}
\newcommand{\eq}{\end{question}}

\newtheorem*{thn}{Theorem}
\newcommand{\bthn}{\begin{thn}}
\newcommand{\ethn}{\end{thn}}

\newtheorem{exo}{Exercise}
\newcommand{\bex}{\begin{exo}}
\newcommand{\eex}{\end{exo}}

\newtheorem{sol}{Solution}
\newcommand{\bsol}{\begin{sol}}
\newcommand{\esol}{\end{sol}}

\newtheorem{pro}[thm]{Proposition}
\newcommand{\bpro}{\begin{pro}}
\newcommand{\epro}{\end{pro}}

\newtheorem{cor}[thm]{Corollary}
\newcommand{\bcor}{\begin{cor}}
\newcommand{\ecor}{\end{cor}}

\newtheorem{lem}[thm]{Lemma}
\newcommand{\blem}{\begin{lem}}
\newcommand{\elem}{\end{lem}}

\theoremstyle{definition}

\newtheorem{defi}[thm]{Definition}
\newcommand{\bdf}{\begin{defi}}
\newcommand{\edf}{\end{defi}}

\newtheorem*{defis}{Definition}
\newcommand{\bdfs}{\begin{defis}}
\newcommand{\edfs}{\end{defis}}

\newtheorem*{rmk}{Remark}
\newcommand{\brk}{\begin{rmk} \upshape}
\newcommand{\erk}{\end{rmk}}

\newtheorem*{rmks}{Remarks}
\newcommand{\brks}{\begin{rmks} \upshape}
\newcommand{\erks}{\end{rmks}}

\newtheorem*{exe}{Example}
\newcommand{\bexe}{\begin{exe} \upshape}
\newcommand{\eexe}{\end{exe}}

\newtheorem*{exes}{Examples}
\newcommand{\bexes}{\begin{exes} \upshape}
\newcommand{\eexes}{\end{exes}}

\newtheorem*{pre}{Proof}
\newcommand{\bp}{\begin{pre} \upshape}
\newcommand{\ep}{\hfill \qed \end{pre}}
\newcommand{\epp}{\end{pre}}

\newcommand{\beq}{\begin{eqnarray*}}
\newcommand{\eeq}{\end{eqnarray*}}

\newcommand{\beqn}{\begin{equation}}
\newcommand{\eeqn}{\end{equation}}

\newcommand{\ben}{\begin{enumerate}}
\newcommand{\een}{\end{enumerate}}
\newcommand{\bit}{\begin{itemize} \renewcommand{\labelitemi}{$\bullet$} \renewcommand{\labelitemii}{$\star$}}
\newcommand{\eit}{\end{itemize}}

\newcommand{\bfg}{
\begin{figure}[H]
\begin{center}}
\newcommand{\efg}{
\end{center}
\end{figure}
\FloatBarrier}

\newcolumntype{M}[1]{>{\raggedright}m{#1}}

\newcommand{\R}{\mathbb{R}}

\newcommand{\N}{\mathbb{N}}
\newcommand{\Z}{\mathbb{Z}}

\newcommand{\E}{\mathbb{E}}

\newcommand{\K}{\mathbb{K}}

\renewcommand{\SS}{\mathbb{S}}
\renewcommand{\H}{\mathbb{H}}

\newcommand{\F}{\mathbb{F}}

\newcommand{\bs}{\symbol{92}}

\newcommand{\ov}{\overline}

\renewcommand{\tilde}{\widetilde}

\newcommand{\CAT}{\operatorname{CAT}}

\newcommand{\Supp}{\operatorname{Supp}}

\renewcommand{\max}{\operatorname{max}}

\newcommand{\st}{\, | \,}

\newcommand{\ra}{\rightarrow}

\renewcommand{\geq}{\geqslant}
\renewcommand{\leq}{\leqslant}

\newcommand{\GL}{\operatorname{GL}}

\newcommand{\<}{\langle}
\renewcommand{\>}{\rangle}

\newcommand{\mk}{\medskip}

\makeatletter
\def\Ddots{\mathinner{\mkern1mu\raise\p@
\vbox{\kern7\p@\hbox{.}}\mkern2mu
\raise4\p@\hbox{.}\mkern2mu\raise7\p@\hbox{.}\mkern1mu}}
\makeatother

\makeatletter
\def\maketitles{%
  \null
  \thispagestyle{empty}%
  \vfill
  \begin{center}\leavevmode
    \normalfont
    {\LARGE \@title\par}%
    \vskip 1.2cm
    {\large \@author\par}%
    \vskip 1.2cm
    {\large \@subtitle\par}%
    \vskip 0.8cm
    {\large \@date\par}%
  \end{center}%
  \vfill
  \null
  \cleardoublepage
  }

\makeatother

\usetikzlibrary{arrows}

\title{New Garside structures and applications to Artin groups}
\author{Thomas Haettel and Jingyin Huang}
\date{\today}

\newcommand{\supp}{\operatorname{Supp}}

\begin{document}

\selectlanguage{english}

\maketitle
\begin{center}
\begin{minipage}{0.8\textwidth}
\textsc{Abstract.} 
Garside groups are combinatorial generalizations of braid groups which enjoy many nice algebraic, geometric, and algorithmic properties. In this article we propose a method for turning the direct product of a group $G$ by $\Z$ into a Garside group, under simple assumptions on $G$. This method gives many new examples of Garside groups, including groups satisfying certain small cancellation condition (including surface groups) and groups with a systolic presentation. 

Our method also works for an infinite class of Artin groups, leading to many new group theoretic, geometric and topological consequences for them. In particular, we prove new cases of $K(\pi,1)$-conjecture for some hyperbolic type Artin groups.
\end{minipage}
\end{center}

\let\thefootnote\relax\footnotetext{Thomas Haettel, thomas.haettel@umontpellier.fr, IMAG, Univ Montpellier, CNRS, France, and IRL 3457, CRM-CNRS, Universit\'{e} de Montr\'{e}al, Canada.\\
Jingyin Huang, Department of Mathematics, The Ohio State University, 100 Math Tower, 231 W 18th Ave, Columbus, OH 43210, U.S.}

%\tableofcontents

\section{Introduction}

The notion of Garside group originated in Garside's work on word problems and conjugacy problems for braid groups \cite{garside1969braid}. It turns out the key structure needed in Garside's argument also appears in more general groups later, notably in spherical Artin groups \cite{brieskorn_saito} and fundamental groups of complexified central simplicial arrangement complements \cite{deligne}. An axiomatic setup was provided in \cite{dehornoy_paris_gaussian,dehornoy2002groupes}, to study groups that share a similar structure as a class, called \emph{Garside groups}. Since then, other important classes of groups were proven to be Garside groups, including but not limited to some semi-direct products \cite{crisp2005representations}, some complex braid groups \cite{bessis_finite_complex_Kpi1,corran2011new,corran2015braid}, structure groups of non-degenerate, involutive and braided set-theoretical
solutions of the quantum Yang-Baxter equation \cite{chouraqui2010garside}, affine Artin groups of $\widetilde A_n$ and $\widetilde C_n$ \cite{digne_garside_An,digne_garside_Cn}, 
crystallographic braid groups \cite{mccammond_sulway}, rank three Artin groups \cite{delucchi_paolini_salvetti_rank_3} etc. Garside groups are also known to be closed under certain kinds of amalgamation products and HNN extensions \cite{picantin2022cyclic}, as well as Zappa--Sz\'ep products \cite{gebhardt2016zappa}. There are also several variations and generalizations of Garside groups, applying to more natural examples - we refer to the book \cite{garside} for a comprehensive review.

Garside group in this article always means Garside group of finite type, i.e. the Garside element has finitely many divisors. If it has infinitely many, then we will call it a \emph{quasi-Garside} group. Garside groups are known to enjoy a long list of nice geometric, group theoretic and topological properties - they are biautomatic \cite{charney1992artin,dehornoy2002groupes}, hence have solvable word problems and conjugacy problems, they are torsion-free, and admit finite $K(\pi,1)$ spaces \cite{charney2004bestvina,dehornoy2003homology}, they act geometrically on Helly graphs and injective metric spaces \cite{huang_osajda_helly,haettel_helly_kpi1}, hence satisfy the Farrell-Jones conjecture and coarse Baum-Connes conjecture \cite{helly_groups} etc. Also, Garside groups of finite or infinite type play a central role in the proof of the $K(\pi,1)$-conjecture for different classes of complex hyperplane complements, see \cite{deligne,bessis_finite_complex_Kpi1,mccammond_sulway,paolini_salvetti}.

\subsection{New Garside groups}
While Garside groups enjoy nice properties, they have a very strong algebraic constraint: since a power of the Garside element is central, they have infinite center. This explains why the list of known examples of Garside groups is somewhat limited. In this article, we propose a simple approach to circumvent this obstruction and use Garside theory to study some groups with possibly trivial center. 
Namely, given a group $G$, we will consider the direct product of $G$ with $\Z$ to artificially create a center, which will serve as the Garside axis. Then we work backward to find necessary conditions on $G$ to make sure $G\times \Z$ is a Garside group, leading to the following simple criterion.

\bmthm(=Theorem~\ref{thm:simple_criterion_U_garside})
\label{thm:U}
Let $U$ be a finite set, endowed with a positive partial multiplication (see Definition~\ref{def:partial mul}), and associated prefix order $\leq_L$ and suffix order $\leq_R$. Assume that the following requirements hold:
\bit
\item $(U,\leq_L)$ and $(U,\leq_R)$ are meet-semilattices.
\item For any $a,u,v,w \in U$ such that $a \cdot u,a \cdot v \in U$ and $w$ is the join for $\leq_L$ of $u$ and $v$, then $a \cdot w \in U$.
\item For any $a,u,v,w \in U$ such that $u \cdot a,v \cdot a \in U$ and $w$ is the join for $\leq_R$ of $u$ and $v$, then $w \cdot a \in U$.
\item For any $a,b,u,v \in U$ such that $a \cdot u,a \cdot v,b \cdot u,b \cdot v \in U$, either $a,b$ have a join for $\leq_R$, or $u,v$ have a join for $\leq_L$.
\eit
Consider the group $G_U$ given by the following presentation:
$$G_U = \<U \st \forall u,v,w \in U \mbox{ such that } u \cdot v =w, \mbox{ we have } uv=w\>.$$
Then the group $G_U \times \Z$ is a Garside group, with Garside element $(e,1)$.
\emthm

First note that if a group $G$ is such that $G \times \Z$ is a Garside group, then we can deduce an impressive list of consequences for $G$, see Theorem~\ref{thm:consequences_product_Z_garside} below.

This method can be applied to several classes of groups that we discuss in this article. For instance, it applies to some groups given by a $T(5)$ positive presentation, see Theorem~\ref{thm:presentation_girth_Garside} for a precise statement. In particular, we deduce the following nice consequence.

\bmcor(=Corollary~\ref{cor:surface_group_garside})
For any surface $S$ of finite type (possibly non-orientable), except the projective plane, $\pi_1(S) \times \Z$ is a Garside group.
\emcor

Another interesting family of examples comes from groups given by a presentation such that the associated flag Cayley complex is systolic, called \emph{systolic restricted presentation} by Soergel in~\cite{soergel_systolic}, where they are defined and studied. We refer to Definition~\ref{def:systolic}. Examples include some amalgams of Garside groups and some $2$-dimensional Artin groups. For these groups, we prove the following.

\bmcor(=Corollary~\ref{cor_systolic_garside})
Let $G$ denote a group with a systolic restricted presentation. Then $G \times \Z$ is a Garside group.
\emcor

Theorem~\ref{thm:U} also applies to a class of groups with positive square presentations in the sense of Definition~\ref{def:square presentation}, where a criterion for such groups times $\Z$ to be Garside is provided in Theorem~\ref{thm:square garside}. This applies to a subclass of groups arising from word labeled oriented graphs in the sense of \cite{harlander2015aspherical}, as well as some of the mock right-angled Artin groups defined in \cite{scott_mock}.

In order to motivate the study of groups $G$ for which $G \times \Z$ is a Garside group, we record here a list of direct consequences. We recall the definition of Garside groups in Section~\ref{subsec:garside}, and we recall various nonpositive curvature notions in Section~\ref{subsec:nonpositive_curvature}.

\bmthm \label{thm:consequences_product_Z_garside}
Assume that $G$ is a group such that $G \times \Z$ is Garside. Then the following hold:
\ben
\item The group $G \times \Z$ is Helly.
\item The group $G$ is torsion-free.
\item The group $G$ is CUB, more precisely it acts geometrically on a finite-dimensional metric space with a unique convex geodesic bicombing. Moreover, this metric space is a simplicial complex such that each simplex is equipped with a polyhedral norm.
\item The group $G$ acts geometrically on a weakly modular graph.
\item The group $G$ is biautomatic, and in particular:
\bit
\item The centralizer of a finite set of elements of $G$ is biautomatic.
\item $G$ has solvable word and conjugacy problems.
\item Any polycyclic subgroup of $G$ is virtually abelian, finitely generated and undistorted.
\item $G$ has quadratic Dehn function, as well as Euclidean higher dimensional Dehn function.
\eit
\item Any element of $G$ has rational translation length, with a uniformly bounded denominator.
\item The group $G$ has contractible asymptotic cones.
\item The group $G$ satisfies the Farrell-Jones conjecture with finite wreath products.
\item The group $G$ satisfies the coarse Baum-Connes conjecture.
\item The group ring $\K[G]$ satisfies Kaplansky's idempotent conjecture, if $\K$ is a field with characteristic zero.
\een
\emthm

We defer the references for this theorem to Section~\ref{subsec:proof_product_Z_garside}.

\subsection{Applications to Artin groups}
One of the main motivation for our work comes from Artin groups, see Section~\ref{subsec:Coxeter Artin} for basic definitions. To each Coxeter group, there is an associated Artin group, in the same fashion that the $n$-strand braid group is associated to the symmetric group of an $n$-element set.  Here is the precise definition. 

For every finite simple graph $\Gamma$ with vertex set $S$ and with edges labeled by some integer in $\{2,3,\dots\}$, one associates the Coxeter group $W_\Gamma$ with the following presentation:
$$W_\Gamma = \<S \st \forall \{s,t\} \in \Gamma^{(1)}, \forall s \in S, s^2=1, [s,t]_m=[t,s]_m \mbox{ if the edge $\{s,t\}$ is labeled $m$}\>,$$
where $[s,t]_m$ denotes the word $ststs\dots$ of length $m$. Such a graph $\Gamma$ may be called a \emph{Coxeter presentation graph}, emphasizing the fact that edges correspond to relations. We will also use $W_S$ to denote $W_\Gamma$.

\mk

We will also be using a graph closely related to $\Gamma$, the \emph{Dynkin diagram} $\Gamma_D$: it has the same vertex set $S$, with some edges labeled in $\{4,5,\dots,\infty\}$, with the following edges between vertices $s,t \in S$:
\bit
\item If there is an edge labeled $2$ between $s$ and $t$ in $\Gamma$, there is no edge between $s$ and $t$ in $\Gamma_D$.
\item If there is an edge labeled $3$ between $s$ and $t$ in $\Gamma$, there is an unlabeled edge between $s$ and $t$ in $\Gamma_D$.
\item If there is an edge labeled by $m \geq 4$ between $s$ and $t$ in $\Gamma$, there is the same edge between $s$ and $t$ in $\Gamma_D$ labeled $m$.
\item If there is no edge between $s$ and $t$ in $\Gamma$, there is an edge between $s$ and $t$ in $\Gamma_D$ labeled $\infty$.
\eit

\mk

The associated Artin group $A_\Gamma$ is defined by a similar presentation:
$$A_\Gamma = \<S \st \forall \{s,t\} \in \Gamma^{(1)}, [s,t]_m=[t,s]_m \mbox{ if the edge $\{s,t\}$ is labeled $m$}\>.$$
The groups $A_\Gamma$ are also called Artin-Tits groups, since they have been defined by Tits in~\cite{tits_artin}. We will also use $A_S$ to denote $A_\Gamma$.

\begin{rmk}[Remark on convention]
	When the Coxeter presentation graph $\Gamma$ of an Artin group is not a complete graph, we tend to work with Coxeter presentation graph. When the Coxeter presentation graph $\Gamma$ of an Artin group is a complete graph, we tend to work with the associated Dynkin diagram $\Gamma_D$, which is often more informative. For example, for the braid group, its Coxeter presentation graph is a complete graph, but its Dynkin diagram is a linear graph. Also, note that the Coxeter presentation graph of an Artin group is complete if and only if its Dynkin diagram does not have edges labeled by $\infty$.
\end{rmk}

General Artin groups are largely mysterious, and even basic questions such as the following are still wide open (see~\cite{godelle_paris}, \cite{charney_problems}, \cite{mccammond_mysterious}).
\ben
\item Are Artin groups torsion-free?
\item What is the center of Artin groups?
\item Do Artin groups have solvable word problem?
\item Is the natural hyperplane complement a classifying space for Artin groups (the $K(\pi,1)$ conjecture, see Section~\ref{subsec:Kpi1})?
\een

Note that a positive answer to the $K(\pi,1)$ conjecture implies that the corresponding Artin group is torsion-free, and also that its center is known (see~\cite{jankiewicz_schreve_center_kpi1}).

\mk

For Artin groups of spherical type, i.e. when the associated Coxeter group is finite, all these questions have a precise answer, which all rely on the existence of Garside structures. In fact, Artin groups of spherical type enjoy two different Garside structures: the standard one, associated with the longest element in the associated finite Coxeter group, and the dual one, associated with a Coxeter element. For an Artin group of non-spherical type, only the dual structure could be studied. In this case, the dual interval is always infinite, so one can only hope for a quasi-Garside structure, which has much fewer consequences. Nevertheless, it is known that for an Artin group of affine type $\tilde{A_n}$, $\tilde{C_n}$ or $\tilde{G_2}$ (\cite{digne_garside_An,digne_garside_Cn,mccammond_failure_lattice_property}), or for an Artin group of rank $3$ \cite{delucchi_paolini_salvetti_rank_3}, this dual structure turns the Artin group into a quasi-Garside group. In fact, for every Artin group of affine type, McCammond and Sulway manage to provide a natural embedding of the Artin group into a quasi-Garside crystallographic braid group, which is central in the proof of the $K(\pi,1)$ conjecture by Paolini and Salvetti (\cite{paolini_salvetti}).

\mk

However, even though a quasi-Garside structure might be sufficient to find classifying spaces, we already mentioned that a Garside structure on the direct product with $\Z$ is also useful, see Theorem~\ref{thm:consequences_product_Z_garside}. 

%In order to state our results concerning Artin groups, let us first recall some notations, we refer to Section~\ref{subsec:Coxeter Artin} for more details on our notations on Artin groups and their associated Coxeter groups. In particular, each Artin group or Coxeter group has a Coxeter presentation graph $\Gamma$, and a Dynkin diagram $\Lambda$. We will write $A_\Gamma$ (resp. $W_\Gamma$) to denote the Artin group (resp. Coxeter group) with Coxeter presentation graph $\Gamma$.

We say an Artin group is of \emph{cyclic type} if its Dynkin diagram is a cycle without $\infty$-labeled edges, and any proper parabolic subgroup is spherical. Note that cyclic type Artin groups all have complete Coxeter presentation graphs.
We refer to Table~\ref{tab:cyclic} for a complete list of cyclic type Artin groups, in terms of their Dynkin diagrams. In particular, it contains some Artin groups that are associated with certain Coxeter groups  acting on the hyperbolic spaces $\H^3$ or $\H^4$ - all of the four basic questions are open for these Artin groups.

\bmthm(=Proposition~\ref{prop_cyclic})
\label{thm:main1}
Suppose $A_\Gamma$ is of cyclic type. Then $A_\Gamma\times \Z$ is a Garside group.
\emthm

As we will see later (Corollary~\ref{cor:K(pi,1)}), Theorem~\ref{thm:main1} gives rise to new examples of Artin groups satisfying the $K(\pi,1)$-conjecture. We emphasize that an advantage of the method here is that it not only gives the $K(\pi,1)$-conjecture, but also it implies a long list of highly nontrivial algorithmic, geometric and topological consequences as in Theorem~\ref{thm:consequences_product_Z_garside}.

We can also treat a more general class of Artin groups which are obtained by gluing cyclic Artin groups and spherical Artin groups in the following way. To state the result, we will use the Coxeter presentation graphs.

Given a 4-cycle $\omega\subset \Gamma$ with consecutive vertices $\{x_i\}_{i=1}^4$, a pair of antipodal vertices in $\omega$ means either the pair $\{x_1,x_3\}$, or the pair $\{x_2,x_4\}$. A 4-cycle in $\Gamma$ has a \emph{diagonal} means it has a pair of antipodal vertices of $\omega$ which are connected by an edge in $\Gamma$. We say an induced subgraph of $\Gamma$ is of \emph{cyclic type} or \emph{spherical type} if the Artin group defined on this subgraph is of cyclic type or spherical type. An edge of $\Gamma$ is \emph{large} if it has label $\ge 3$. For an induced subgraph $\Lambda$ of $\Gamma$, let $\Lambda^\perp$ be the induced subgraph of $\Gamma$ spanned by vertices of $\Gamma\setminus \Lambda$ that commute with each vertex of $\Lambda$. We say the Coxeter presentation graph $\Gamma$ is a \emph{join} of two Coxeter presentation subgraphs $\Gamma_1$ and $\Gamma_2$ if $\Gamma$ is a join of $\Gamma_1$ and $\Gamma_2$ as graphs (i.e. $\Gamma$ is obtained from a disjoint union of $\Gamma_1$ and $\Gamma_2$ by adding a single edge between each vertex of $\Gamma_1$ and each vertex of $\Gamma_2$), and each vertex of $\Gamma_1$ commute with every vertex of $\Gamma_2$. We say a Coxeter presentation graph $\Gamma$ is spherical (resp. of cyclic type) if $A_\Gamma$ is a spherical Artin group (resp. an Artin of cyclic type). 

\bmthm(=Theorem~\ref{thm_main})
\label{thm:intro garside}
Let $\Gamma$ be a Coxeter presentation graph such that
\bit
\item each complete subgraph of $\Gamma$ is a join of a cyclic type graph and a spherical type graph (we allow one of the join factors to be empty);
\item for any cyclic type induced subgraph $\Lambda\subset \Gamma$, $\Lambda^{\perp}$ is spherical.
\eit
We assume in addition that there exists an orientation of all large edges of $\Gamma$ such that
\ben
\item the orientation restricted to each cyclic type subgraph of $\Gamma$ gives a consistent orientation on the associated circle;
\item if $\omega$ is a 4-cycle in $\Gamma$ with a pair of antipodal points $x_1$ and $x_2$ such that each edge of $\omega$ containing $x_i\in \{x_1,x_2\}$ is either not large or oriented towards $x_i$, then the cycle has a diagonal.
\een Then $A_\Gamma\times \Z$ is a Garside group. 
\emthm

\begin{figure}
	\begin{center}
		\begin{tabular}{cc}
			\begin{tikzpicture}
				\def \p {0.08}
				\def \l {6}
				\def \w {6}
				\draw[fill] (0,2) circle (\p) node(a) {};
				\draw[fill] (0.3,0) circle (\p) node(b){};
				\draw[fill] (0,-2) circle (\p) node(c){};
				\draw[fill] (2,-0.8) circle (\p) node(d){};
				\draw[fill] (2,0.8) circle (\p) node(e){};
				\draw[fill] (-2,-0.8) circle (\p) node(d'){};
				\draw[fill] (-2,0.8) circle (\p) node(e'){};
				\draw[-{Latex[length=\l, width=\w]}] (a) edge node[right] {$3$} (b);
				\draw[-{Latex[length=\l, width=\w]}] (b) edge node[right] {$3$} (c);
				\draw[-{Latex[length=\l, width=\w]}] (c) edge node[below right] {$3$} (d);
				\draw[-{Latex[length=\l, width=\w]}] (d) edge node[right] {$3$} (e);
				\draw[-{Latex[length=\l, width=\w]}] (e) edge node[above right] {$3$} (a);
				\draw[-{Latex[length=\l, width=\w]}] (c) edge node[below left] {$3$} (d');
				\draw[-{Latex[length=\l, width=\w]}] (d') edge node[left] {$3$} (e');
				\draw[-{Latex[length=\l, width=\w]}] (e') edge node[above left] {$3$} (a);
				\draw (a) edge (c);
				\draw (b) edge (d);
				\draw (c) edge (e);
				\draw (d) edge (a);
				\draw (e) edge (b);
				\draw (b) edge (d');
				\draw (c) edge (e');
				\draw (d') edge (a);
				\draw (e') edge (b);
			\end{tikzpicture}
			&
			\begin{tikzpicture}
				\def \p {0.08}
				\def \l {6}
				\def \w {6}
				\def \r {1.5}
				\def \R {3.5}
				\draw[fill] (0:\r) circle (\p) node(a) {};
				\draw[fill] (60:\r) circle (\p) node(b) {};
				\draw[fill] (120:\r) circle (\p) node(c) {};
				\draw[fill] (180:\r) circle (\p) node(d) {};
				\draw[fill] (240:\r) circle (\p) node(e) {};
				\draw[fill] (300:\r) circle (\p) node(f) {};
				\draw[fill] (0:\R) circle (\p) node(a') {};
				\draw[fill] (60:\R) circle (\p) node(b') {};
				\draw[fill] (120:\R) circle (\p) node(c') {};
				\draw[fill] (180:\R) circle (\p) node(d') {};
				\draw[fill] (240:\R) circle (\p) node(e') {};
				\draw[fill] (300:\R) circle (\p) node(f') {};
				\draw[fill] (67:2.3) circle (\p) node(B) {};
				\draw[fill] (247:2.3) circle (\p) node(E) {};
				\draw[-{Latex[length=\l, width=\w]}] (a) edge node[below left] {$3$} (b);
				\draw[-{Latex[length=\l, width=\w]}] (c) edge node[below ] {$3$} (b);
				\draw[-{Latex[length=\l, width=\w]}] (c) edge node[below right] {$5$} (d);
				\draw[-{Latex[length=\l, width=\w]}] (e) edge node[above right] {$4$} (d);
				\draw[-{Latex[length=\l, width=\w]}] (e) edge node[above] {$3$} (f);
				\draw[-{Latex[length=\l, width=\w]}] (a) edge node[above left] {$4$} (f);
				\draw[-{Latex[length=\l, width=\w]}] (b') edge node[above right] {$3$} (a');
				\draw[-{Latex[length=\l, width=\w]}] (b') edge node[above] {$4$} (c');
				\draw[-{Latex[length=\l, width=\w]}] (d') edge node[above left] {$5$} (c');
				\draw[-{Latex[length=\l, width=\w]}] (d') edge node[below left] {$3$} (e');
				\draw[-{Latex[length=\l, width=\w]}] (f') edge node[below] {$3$} (e');
				\draw[-{Latex[length=\l, width=\w]}] (f') edge node[below right] {$3$} (a');
				\draw[-{Latex[length=\l, width=\w]}] (a') edge node[below] {$3$} (a);
				\draw[-{Latex[length=\l, width=\w]}] (c') edge node[below left] {$3$} (c);
				\draw[-{Latex[length=\l, width=\w]}] (d) edge node[above] {$3$} (d');
				\draw[-{Latex[length=\l, width=\w]}] (f) edge node[above right] {$3$} (f');
				\draw[-{Latex[length=\l, width=\w]}] (b) edge node {$3$} (B);
				\draw[-{Latex[length=\l, width=\w]}] (B) edge node {$3$} (b');
				\draw[-{Latex[length=\l, width=\w]}] (e') edge node {$3$} (E);
				\draw[-{Latex[length=\l, width=\w]}] (E) edge node {$3$} (e);
				
				\draw (a) edge (B);
				\draw (a) edge (b');
				\draw (a') edge (b);
				\draw (a') edge (B);
				\draw (b) edge (b');
				\draw (b) edge (c');
				\draw (c) edge (b');
				\draw (B) edge (c);
				\draw (B) edge (c');
				\draw (d') edge (c);
				\draw (d) edge (c');
				\draw (d') edge (e);
				\draw (d') edge (E);
				\draw (d) edge (e');
				\draw (d) edge (E);
				\draw (e) edge (e');
				\draw (f') edge (e);
				\draw (f) edge (e');
				\draw (E) edge (f);
				\draw (E) edge (f');
				\draw (a) edge (f');
				\draw (a') edge (f);
			\end{tikzpicture}
		\end{tabular}
		\caption{Examples of Coxeter presentation graphs of Artin groups to which Theorem~\ref{thm:intro garside} applies.}
		\label{fig:example_graphs}
	\end{center}
\end{figure}
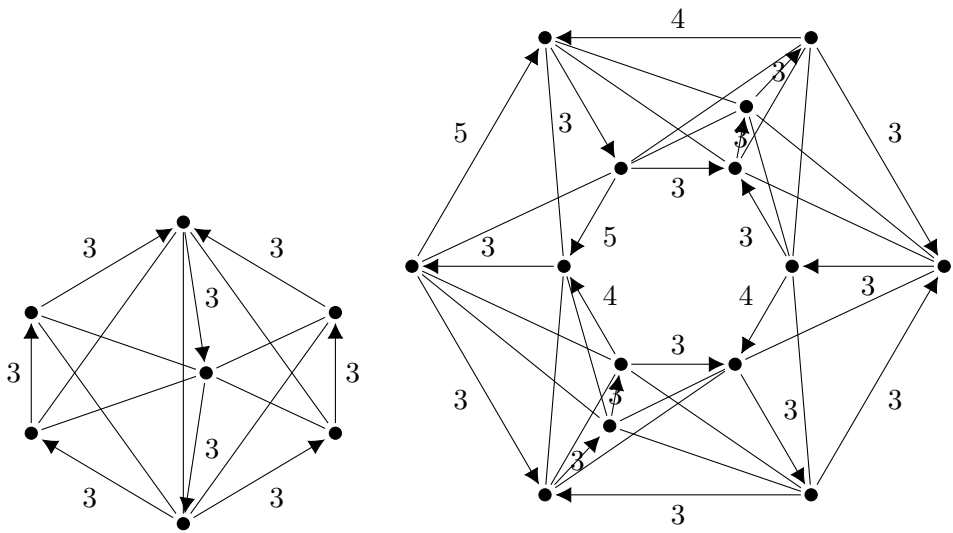

Below we include two simple examples of Coxeter presentation graph $\Gamma$ where the Theorem~\ref{thm:intro garside} applies, see Figure~\ref{fig:example_graphs}. The first is an amalgamation of two Artin groups of type $\widehat A_4$ along a spherical parabolic subgroup of type $A_3$. The second example is a bit more complicated, made of a few cyclic type Artin groups glued together in a cyclic way. Note that the edges without labels are assumed to be labeled by $2$.

In particular, all consequences listed in Theorem~\ref{thm:consequences_product_Z_garside} hold for this class of Artin groups. All of these consequences are new for this class, including the solvability of word problem. 
As a more precise comparison to previous results, we view the class of Artin groups in the above theorem as a combination of basic building blocks made of cyclic type Artin groups and spherical Artin groups. Then
\begin{enumerate}
	\item All consequences listed in Theorem~\ref{thm:consequences_product_Z_garside} were known before for spherical Artin groups \cite{charney_biautomatic,huang_osajda_helly,haettel_huang_weakly_modular}, hence also known for the Artin group of type $\widetilde A_n$, as the direct product of this Artin group with $\mathbb Z$ has finite index in a spherical Artin groups \cite{kent2002geometric};
	\item All consequences of Theorem~\ref{thm:consequences_product_Z_garside} except the first one (acting geometrically on a Helly graph) were known before for cyclic Artin groups with at most three generators - as these groups act geometrically on $\CAT(0)$ complexes made of equilateral triangles \cite{brady2000three};
	\item All consequences of Theorem~\ref{thm:consequences_product_Z_garside} are new for the remaining cyclic type Artin groups.
	\item To the best of our knowledge, for each of the properties in the list of consequences of Theorem~\ref{thm:consequences_product_Z_garside}, there does not exist a combination theorem \footnote{By a combination theorem, we mean a theorem of form that if a graph of groups $G$ has each of its vertex groups and edge groups satisfying property $X$, then $G$ also satisfies property $X$, under possibly some extra assumptions on $G$.} which is powerful enough to cover the pattern of combination of cyclic type and spherical type Artin groups in Theorem~\ref{thm:intro garside}, thus all consequences are new for the class of Artin groups in Theorem~\ref{thm:intro garside}. For example, the most recent combination theorem for Farrell-Jones conjecture \cite{knopf2019acylindrical} requires an acylindrical action of the group on a tree, which is not satisfied in our situation.
	\item Artin groups in Theorem~\ref{thm:intro garside} are in general not of type FC, so consequences of Theorem~\ref{thm:consequences_product_Z_garside} for this class do not follow from \cite{huang_osajda_helly}.
\end{enumerate}

\mk

These conditions are the most general that we can deal with. In particular, we isolate simple families of Artin groups to which this result applies.

\bmcor
Assume that $A_\Gamma$ is one of the following Artin groups:
\bit
\item $A_\Gamma$ is of large type (i.e. no commuting relations in its standard presentation), and has rank at most $3$.
\item $A_\Gamma$ is right-angled (i.e. only commuting relations in its standard presentation), without any induced squares.
\eit
Then $A_\Gamma \times \Z$ is Garside.
\emcor

We emphasize that even for one of the simplest and most extensively studied classes of Artin groups, namely the class of right-angled Artin groups, not much is known about the connection to Garside groups. For the free group $\F_r$ or rank $r$, Bessis has defined a quasi-Garside structure on $\F_r$ (\cite{bessis_free}). In comparison, it is interesting that we are able to endow the direct product $\F_r \times \Z$ with an actual Garside structure, and not a mere quasi-Garside structure.

\mk

Assumptions of Theorem~\ref{thm:intro garside} have a close connection to an existing result for a class of 2-dimensional Artin groups by \cite{brady2000three}. More precisely, \cite{brady2000three} studied the class of \emph{large type} Artin groups, i.e. each edge in the Coxeter presentation graph has label $\ge 3$. A dihedral subgroup of $A_\Gamma$ is a subgroup generated by two vertices in an edge of $\Gamma$. Interestingly, if we restrict Theorem~\ref{thm:intro garside} within the class of large type Artin groups,  the left forbidden configuration in \cite[Figure 5]{brady2000three} corresponds exactly to Assumption 1 in Theorem~\ref{thm:intro garside}, and the right forbidden configuration in \cite[Figure 5]{brady2000three} corresponds exactly to Assumption 2 in Theorem~\ref{thm:intro garside}. There is an interesting geometric phenomenon behind this.

The strategy in \cite{brady2000three} is to consider a dual Garside structure of each dihedral subgroup (choosing a dual Garside structure amounts to choosing an orientation of the associated edge), metrize each triangle in the presentation complex with respect to the dual Garside structure as flat equilateral triangles, and gluing these presentation complexes for dihedral subgroup in a natural way to obtain a complex with fundamental group $A_\Gamma$. Then the result \cite[Thereom 7]{brady2000three} implies that as long as the presentation graph $\Gamma$ avoids two configurations in \cite[Figure 5]{brady2000three}, then the resulting complex is locally $\CAT(0)$. 

Theorem~\ref{thm:intro garside} has a geometric counterpart (cf. Corollary~\ref{cor:complex1}). More precisely, given an Artin group $A_\Gamma$, we can choose a dual Garside structure on each standard spherical parabolic subgroup in a consistent way (again such information can be encoded as an appropriate orientation of all large edges of $\Gamma$). The dual Garside structure on each spherical parabolic subgroup $H$ gives an associated Garside complex (Definition~\ref{def:interval_group}) with fundamental group $H$. By gluing these Garside complexes in a natural way (as in \cite[Section 5]{paolini_salvetti}), we obtain a complex $X_\Gamma$ with fundamental group $A_\Gamma$. Here we metrize each simplex in $X_\Gamma$ by a polyhedral norm which is related to the $\widetilde A_n$-geometry, see \cite{haettel_simplicial_npc} (the norm here is not Euclidean), which echoes the work of \cite{brady2000three} where they metrize triangles with Euclidean $\widetilde A_2$ shape. The assumptions in Theorem~\ref{thm:intro garside} will ensure that the universal cover of $X_\Gamma$ with such metric is a space with convex geodesic bicombing (see Definition~\ref{def:bicombing}), which can be viewed as a form of non-positive curvature, and echoes the CAT$(0)$ metric in \cite{brady2000three}. 

It is natural to ask, if we metrize each simplex in $X_\Gamma$ by Euclidean simplices with $\widetilde A_n$-shape, whether the complex we obtain is locally $\CAT(0)$. However, it is notoriously difficult to verify local $\CAT(0)$-ness in high dimension. Here this issue is bypassed by metrizing the simplices with different kinds of norms rather than the Euclidean norm. While the resulting metric is not locally $\CAT(0)$, it is almost as good as $\CAT(0)$ in the sense that it implies most of the consequences of $\CAT(0)$ groups. We refer to \cite{haettel_simplicial_npc,haettel_helly_kpi1}, as well as \cite{descombes_lang_hyperbolicity,descombes_lang_flats} for more discussion in this direction.

\mk

Interestingly, for every Artin group as in Theorem~\ref{thm:main1}, we have an answer to all four questions stated above for general Artin groups. In particular, we can deduce new cases of $K(\pi,1)$-conjecture from Theorem~\ref{thm:main1}.

\bmcor(=Corollary~\ref{cor:K(pi,1)1})
\label{cor:K(pi,1)}
Assume that $A_\Gamma$ is of hyperbolic cyclic type. Then $A_\Gamma$ satisfies the $K(\pi,1)$ conjecture and has a trivial center.
\emcor

More precisely, the $K(\pi,1)$-conjecture is new for 6 examples of Artin groups whose Coxeter groups act cocompactly on $\mathbb H^3$ or $\mathbb H^4$. These examples seem to be challenging from the viewpoint of other approaches to $K(\pi,1)$-conjecture. Though the $K(\pi,1)$-conjecture when $A_\Gamma$ is 2-dimensional hyperbolic cyclic type follows from previous work \cite{charney_davis_kpi1}, and there is also a more recent proof in \cite{delucchi_paolini_salvetti_rank_3} using dual quasi-Garside structures.

To put Corollary~\ref{cor:K(pi,1)} in another context, note that the $K(\pi,1)$-conjecture is proved by Artin groups associated with reflection groups acting on $\mathbb S^n$ by Deligne \cite{deligne}, and Artin groups associated with reflection groups acting on $\E^n$ by Paolini and Salvetti \cite{paolini_salvetti}. The next step is to look at Artin groups associated with reflection groups acting on $\H^n$ (we call them hyperbolic type Artin groups), whose $K(\pi,1)$-conjecture is wide open. A fundamental subclass of hyperbolic type Artin groups is those associated with hyperbolic reflection groups whose fundamental domain is a compact simplex. This subclass is classified by Lanner \cite{lanner1950complexes}, which consists of infinitely many members in dimension 2 (whose $K(\pi,1)$-conjecture is already understood \cite{charney_davis_kpi1,delucchi_paolini_salvetti_rank_3}), and in 14 remaining cases in higher dimension. From this perspective, Corollary~\ref{cor:K(pi,1)} treats 6 out of these 14 remaining cases.

Corollary~\ref{cor:K(pi,1)} also follows from
another article of the second named author \cite[Theorem 1.4]{huangfourcycle}, via an alternative approach to the $K(\pi,1)$-conjecture. However, the method here establishes all the properties in Theorem~\ref{thm:consequences_product_Z_garside} for hyperbolic cyclic type Artin groups, which are not consequences of \cite{huangfourcycle}.

\subsection{Structure of the article}
In Section~\ref{sec:background}, we collect some background, notably on Garside groups, Artin groups and nonpositively curved spaces. In Section~\ref{sec:general_construction}, we discuss the general criterion of making $G\times \Z$ a Garside group and prove Theorem~\ref{thm:U}. Then we discuss examples of $T(5)$ and systolic restricted presentation groups. In Section~\ref{sec:Artin}, we adapt Theorem~\ref{thm:U} to the special situation of Artin groups, and produce a criterion for when an Artin group times $\Z$ is Garside, see Proposition~\ref{pro:E_lattice} and Corollary~\ref{cor:Artin_nice_U_garside}. In Section~\ref{sec:cyclic_type}, we verify the criterion in Proposition~\ref{pro:E_lattice} and Corollary~\ref{cor:Artin_nice_U_garside} for cyclic type Artin groups. In Section~\ref{sec:combination} we treat more general Artin groups and prove Theorem~\ref{thm:intro garside}.

\subsection{Acknowledgement}

The authors would like to thank Anthony Genevois, Jon McCammond, Alex Martin, Damian Osajda and Dani Wise for interesting discussions. The authors thank the Centre de recherches math\'ematiques de Montr\'eal for hospitality. The authors would like to thank the anonymous referees for their careful comments which greatly helped to improve this article.

Thomas Haettel was partially supported by French project ANR-22-CE40-0004 GOFR. Jingyin Huang was partially supported by a Sloan fellowship. 
\section{Background}
\label{sec:background}

We start by giving references for Theorem~\ref{thm:consequences_product_Z_garside} in the introduction, then we collect background definitions and results concerning Artin groups, Garside groups, and nonpositively curved spaces.

\subsection{Proof of Theorem~\ref{thm:consequences_product_Z_garside}} \label{subsec:proof_product_Z_garside}

We now give precise references for the various items of Theorem~\ref{thm:consequences_product_Z_garside} from the introduction, listing consequences for a group $G$ such that $G \times \Z$ is Garside.

\bp
\ben
\item This is a consequence of~\cite{huang_osajda_helly}, see also~\cite{haettel_helly_kpi1}.
\item This is a consequence of~\cite[Proposition~3.25]{garside}.
\item This is a consequence of~\cite[Corollary~9.8]{haettel_simplicial_npc}.
\item This is a consequence of~\cite{haettel_huang_weakly_modular}.
\item This is a consequence of~\cite{mosher_biautomatic}. For consequences of biautomaticity, see for instance~\cite{bridson_haefliger,wenger2003isoperimetric,behrstock2019combinatorial}.
\item This is a consequence of~\cite[Proposition~7.10]{helly_groups}.
\item This is a consequence of~\cite{lee_lee_garside_translation}, see also~\cite{haettel_osajda_locally_elliptic}.
\item Since $G$ acts geometrically on a metric space with a convex geodesic bicombing, according to~\cite[Theorem~6.1]{kasprowski_rueping}, it satisfies the Farrell-Jones conjecture with finite wreath products.
\item Since $G$ acts geometrically on a metric space with a convex geodesic bicombing, according to~\cite{fukaya_oguni}, it satisfies the coarse Baum-Connes conjecture.
\item This is a consequence of the Farrell-Jones conjecture and~\cite[Theorem~0.12]{bartels_luck_holger_farrell_jones}.
\een
\ep

\subsection{Coxeter groups and Artin groups}
\label{subsec:Coxeter Artin}

We recall the definitions of Coxeter groups and Artin groups.

\mk

For every finite simple graph $\Gamma$ with vertex set $S$ and with edges labeled by some integer in $\{2,3,\dots\}$, one associates the Coxeter group $W_\Gamma$ with the following presentation:
$$W_\Gamma = \<S \st \forall \{s,t\} \in \Gamma^{(1)}, \forall s \in S, s^2=1, [s,t]_m=[t,s]_m \mbox{ if the edge $\{s,t\}$ is labeled $m$}\>,$$
where $[s,t]_m$ denotes the word $ststs\dots$ of length $m$. Such a graph $\Gamma$ may be called a \emph{Coxeter presentation graph}, emphasizing the fact that edges correspond to relations. We will also use $W_S$ to denote $W_\Gamma$.

\mk

We will also be using a graph closely related to $\Gamma$, the \emph{Dynkin diagram} $\Gamma_D$: it has the same vertex set $S$, with some edges labeled in $\{4,5,\dots,\infty\}$, with the following edges between vertices $s,t \in S$:
\bit
\item If there is an edge labeled $2$ between $s$ and $t$ in $\Gamma$, there is no edge between $s$ and $t$ in $\Gamma_D$.
\item If there is an edge labeled $3$ between $s$ and $t$ in $\Gamma$, there is an unlabeled edge between $s$ and $t$ in $\Gamma_D$.
\item If there is an edge labeled by $m \geq 4$ between $s$ and $t$ in $\Gamma$, there is the same edge between $s$ and $t$ in $\Gamma_D$ labeled $m$.
\item If there is no edge between $s$ and $t$ in $\Gamma$, there is an edge between $s$ and $t$ in $\Gamma_D$ labeled $\infty$.
\eit

\mk

The associated Artin group $A_\Gamma$ is defined by a similar presentation:
$$A_\Gamma = \<S \st \forall \{s,t\} \in \Gamma^{(1)}, [s,t]_m=[t,s]_m \mbox{ if the edge $\{s,t\}$ is labeled $m$}\>.$$
The groups $A_\Gamma$ are also called Artin-Tits groups since they have been defined by Tits in~\cite{tits_artin}. We will also use $A_S$ to denote $A_\Gamma$.

Note that only the relations $s^2=1$ have been removed, so that there is a natural surjective morphism from $A_\Gamma$ to $W_\Gamma$. Also note that when $m=2$, then $s$ and $t$ commute, and when $m=3$, then $s$ and $t$ satisfy the classical braid relation $sts=tst$.

For a subset $S'$ of the generating $S$, the subgroup of $A_\Gamma$ or $W_\Gamma$ generated by $S'$ is called a \emph{standard parabolic subgroup}. A standard parabolic subgroup of an Artin group is itself an Artin group \cite{vanderlek}. A similar statement is true for Coxeter groups \cite{bourbaki_lie_456}. A \emph{parabolic subgroup} is a conjugate of a standard parabolic subgroup.

\mk

Most results about Artin-Tits groups concern particular classes. The Artin group $A_\Gamma$ is called:
\bit
\item of \emph{spherical type} if its associated Coxeter group $W_\Gamma$ is finite, i.e. may be realized as a reflection group of a sphere.
\item of \emph{Euclidean type} if its associated Coxeter group $W_\Gamma$ may be realized as a reflection group of an Euclidean space.
\item of \emph{hyperbolic type} if its associated Coxeter group $W_\Gamma$ may be realized as a reflection group of a real hyperbolic space.
\eit

\mk

We say a Coxeter group $W_S$ is of \emph{cyclic type} if the associated Dynkin diagram is a cycle, and the parabolic subgroup generated by $S\setminus\{s\}$ is spherical for any vertex $s\in\Gamma$. We list in Table~\ref{tab:cyclic} the Dynkin diagrams of cyclic type. Note that we use in this table the convention of Dynkin diagrams: vertices that are not joined by an edge commute, and we drop the label $3$ from edges. Note that cyclic type Coxeter groups are either of Euclidean type or of hyperbolic type.

\begin{table}
	\begin{center}
		\begin{tabular}{|c|c|c|c|c|}
			\hline
			Name 
			&
			$\tilde{A_n}$, for $n \geq 3$
			&
			Triangle
			&
			$3-3-3-4$
			&
			$3-3-3-5$
			\\
			\hline
			Dynkin diagram
			&
			
			\begin{tikzpicture}
				\def \p {0.05}
				\def \r {1}
				\def \op {1}
				\def \gris {black!10}
				
				\draw[fill] (0:\r) circle (\p) node(s1) {};
				\draw[fill] (60:\r) circle (\p) node(s2) {};
				\draw[fill] (120:\r) circle (\p) node(s3) {};
				\draw[fill] (180:\r) circle (\p) node(s4) {};
				\draw[fill] (240:\r) circle (\p) node(s5) {};
				\draw[fill] (300:\r) circle (\p) node(s6) {};
				
				\node (dots) at (0,-0.9) {\bfseries $\dots$};
				
				\draw [-] (s1) edge (s2) (s2) edge (s3) (s3) edge (s4) (s4) edge (s5) (s6) edge (s1);
				
			\end{tikzpicture}
			&
			
			\begin{tikzpicture}
				\def \p {0.05}
				\def \r {1}
				\def \op {1}
				\def \gris {black!10}
				
				\draw[fill] (0:\r) circle (\p) node(s1) {};
				\draw[fill] (120:\r) circle (\p) node(s2) {};
				\draw[fill] (240:\r) circle (\p) node(s3) {};
				
				\draw [-] (s1) edge (s2) (s2) edge (s3) (s3) edge (s1);
				
				\node (label) at (-0.9,0) {$\geq 3$};
				\node (label) at (0.6,0.6) {$\geq 3$};
				\node (label) at (0.6,-0.6) {$\geq 3$};
				
			\end{tikzpicture}
			&
			
			\begin{tikzpicture}
				\def \p {0.05}
				\def \r {1}
				\def \op {1}
				\def \gris {black!10}
				
				\draw[fill] (45:\r) circle (\p) node(s1) {};
				\draw[fill] (135:\r) circle (\p) node(s2) {};
				\draw[fill] (-135:\r) circle (\p) node(s3) {};
				\draw[fill] (-45:\r) circle (\p) node(s4) {};
				
				\draw [-] (s1) edge (s2) (s2) edge (s3) (s3) edge (s4) (s4) edge (s1);
				
				\node (label) at (-1,0) {$4$};
			\end{tikzpicture}
			&
			
			\begin{tikzpicture}
				\def \p {0.05}
				\def \r {1}
				\def \op {1}
				\def \gris {black!10}
				
				\draw[fill] (45:\r) circle (\p) node(s1) {};
				\draw[fill] (135:\r) circle (\p) node(s2) {};
				\draw[fill] (-135:\r) circle (\p) node(s3) {};
				\draw[fill] (-45:\r) circle (\p) node(s4) {};
				
				\draw [-] (s1) edge (s2) (s2) edge (s3) (s3) edge (s4) (s4) edge (s1);
				
				\node (label) at (-1,0) {$5$};
				
			\end{tikzpicture}\\
			
			\hline
			
			Name 
			&
			$3-4-3-4$
			&
			$3-4-3-5$
			&
			$3-5-3-5$
			&
			$3-3-3-3-4$
			\\
			\hline
			
			Dynkin diagram
			&

			\begin{tikzpicture}
				\def \p {0.05}
				\def \r {1}
				\def \op {1}
				\def \gris {black!10}
				
				\draw[fill] (45:\r) circle (\p) node(s1) {};
				\draw[fill] (135:\r) circle (\p) node(s2) {};
				\draw[fill] (-135:\r) circle (\p) node(s3) {};
				\draw[fill] (-45:\r) circle (\p) node(s4) {};
				
				\draw [-] (s1) edge (s2) (s2) edge (s3) (s3) edge (s4) (s4) edge (s1);
				
				\node (label) at (-1,0) {$4$};
				\node (label) at (1,0) {$4$};
			\end{tikzpicture}
			
			&
			
			\begin{tikzpicture}
				\def \p {0.05}
				\def \r {1}
				\def \op {1}
				\def \gris {black!10}
				
				\draw[fill] (45:\r) circle (\p) node(s1) {};
				\draw[fill] (135:\r) circle (\p) node(s2) {};
				\draw[fill] (-135:\r) circle (\p) node(s3) {};
				\draw[fill] (-45:\r) circle (\p) node(s4) {};
				
				\draw [-] (s1) edge (s2) (s2) edge (s3) (s3) edge (s4) (s4) edge (s1);
				
				\node (label) at (-1,0) {$4$};
				\node (label) at (1,0) {$5$};
			\end{tikzpicture}
			
			&
			
			\begin{tikzpicture}
				\def \p {0.05}
				\def \r {1}
				\def \op {1}
				\def \gris {black!10}
				
				\draw[fill] (45:\r) circle (\p) node(s1) {};
				\draw[fill] (135:\r) circle (\p) node(s2) {};
				\draw[fill] (-135:\r) circle (\p) node(s3) {};
				\draw[fill] (-45:\r) circle (\p) node(s4) {};
				
				\draw [-] (s1) edge (s2) (s2) edge (s3) (s3) edge (s4) (s4) edge (s1);
				
				\node (label) at (-1,0) {$5$};
				\node (label) at (1,0) {$5$};
			\end{tikzpicture}
			
			&
			
			\begin{tikzpicture}
				\def \p {0.05}
				\def \r {1}
				\def \op {1}
				\def \gris {black!10}
				
				\draw[fill] (0:\r) circle (\p) node(s1) {};
				\draw[fill] (72:\r) circle (\p) node(s2) {};
				\draw[fill] (144:\r) circle (\p) node(s3) {};
				\draw[fill] (-144:\r) circle (\p) node(s4) {};
				\draw[fill] (-72:\r) circle (\p) node(s5) {};
				
				\draw [-] (s1) edge (s2) (s2) edge (s3) (s3) edge (s4) (s4) edge (s5) (s5) edge (s1);
				\node (label) at (-1,0) {$4$};
				
			\end{tikzpicture}
			\\
			
			\hline

		\end{tabular}
		\caption{Dynkin diagrams of cyclic type}
		\label{tab:cyclic}
	\end{center}
\end{table}

For an element $g$ in Coxeter group $W_S$, we can represent $g$ as a word in the free monoid on $S$. Such representation is
\emph{reduced} if its length is the shortest possible among words in the free monoid that represent $g$. It is known that any two reduced words representing the same element in $W_S$ differ by a finite sequence of moves applying the relation in $W_S$, see e.g. \cite{humphreys1990reflection}. Thus each element in $W_S$ has a well-defined \emph{support}, which is the collection of elements in $S$ which appears in a reduced word representing this element.

A subset $S'\subset S$ is \emph{irreducible} if it spans a connected subgraph of the Dynkin diagram, otherwise $S'$ is \emph{reducible}.

A conjugate of an element of the standard generating set $S$ is called a \emph{reflection of $W_S$}.
\begin{lem}
	\label{lem:irreducible}
	The support of each reflection is irreducible.
\end{lem} 

\bp
Let $r=wsw^{-1}$ be a reflection in $W_S$ with $s\in S$ and $w\in W_S$. If $\supp(r)$ is reducible, then $\supp(r)=I_1\sqcup I_2$ with elements in $I_1$ commuting with elements in $I_2$. As $\langle r\rangle$ is a parabolic subgroup of $W_S$ which is contained in the standard parabolic subgroup $W_{I_1\cup I_2}$, by the remark after \cite[Lemma 5.3.6]{davis_coxeter}, there exists $w'\in W_{I_1\cup I_2}$ and $s'\in I_1\cup I_2$ such that $r=w's'(w')^{-1}$. We assume without loss of generality that $s'\in I_1$. Write $w'=w'_1w'_2$ with $w'_i\in W_{I_i}$ for $i=1,2$. Then $r=w'_1s'(w'_1)^{-1}$ and $\supp(r)\subset I_1$, which is a contradiction. Thus the lemma is proved.
\ep

\subsection{The $K(\pi,1)$-conjecture} \label{subsec:Kpi1}

Artin groups are closely related to hyperplane complements, which can be presented in a simple way in spherical, Euclidean, and hyperbolic types. Fix a Coxeter group $W=W_\Gamma$ of spherical type, Euclidean or hyperbolic type acting by isometries on a sphere $\SS^{n-1}$, Euclidean space $\R^{n-1}$ or a real hyperbolic space $\H^{n-1}$, where the standard generators act by reflections.

In the case of $\SS^{n-1}$, we will consider $W$ as a subgroup of $O(n)$ acting by linear transformations on $\Omega=\R^n$.
In the case of $\R^{n-1}$, we will consider $W$ as a subgroup of $\GL(n)$ acting by linear transformations on $\R^n$, preserving the hyperplane $\{x_n=1\}$ and acting by isometries on it. The group $W$ preserves the open cone $\Omega=\{x_n>0\}$ of $\R^n$.
In the case of $\H^{n-1}$, we will consider $W$ as a subgroup of $O(n-1,1)$ acting linearly on $\R^n$, and preserving the open cone $\Omega=\H^{n-1}$. Let ${\mathcal R}$ denote the set of reflections of $W$. Consider the family of linear hyperplanes of $\R^n$
$${\mathcal H} = \{H_r \st r \in {\mathcal R}\},$$
where $H_r \subset \R^n$ denotes the fixed point set of the reflection $r$.

\mk

The analogue of the complement of the complexified hyperplane arrangement is
$$M(\Gamma) = (\Omega \times \Omega) \mathbin{\big\backslash} \bigcup_{r \in {R}} (H_r \times H_r),$$
see~\cite{paris_kpi1} for more details. Note that $W$ acts naturally on $M$, and we have the following (see~\cite{vanderlek}):
$$\pi_1(W_\Gamma \bs M(\Gamma)) \simeq A_\Gamma.$$
So the Artin group $A_\Gamma$ appears as the fundamental group of (a quotient of) the complement of a complexified hyperplane arrangement. One very natural question is to decide whether it is a classifying space. This is the statement of the following conjecture.

\bconj[$K(\pi,1)$ conjecture]
The space $M(\Gamma)$ is aspherical.
\econj

This conjecture has been proved for spherical type Artin groups by Deligne in~\cite{deligne}, for 2-dimensional and type FC Artin groups by Charney and Davis in \cite{charney_davis_kpi1}, and for Euclidean type Artin groups by Paolini and Salvetti in~\cite{paolini_salvetti} very recently.

\mk

\subsection{Interval groups and Garside groups} \label{subsec:garside}

We will follow McCammond's article~\cite{mccammond_intro_garside} for the description of interval groups.

\bdf[Posets]
A poset $P$ is called \emph{bounded} if it has a minimum, denoted $0$, and a maximum, denoted $1$. 

For $x\le y$ in a poset $P$, the \emph{interval} between $x$
and $y$ is the restriction of the poset to those elements $z$ with $x\le z\le y$. We denote this interval by $[x,y]$.
A poset $P$ is called \emph{graded} if for any $x\le y$ in $P$, any chain in $[x,y]$ belongs to a maximal chain and all maximal chains have the same finite length.

A poset $P$ is called \emph{weakly graded} if there is a poset map $r:P \ra \Z$, i.e. such that for every $x<y$ in $P$, we have $r(x)< r(y)$: the map $r$ is called a \emph{rank map}. A poset $P$ is called \emph{weakly boundedly graded} if there is a rank map $r:P \ra \Z$ with finite image. Note that a finite weakly graded poset is always weakly boundedly graded.

An \emph{upper bound} for a pair of elements $a,b\in P$ is an element $c\in P$ such that $a\le c,b\le c$. A \emph{minimal upper bound} for $a,b$ is an upper bound $c$ such that there does not exist upper bound $c'$ of $a,b$ such that $c'<c$.
The \emph{join} of two elements $a,b$ in $P$ is an upper bound $c$ of them such that for any other upper bound $c'$ of $a,b$, we have $c\le c'$. We define \emph{lower bound}, \emph{maximal lower bound}, and \emph{meet} similarly.
In general, the meet or join of two elements in $P$ might not exist.
A poset $P$ is a \emph{lattice} if any pair of elements have a meet and a join.

A poset $P$ is a \emph{meet-semilattice} (resp. join-semilattice) if any pair of elements have a meet (resp. a join).
\edf

\bdf
We say that a poset $P$ contains a \emph{bowtie} if there
exist pairwise distinct elements $a,b,c$ and $d$ such that $a,b < c,d$, and there exists no $x \in P$ such that $a,b \leq x \leq c,d$.
\edf

It turns out that bowties are the only obstruction to being a lattice, for a weakly graded poset. This is proved in~\cite[Proposition 1.5]{brady_mccammond} for bounded graded lattices. This also holds for weakly graded lattices, so we give a proof here for the convenience of the reader.

\begin{pro} \label{pro:bowtie}
	Let $L$ denote a weakly graded poset. Then $L \cup \{0,1\}$ is a lattice if and only if $L$ has no bowtie.
\end{pro}

\bp
Assume that $L \cup \{0,1\}$ is a lattice, and consider $a,b<c,d$ in $L$. Then the meet $x$ of $c,d$ is such that $a,b \leq x \leq c,d$. So $L$ has no bowties.

\mk

Conversely, assume that $L$ has no bowtie. Note that $L \cup \{0,1\}$ has no bowtie either. Fix $a,b \in L$, and let $M$ denote the set of upper bounds of $a$ and $b$ in $L \cup \{0,1\}$: we have $1 \in M$, so $M$ is not empty. Let us consider a sequence $(x_n)_{n \in \N}$ in $M$ such that for each $n \in \N$, we have $x_n \geq x_{n+1}$. Let $r:P \ra \Z$ denote a weak grading on $P$. Then the sequence $(r(x_n))_{n \in \N}$ in $\Z$ is non-increasing and bounded below by $r(a)$, so it is eventually constant. This implies that the sequence $(x_n)_{n \in \N}$ itself is eventually constant.

We may therefore consider a minimal element $x$ of $M$. We will prove that $x$ is  unique: by contradiction, assume that $y \in M$ is a minimal element distinct from $x$. Then $a,b<x,y$ form a bowtie. Hence $x$ is the unique minimal element of $M$, and it is the join of $a$ and $b$ in $L \cup \{0,1\}$.

Similarly, any two elements of $L$ has a meet in $L \cup \{0,1\}$. So $L \cup \{0,1\}$ is a lattice.
\ep

Here is one definition of Garside groups. We refer the reader to~\cite{garside} and \cite{mccammond_intro_garside} for more background on Garside groups. We also refer the reader to~\cite{haettel_huang_weakly_modular} to equivalent definitions of Garside groups, which are more geometric in flavour.

\bdf[Garside group]
\label{def:garside}
Let $G$ denote a group, $S \subset G$ a finite subset and $\Delta \in G$. The triple $(G,S,\Delta)$ is called a \emph{Garside structure} if the following conditions hold. Let $G^+$ denote the submonoid of $G$ generated by $S$.
\ben
\item The group $G$ is generated by $S$.
\item For any element $g \in G^+$, there is a bound on the length of expressions $g=s_1 \dots s_n$, where $s_1,\dots,s_n \in S \bs \{1\}$.
\item We define the partial orders $\leq_L$, $\leq_R$ on $G^+$ by $a\leq_L b$ if and only if $b=ac$ for some $c\in G^+$ and $a\leq_R b$ if and only if $b=ca$ for some $c\in G^+$.
The left $\leq_L$ and right $\leq_R$ orders on $G^+$ are lattices.
\item The element $\Delta\in G^+$, and the set $S$ is a balanced interval between $1$ and $\Delta$, i.e.
$$S = \{g \in G^+ \st 1 \leq_L g \leq_L \Delta\} =  \{g \in G^+ \st 1 \leq_R g \leq_R \Delta\}.$$
\een
A group is called \emph{Garside} if it admits such a Garside structure, and $\Delta$ is called the \emph{Garside element}.
If the set $S$ is allowed to be infinite, but with the additional assumption that $(S,\le_L)$ and $(S,\le_R)$ are weakly boundedly graded, we may say that $(G,S,\Delta)$ is a \emph{quasi-Garside structure}.
\edf

\bdf[Labeled posets]
If $P$ is a poset, the \emph{set of intervals} is $I(P)=\{(x,y) \in P^2 \st x \leq y\}$.

Let $P$ denote a bounded poset, and let $S$ denote a labeling set.

An \emph{interval-labeling} of $P$ is a map $\lambda:I(P) \ra S$.

An interval-labeling $\lambda$ is \emph{group-like} if, for any two chains $x \leq y \leq z$ and $x' \leq y' \leq z'$ having two pairs of corresponding labels in common, the third pair of labels are equal.

An interval-labeling $\lambda$ is \emph{balanced} if
$$\{\lambda(0,x) \st x \in P\} = \{\lambda(x,1) \st x \in P\} = \{\lambda(x,y) \st (x,y) \in I(P)\}.$$
\edf

Note that McCammond's definition of balanced interval labeling (\cite[Definition~1.11]{mccammond_intro_garside}) only requires the first equality to hold. However, McCammond states that the second inequality is a consequence of being balanced and group-like, which does not seem obvious. We therefore chose to strengthen the definition of a balanced labeling, to ensure that all consequences of a combinatorial Garside structure hold.

\bdf[Interval complex and interval group] \label{def:interval_group}
Let $P$ denote a poset with a group-like interval-labeling $\lambda$.

Let us consider the quotient $K_P$ of the geometric realization $|P|$ of $P$, where the $k$-simplices corresponding to two $k$-chains $(x_0<x_1< \dots <x_k)$ and $(x'_0<x'_1< \dots <x'_k)$ are identified if and only if $\lambda(x_0,x_1)=\lambda(x'_0,x'_1),\dots,\lambda(x_{k-1},x_k)=\lambda(x'_{k-1},x'_k)$. It is called the \emph{interval complex} of $P$.

The fundamental group $G_P$ of $K_P$ is called the \emph{interval group} of $P$, it is naturally a quotient of the free group over the labeling set $S$.
\edf

\bexe
Let us consider the Boolean lattice $P={\mathcal P}(S)$ consisting of all subsets of a finite set $S$. The geometric realization $|P|$ of $P$ is isomorphic to a simplicial subdivision of the cube $[0,1]^S$.

For each $x \subset y \subset S$, let us consider the labeling $\lambda(x,y)=y\setminus x \in P$. The corresponding quotient $K_P$ is isomorphic to a simplicial subdivision of the torus $(\SS^1)^S$. The interval group $G_P$ is isomorphic to the free abelian group $\Z^S$, with the following presentation:
$$G_P = \<P \st \forall x \subset y \subset z \subset S, (y\setminus x) \cdot (z\setminus y) = (z\setminus x)\> \simeq \Z^S.$$ 
\eexe

\bdf[Combinatorial Garside structure]
\label{def:comb Gar}
A \emph{combinatorial Garside structure} is a poset $P$ with an interval-labeling $\lambda:I(P) \ra S$ such that:
\bit
\item $P$ is a bounded, weakly graded lattice.
\item $\lambda$ is group-like and balanced.
\eit
We say the combinatorial Garside structure is \emph{finite} if $P$ consists of finitely many elements.
\edf

Combinatorial Garside structures are just an explicit combinatorial way to describe arbitrary Garside groups, as explained by McCammond.

\bthm \cite[Theorem~1.17]{mccammond_intro_garside}
\label{thm:gar}
A group $G$ is a quasi-Garside group if and only if $G$ is isomorphic to the interval group of a combinatorial Garside structure. A group $G$ is a Garside group if and only if $G$ is isomorphic to the interval group of a finite combinatorial Garside structure.
\ethm

\begin{proof}[Sketch of proof]
	As the proof of this theorem is not found in \cite{mccammond_intro_garside}, we provide a sketch of how this can be deduced.
	
	\mk
	
	If $G$ is a quasi-Garside group, then we can consider the interval $P=([1,\Delta],\le_L)$. If $a\le_L b$ in $P$, then the interval $[a,b]$ is labeled by $a^{-1}b$. This makes $P$ a combinatorial Garside structure, and it is clear that the relations in the presentation of the interval group of $P$ hold in $G$. Conversely, let us show that these relations are sufficient to define $G$. Let us consider $s_1,\dots,s_n$ in $P \cup P^{-1}$ such that $s_1 s_2 \dots s_n=1$. Using $\Delta$ and the fact that $P$ is balanced, we may assume that $s_1,\dots,s_n$ are in $P$ and satisfy $s_1 s_2 \dots s_n=\Delta^k$, for some $k \geq 0$. Without loss of generality, we may assume that $s_1,s_n \neq 1,\Delta$. Therefore the word $s_1 s_2 \dots s_n$ is not the left Garside normal form of this element, which is $\Delta^k$. Hence there exists $1 \leq i \leq n-1$ such that $\Delta \wedge (s_is_{i+1}) = s'_i \neq s_i$. Let us write $s'_i=s_iu$, with $u \in P$, and $s_{i+1}=us'_{i+1}$, with $s'_{i+1} \in P$. Both relations hold in the interval group of $P$, hence we can replace $s_is_{i+1}$ by $s'_is'_{i+1}$ in the word $s_1s_2 \dots s_n$. By iterating this process, we arrive at the left Garside normal form $\Delta^k$. We conclude that the equality $s_1s_2 \dots s_n = \Delta^k$ holds in the interval group of $P$. Hence $G$ is isomorphic to the inteval group of $P$.
	
	\mk 
	
	Conversely, given $P$ being a combinatorial Garside structure, then the collection of interval labels of $P$ form an associative Noetherian germ in the sense of \cite[Section 2.1]{digne_garside_Cn}. This germ has a special element, $\Delta$, which is the interval label associated to the interval between the minimum and maximum of $P$. This germ defines a monoid $M(P)$, as explained in \cite[Section 2.1]{digne_garside_Cn}, and this germ satisfies the assumptions of \cite[Proposition 2.4]{digne_garside_Cn}. Hence \cite[Proposition 2.4]{digne_garside_Cn} implies that the group of fraction of $M(P)$ is a Garside group. This Garside group is isomorphic to the interval group of this combinatorial Garside structure.
\end{proof}
%\brk
%More generally, a group is quasi-Garside if and only if it is isomorphic to the interval group of an arbitrary combinatorial Garside structure.
%\erk

\subsection{Nonpositive curvature: Helly graphs and CUB spaces} \label{subsec:nonpositive_curvature}

We will present briefly several notions of metric spaces and graphs of nonpositive curvature which are relevant to Garside groups.

Let us start with Helly graphs: we refer the reader to~\cite{helly_groups} for more details. Let $\Gamma$ be a simplicial graph. We endow $\Gamma$ with a path metric $d$ such that each edge has length 1. A \emph{combinatorial ball} in $\Gamma$ is the collection of all vertices in some metric ball $B(x,r)$ centered at a vertex of $\Gamma$. Recall that the action of a group $G$ on a metric space $X$ is \emph{geometric}, if  the action is isometric, proper, and cocompact (i.e. $X/G$ is compact).

\bdf[Helly graph, Helly group]
Let $\Gamma$ be a simplicial graph, which is connected. Then $\Gamma$ is called \emph{Helly} if any family of pairwise intersecting combinatorial balls have a non-empty total intersection.

A group is called \emph{Helly} if it acts geometrically by automorphisms on a Helly graph. 
\edf

Helly graphs are combinatorial analogues of injective metric spaces, defined as follows.

\bdf[Injective metric spaces]
\label{def:injective}
A geodesic metric space $X$ is \emph{injective} if any family of pairwise intersecting closed metric balls in $X$ have a non-empty total intersection.
\edf

Helly groups enjoy many properties which are typical of nonpositive curvature, see for instance~\cite{helly_groups}, \cite{lang} and \cite{haettel_osajda_locally_elliptic} and also Theorem~\ref{thm:consequences_product_Z_garside}.

\mk

A much weaker, but way broader notion is that of weakly modular graphs, see~\cite{chalopin_chepoi_hirai_osajda}. These graphs encompass many "nonpositive curvature type" graphs, such as Helly graphs, (weakly) systolic graphs, median and quasi-median graphs, modular graphs.

\bdf[Weakly modular graph]
A connected graph $\Gamma$ is called \emph{weakly modular} if it satisfies the triangle condition (TC) and the quadrangle condition (QC):
\bit
\item[(TC)] For any $x,y,z \in \Gamma^{(0)}$ such that $d(y,z)=1$ and $d(x,y)=d(x,z)=n \geq 2$, there exists $t \in \Gamma^{(0)}$ such that $d(t,y)=d(t,z)=1$ and $d(x,t)=n-1$.
\item[(QC)] For any $x,y,z,u \in \Gamma^{(0)}$ such that $d(y,u)=d(z,u)=1$, $d(y,z)=2$, $d(x,u)=n \geq 3$ and $d(x,y)=d(x,z)=n-1$, there exists $t \in \Gamma^{(0)}$ such that $d(t,y)=d(t,z)=1$ and $d(x,t)=n-2$. 
\eit
\edf

\mk

Many of the consequences for Helly groups rely simply on the existence of a convex geodesic bicombing, whose definition we recall here. We also recall the definition of CUB spaces and groups, defined in~\cite{haettel_simplicial_npc}.

\bdf[Bicombing, CUB]
\label{def:bicombing}
A \emph{convex geodesic bicombing} on a metric space $X$ is a continuous map $\sigma:X \times X \times [0,1] \ra X$ such that:
\bit
\item For each $x,y \in X$, the map $t \in [0,1] \mapsto \sigma(x,y,t)$ is a constant speed reparametrized geodesic from $x$ to $y$.
\item For each $x,x',y,y' \in X$, the map $t \in [0,1] \mapsto d(\sigma(x,y,t),\sigma(x',y',t))$ is convex.
\eit
A metric space is called \emph{CUB}, for Convexly Uniquely Bicombable, if it admits a unique convex geodesic bicombing.
A group is called \emph{CUB} if it acts geometrically by isometries on a CUB space.
\edf

Groups acting on spaces with convex bicombings enjoy many properties, see for instance~\cite{descombes_lang_flats} and \cite{descombes_lang_hyperbolicity}. Furthermore, CUB groups satisfy some extra properties presented in~\cite{haettel_simplicial_npc}, see also Theorem~\ref{thm:consequences_product_Z_garside}.

\mk

One major incarnation of the nonpositive curvature properties of Garside groups is the following.

\bthm[\cite{huang_osajda_helly}, see also~\cite{haettel_helly_kpi1}]
Any Garside group acts geometrically by automorphisms on a Helly graph.
\ethm

The quotient of a Garside group by the cyclic subgroup generated by the Garside element also has nonpositive curvature in the following sense.

\bthm[\cite{haettel_simplicial_npc},\cite{haettel_huang_weakly_modular}]
Let $G$ denote a Garside group, with Garside element $\Delta$. Then the group $G/\<\Delta\>$ acts geometrically by isometries on a CUB space, and it acts geometrically by automorphisms on a weakly modular graph.
\ethm

\subsection{Dual Garside structures on spherical type Artin groups} \label{subsec:dual_artin_garside}

Dual Garside structure on spherical type Artin groups have been studied notably by Birman-Ko-Lee (\cite{birman_ko_lee}) and Bessis (\cite{bessis}), see also~\cite{paolini_dual_approach} for an overview of dual Garside structures on general Artin groups. We also refer the reader to~\cite{mccammond_intro_garside} for the point of view of interval groups that we are presenting here.

\mk

Let $\Gamma$ denote a Coxeter presentation graph, with vertex set $S$. Given any linear ordering $S=\{s_1,\dots,s_n\}$ of $S$, we have an associated \emph{Coxeter element} $\delta=s_1 s_2 \dots s_n$ in the Coxeter group $W=W_\Gamma$. 

\mk

Let $R$ denote the set of \emph{reflections} of $W$, i.e. the set of all conjugates of elements of $S$. Since $R$ generates $W$, we may consider its associated word norm $\|\cdot\|_R$. In the Cayley graph of $W$ with respect to $R$, let us consider the interval $P$ between $e$ and $\delta$: more precisely
$$P=\{u \in W \st \|u\|_R + \|u^{-1}\delta\|_R = \|\delta\|_R=n\}.$$
The set $P$ has a natural partial (prefix) order $\leq_L$: if $u,v \in P$, we say that $u \leq_L v$ if $\|u\|_R + \|u^{-1}v\|_R=\|v\|_R$. Equivalently, $u$ is a prefix of a minimal expression of $v$ as a product of reflections. Also equivalently, $u$ lies on a geodesic in the Cayley graph between $e$ and $v$.

\mk

The poset $P$ is easily seen to be bounded and graded. Let us define an interval-labeling $\lambda : I(P) \ra W$ by $\lambda(u,v) = u^{-1}v \in W$: this labeling is group-like and balanced. The poset $P$ is finite if and only if $W$ is finite, i.e. if $\Gamma$ is of spherical type.

\bdf[Dual Artin group]
The \emph{dual Artin group} associated to $\Gamma$ and $\delta$ is the interval group $A_\delta(\Gamma)$ of the poset $P$.
\edf

\bthm[Birman-Ko-Lee \cite{birman_ko_lee}, Bessis \cite{bessis}]
\label{thm:bessis}
If $\Gamma$ is of spherical type, for any Coxeter element $\delta$, the dual Artin group $A_\delta(\Gamma)$ is isomorphic to the standard Artin group $A_\Gamma$. Moreover, the poset $P$ is a lattice: in particular, the Artin group $A_\Gamma$ is a Garside group.
\ethm

\subsection{Complexes associated with Garside groups}
\label{subsec:garside_complex}
Consider a Garside group $G$, with positive monoid $G^+$, Garside element $\Delta$ and Garside generating set $S$ as in Definition~\ref{def:garside}. Let $\le_L$ and $\le_R$ be the orders as in Definition~\ref{def:garside}, which also extend to orders in $G$. More precisely, for $a,b\in G$, $a\le_L b$ if $b=ac$ for some $c\in G^+$, and $a\le_R b$ if $b=ca$ for some $c\in G^+$.

The \emph{Garside complex} of $G$ is the
simplicial complex $\widehat X_G$ with vertex set $G$, and with simplices corresponding to chains $g_1<_Lg_2<_L\cdots <_L g_n$ such that $g_n\le_L g_1\Delta$. Note that $G$ acts properly and cocompactly by simplicial automorphisms on its Garside complex. Alternatively, from the Garside group $G$, we can define an associated combinatorial Garside structure with the underlying poset $P$ being the set $\{e\}\cup S$ equipped with the partial order $\le_L$, and $\lambda(x,y)=x^{-1}y$ for $x,y\in P$. Then the universal cover of the interval complex associated with this combinatorial Garside structure is the Garside complex.

The \emph{Bestvina complex} of $G$ is the simplicial complex $X_G$ whose vertices correspond to left cosets of $\langle\Delta\rangle$ in $G$ (\cite{bestvina_artin}). There is an edge between
two vertices $a\langle\Delta\rangle$ and $b\langle\Delta\rangle$ if we can find $a'\in a\langle\Delta\rangle$ and $b'\in b\langle\Delta\rangle$ such that $b'=a'x$ for some $x\in S\setminus \{\Delta\}$. The Bestvina complex is the flag complex induced by its 1-skeleton. Note that $\bar G=G/\langle\Delta\rangle$ acts properly and cocompactly by simplicial automorphisms on the Bestvina complex. Topologically $\widehat X_G$ is homeomorphic to $X_G\times \mathbb R$.

\begin{thm} (\cite[Theorem~E]{haettel_helly_kpi1})
	\label{thm:garside_cub}
	For a Garside group $G$, if we metrize each simplex in the Garside complex $\widehat X_G$ as orthoschemes with $\ell^\infty$-metric as in \cite{haettel_helly_kpi1}, then $\widehat X_G$ is an injective metric space. In particular, it is CUB. Moreover, the injective metric on $\widehat X_G$ descends to a CUB metric on the Bestvina complex $X_G$, whose simplices are equipped with special polyhedral norms in the sense of \cite{haettel_simplicial_npc}.
\end{thm}

\section{Garside structure on $G\times \Z$}

\label{sec:general_construction}

\subsection{General construction}
\label{subsec:general}
We will now present a general construction of a Garside structure on the direct product $G \times \Z$, where $G$ is a group given by a specific presentation with generating set denoted $U$. We will consider $U$ as an abstract set endowed with a partial multiplication as defined below.

\bdf[Positive partial multiplication]
\label{def:partial mul}
Let $U$ denote a set. A map $\cdot$ defined on a subset of $U \times U$ with range $U$ is called a \emph{positive partial multiplication} if the following hold:
\bit
\item {\bf Left associativity} For any $u,v,w \in U$ such that $u \cdot v$ and $(u \cdot v) \cdot w$ are defined, we require that $v\cdot w$ and $u \cdot (v \cdot w)$ are defined, and that we have the equality $(u \cdot v) \cdot w =u \cdot (v \cdot w)$.
\item {\bf Right associativity} For any $u,v,w \in U$ such that $v\cdot w$ and $u \cdot (v \cdot w)$ are defined, we require that $u \cdot v$ and $(u \cdot v) \cdot w$ are defined, and that we have the equality $u \cdot (v \cdot w) = (u \cdot v) \cdot w$.
\item {\bf Identity} There exists a distinguished element $e \in U$ such that, for every $u \in U$, we have $e \cdot u = u \cdot e = u$.
\item {\bf Positivity} For any $u,v \in U$ such that $u \cdot v = e$, we have $u=v=e$.
\item {\bf Left cancellability} For any $u,v,w \in U$ such that $u \cdot v = u \cdot w$, we have $v=w$.
\item {\bf Right cancellability} For any $u,v,w \in U$ such that $v \cdot u = w \cdot u$, we have $v=w$.
\eit
\edf

Let us define relations $\leq_L,\leq_R$ on $U$ by:
\beq u \leq_L v & \mbox{ if there exists $w \in U$ such that $u \cdot w=v$} \\
u \leq_R v & \mbox{ if there exists $w \in U$ such that $w \cdot u=v$}.\eeq

\brk
Given $u,v \in U$, we will often write in the sequel "$u \cdot v \in U$" in place of "$u \cdot v$ is defined".
\erk

\blem \label{lem:order_on_U}
The relations $\leq_L$, $\leq_R$ are partial orders on $U$.
\elem

\bp
By the existence of $e \in U$, we know that both relations are reflexive.

By the associativity assumption, we know that both relations are transitive.

We will now prove that $\leq_L$ is antisymmetric, the proof for $\leq_R$ is similar. Let us assume that $u,v \in U$ are such that $u \leq_L v$ and $v \leq_L u$. There exists $w,w' \in U$ such that $v=u \cdot w$ and $u=v \cdot w'$, hence $u = (u \cdot w) \cdot w' = u \cdot (w \cdot w')$ by associativity. Since $U$ is cancellable, we deduce that $w \cdot w'=e$. Since $U$ is positive, we conclude that $w=w'=e$, hence $u=v$.
\ep

If $u,v\in (U,\leq_L)$ has a join (resp. meet), then we denote the join (resp. meet) by $u\vee_L v$ (resp. $u\wedge_Lv$). Similarly, we define $u\vee_R v$ and $u\wedge_R v$.

Note that the poset $(U,\leq_L)$ admits an interval-labeling with labels in $U$, i.e. for $u,v\in U$, the label of the interval between $u$ and $u \cdot v$ is $v\in U$. One readily verifies that this interval-labeling is group-like, so it makes sense to define the interval group $G_U$. In particular, $G_U$ has the following presentation:
$$G_U = \<U \st \forall u,v,w \in U \mbox{ such that } u \cdot v =w, \mbox{ we have } uv=w\>.$$

\mk

We will now describe the construction of a bounded poset $E$ consisting of one copy of $U$ and another ``inverted'' copy of $U$ as follows, which will be such that $G_E$ is isomorphic to $G_U \times \Z$.

\mk

Let $\bar U$ be another copy of $U$, and we denote $\bar u\in\bar U$ to be the element associated with $u\in U$. We will think $\bar u$ as a formal inverse of $u$.

\mk

Consider the set $E= (U,0) \sqcup (\bar U,1)$, with the following relation $\preceq$:
\bit
\item $(u,0) \preceq (v,0)$ if and only if $u \leq_L v$.
\item $(u,0) \preceq (\bar v,1)$ if and only if $v\cdot u \in U$.
\item $(\bar u,1) \preceq (\bar v,1)$ if and only if $v \leq_R u$.
\eit

\blem 
\label{lem:E partial order}
The relation $\prec$ is a partial order on $E$, with minimum $(e,0)$ and maximum $(\bar e,1)$. \elem

\bp
The reflexivity is clear. For transitivity, if $(u,0)\prec (\bar v,1)$ and $(\bar v,1)\prec (\bar w,1)$, then $v\cdot u\in U$ and $w\leq_R v$. Thus $v=w'\cdot w$ for some $w'\in U$. Thus $(w'\cdot w)\cdot u\in U$. By right associativity of the partial multiplication, we know $w\cdot u\in U$. Thus $(u,0)\prec (\bar w,1)$. Other cases of transitivity are similar. The antisymmetry of $\prec$ follows from the antisymmetry of $\leq_L$ and $\leq_R$ as in Lemma~\ref{lem:order_on_U}.
\ep

The poset $E$ is interval-labeled, with labels in $E$:
\bit
\item For $u,v\in U$, the label of the interval between $(u,0)$ and $(u\cdot v,0)$ is $(v,0) \in E$.
\item For $u,v,v\cdot u\in U$, the label of the interval between $(u,0)$ and $(\bar v,1)$ is $(\overline{v \cdot u},1) \in E$.
\item For $u,v,v\cdot u\in U$, the label of the interval between $(\overline{v\cdot u},1)$ and $(\bar u,1)$ is $(v,0) \in E$.
\eit

\blem 
\label{lem:grouplike}
The interval-labeled poset $E$ is group-like. \elem

\bp
Consider a chain with $3$ elements $a \prec b \prec c$ in $E$. Among the three labels $\lambda(a,b)$, $\lambda(a,c)$ and $\lambda(b,c)$, we will show that two of them determine the third one uniquely.

\mk

If $\lambda(a,b)$ and $\lambda(b,c)$ are known, there are three possibilities.
\bit
\item Assume that $\lambda(a,b)=(u,0)$ and $\lambda(b,c)=(v,0)$. Then $\lambda(a,c)=(u\cdot v,0) \in E$.
\item Assume that $\lambda(a,b)=(u,0)$ and $\lambda(b,c)=(\bar v,1)$. Then $a=(w,0)$ for $w\in U$, $b=(w\cdot u, 0)$ and $c=(\bar{x},1)$, where $x \cdot w \cdot u=v \in U$. Then $\lambda(a,c)=(\overline{x \cdot w},1)\in E$, and the value of $x\cdot w$ is determined by $u$ and $v$ through the equality $(x\cdot w)\cdot u=v$.
\item Assume that $\lambda(a,b)=(\bar u,1)$ and $\lambda(b,c)=(v,0)$. Then $a=(w,0)$ for some $w\in U$, $b=(\bar{x},1)$ and $c=(\bar{y},1)$ such that $u=x \cdot w \in U$ and $x= v \cdot y \in U$. Thus $u=v\cdot y\cdot w$. 
Since $a \prec c$, we know that $y \cdot w \in U$, so $\lambda(a,c)=(\overline{y \cdot w},1)$, and the value of $y\cdot w$ is determined by $u$ and $v$ through the equality $v\cdot (y\cdot w)=u$.
\eit

\mk

If $\lambda(a,b)$ and $\lambda(a,c)$ are known, there are three possibilities.
\bit
\item Assume that $\lambda(a,b)=(u,0)$ and $\lambda(a,c)=(v,0)$. Since $b \prec c$, there exists $w \in U$ such that $u \cdot w=v$. Such $w$ is unique by cancellability. Hence $\lambda(b,c)=(w,0) \in E$.
\item Assume that $\lambda(a,b)=(u,0)$ and $\lambda(a,c)=(\bar v,1)$. Then $\lambda(b,c)=(\overline{v \cdot u},1) \in E$.
\item Assume that $\lambda(a,b)=(\bar u,1)$ and $\lambda(a,c)=(\bar v,1)$. Then $a=(w,0)$ for some $w\in U$, $b=(\bar{x},1)$ and $c=(\bar{y},1)$, with $x,y \in U$ such that $x \cdot w=u$ and $y \cdot w=v$. Since $b \prec c$, there exists $z \in U$ such that $z \cdot y=x$. Hence $z \cdot y \cdot w = x \cdot w$, so $z \cdot v=u$. By cancellability, $z$ is uniquely determined by $u,v$. Then $\lambda(b,c)=(z,0)$.
\eit

\mk

By symmetry, the remaining case is similar.
\ep

\blem
\label{lem:balanced}
The interval-labeled poset $E$ is balanced.
\elem

\bp
The interval between $(u,0)$ and $(u \cdot v,0)$ has label $(v,0) \in E$, which is also the label of the interval between $(e,0)$ and $(v,0)$, and also between $(\bar v,1)$ and $(\bar e,1)$.

The interval between $(\ov{v \cdot u},1)$ and $(\bar u,1)$ has label $(v,0) \in E$, which is also the label of the interval between $(e,0)$ and $(v,0)$, and also between $(\bar v,1)$ and $(\bar e,1)$.

The interval between $(u,0)$ and $(\bar v,1)$ has label $(\ov{v \cdot u},1) \in E$, which is also the label of the interval between $(e,0)$ and $(\ov{v \cdot u},1)$, and also between $(v \cdot u,0)$ and $(\bar e,1)$.
\ep

\mk

\bpro
\label{prop:iso}
Let us consider the interval groups $G_U,G_E$ associated with the interval-labeled posets $U,E$. Then the natural map
\beq E & \mapsto & G_U \times \Z \\
(u,0) \in U \times \{0\} \subset E & \mapsto & (u,0) \in G_U \times \Z\\
(\bar{u},1) \in \bar U \times \{1\} \subset E & \mapsto & (u^{-1},1) \in G_U \times \Z\eeq
extends to an isomorphism of groups between $G_E$ and $G_U \times \Z$.
\epro

\bp
Note that $G_U \times \Z$ has generating set $(U \times \{0\}) \cup \{(e,1)\}$, and there are two types of relations:
\begin{enumerate}
	\item $(e,1)(u,0)=(u,0)(e,1)$ for each $u\in U$;
	\item $(u,0)(v,0)=(w,0)$ for any $u,v,w\in U$ with $u\cdot v=w$.
\end{enumerate}
On the other hand, the group $G_E$ has generating set $E$, and the relations are:
\begin{enumerate}
	\item $(u,0)(v,0)=(w,0)$ for any $u,v,w\in U$ with $u\cdot v=w$;
	\item $(u,0)(\bar v,1)=(\bar w,1)$ for any $u,v,w\in U$ with $w\cdot u=v$;
	\item $(\bar u,1)(v,0)=(\bar w,1)$ for $u,v,w\in U$ with $v\cdot w=u$.
\end{enumerate}
One readily checks that the map defined in the proposition extends to a group homomorphism $G_E\to G_U\times \mathbb Z$ as it is compatible with the relations.

We now define the inverse of this map on the standard generators of $G_U \times \Z$:
\beq \left(U \times \{0\}\right) \cup \{(e,1)\} \subset G_U \times \Z & \mapsto & G_E \\
(u,0) & \mapsto & (u,0) \in G_E\\
(e,1) & \mapsto & (\bar{e},1) \in G_E.\eeq
This map is also compatible with the relations of $G_U$ and $G_E$, note that we see relations of $G_U\times \Z$ of type (1) (as in the beginning of the proof) in $G_E$ as follows:
$$(u,0)(\bar e,1)=(u,0)((\bar u,1)(u,0))=((u,0)(\bar u,1))(u,0)=(\bar e,1)(u,0).$$
Note that the composition of these two maps are identity on the generators, thus they are inverses of each other. Then we are done.
\ep

Given $u,v\in U$, a \emph{left upper common bound} for $u,v$ is a upper bound for $\le_L$. A \emph{left join} of $u$ and $v$ is an element $w\in U$ with $u\le_L w$ and $v\le_L w$, such that  $w\le_L w'$ for any other left upper common bound $w'$ of $u,v$. A left join, if exists, must be unique. Similarly, we define right upper common bound and right join for $u,v$.

%A \emph{weak left join} of $u$ and $v$ is an element $w\in U$ with $u\le_L w$ and $v\le_L w$ such that there does not exist a left upper common bound $w'$ of $u,v$ such that $w'<w$. 

\bpro \label{pro:I_lattice}
Assume that $U$ satisfies the following:
\ben
\item $(U,\leq_L)$ and $(U,\leq_R)$ are weakly boundedly graded posets.
\item $(U,\leq_L)$ and $(U,\leq_R)$ are meet-semilattices.
\item For any $a,u,v,w \in U$ such that $a \cdot u,a \cdot v \in U$ and $w$ is the join for $\leq_L$ of $u$ and $v$, then $a \cdot w \in U$.
\item For any $a,u,v,w \in U$ such that $u \cdot a,v \cdot a \in U$ and $w$ is the join for $\leq_R$ of $u$ and $v$, then $w \cdot a \in U$.
\item For any $a,b,u,v \in U$ such that $a \cdot u,a \cdot v,b \cdot u,b \cdot v \in U$, either $a,b$ have a join for $\leq_R$, or $u,v$ have a join for $\leq_L$.
\een
Then $E$ is a lattice.
\epro

\bp
Assumption 1 implies $E$ is a weakly boundedly graded poset. More precisely, we define the rank function on $(U,0)$ using a rank function $f_1$ on $(U,\leq_L)$, and the rank function on $(U,\leq_R)$ using the an (inverted) rank function $f_2$ on $(U,\leq_R)$. As $f_1$ and $f_2$ has finite image, we can assume any image of $f_1$ is strictly smaller than each element in the image of $f_2$. So the new rank function on $E$ gives a weak grading, with finite image.
By Proposition~\ref{pro:bowtie}, it is sufficient to prove that $E$ contains no bowtie.

\mk

Assume that $(u,0),(v,0) \prec (w,0),(x,0)$ is a bowtie in $E$, where $u,v,w,x \in U$: hence $u,v \leq_L w,x$ is a bowtie in $U$, which contradicts that $(U,\leq_L)$ is a meet-semilattice by Proposition~\ref{pro:bowtie}.

Assume that $(\bar{a},1),(\bar{b},1) \leq (\bar{c},1),(\bar{d},1)$ is a bowtie in $E$, where $a,b,c,d \in U$: hence $a \in U \cdot c$, and $c \leq_R a$. So $c,d \leq_R a,b$ is a bowtie in $U$, which contradicts that $(U,\leq_R)$ is a meet-semilattice.

Assume that $(u,0),(v,0) \prec (w,0),(\bar{a},1)$ is a bowtie in $E$, where $u,v,w,a \in U$: hence $u,v \leq_L w$. Then $u$ and $v$ has a left join, denoted $u \vee_L v$, in $(U,\leq_L)$. Indeed, $u \vee_L v$ is the element with the lowest rank in $(U,\leq_L)$ which is an upped bound of $\{u,v\}$, and this element is unique since $(I,\leq_L)$ is a meet-semilatice. Since we are assumed to have a bowtie, $(w,0)$ is a minimal upper bound of $(u,0)$ and $(v,0)$ with respect to $\prec$, so $w$ is a minimal upper bound of $u$ and $v$ with respect to $\leq_L$. Thus $w=u \vee_L v$. As $(u,0),(v,0) \prec (\bar{a},1)$, we know $a \cdot u,a \cdot v \in U$. By assumption (3), this implies that $a \cdot w \in U$, so $(w,0) \prec (\bar{a},1)$, which contradicts that $(u,0),(v,0) \prec (w,0),(\bar{a},1)$ is a bowtie in $E$.

Assume that $(a,0),(\bar{w},1) \prec (\bar{u},1),(\bar{v},1)$ is a bowtie in $E$, where $u,v,w,a \in U$: hence $u,v \leq_R w$. Similarly to the previous paragraph, $u$ and $v$ has a right join, denoted $u \vee_R v$, in $(U,\leq_R)$; and $w=u \vee_R v$. As $(a,0)\prec (\bar{u},1),(\bar{v},1)$, we have $u \cdot a,v \cdot a \in U$. By assumption, this implies that $w \cdot a \in U$, so $(\bar{w},1) \prec (a,0)$, which contradicts that $(a,0),(\bar{w},1) \prec (\bar{u},1),(\bar{v},1)$ is a bowtie in $E$.

Assume that $(u,0),(v,0) \prec (\bar{a},1),(\bar{b},1)$ is a bowtie in $E$, where $u,v,a,b \in U$: hence $a \cdot u,a \cdot v,b \cdot u,b \cdot v \in U$. By assumption, this implies that either $a,b$ have a right join $c$ for $\leq_R$ or $u,v$ have a left join $w$ for $\leq_L$. In the former case, we deduce from assumption (4) of the proposition that $c\cdot u\in U$ and $c\cdot v\in U$. Hence $(\bar c,1)$ is in the middle of the bowtie, i.e. $\{(u,0),(v,0)\}\prec(\bar c,1)\prec {(\bar a,1),(\bar b,1)}$, which is a contradiction. Similarly, in the latter case, we deduce from assumption (3) of the proposition that $(w,0)$ is in the middle of the bowtie, which is a contradiction.
\ep

The following is combination of Definition~\ref{def:comb Gar}, Theorem~\ref{thm:gar}, Theorem~\ref{thm:garside_cub}, Lemma~\ref{lem:grouplike}, Lemma~\ref{lem:balanced} and Proposition~\ref{pro:I_lattice} (we refer to Definition~\ref{def:injective} and Definition~\ref{def:bicombing} for definitions of injective metric spaces and CUB spaces).

\bthm \label{thm:I_lattice_implies_Garside}
Under the assumption of Proposition~\ref{pro:I_lattice}, the group $G_E \simeq G_U \times \Z$ is quasi-Garside with Garside element $(e,1)$ and set of simple elements $E$. If $U$ is finite, then $E$ is finite and $G_E$ is Garside. 

Moreover, suppose $U$ is finite. Let $K_E$ denote the interval complex of the labeled poset $E$, whose universal cover is isomorphic to the Garside complex $\widehat X_{G_E}$ of the Garside group $G_E$ (as discussed in Section~\ref{subsec:garside_complex}). If we metrize each simplex in $\widehat X_G$ as orthoschemes with $\ell^\infty$-metric as in \cite{haettel_helly_kpi1}, then $\widehat X_G$ is an injective metric space. In particular, it is CUB. 
\ethm

In particular, we can deduce the following result stated in the introduction.

\bthm \label{thm:simple_criterion_U_garside}
Let $U$ be a finite set, endowed with a positive partial multiplication, and let $G_U$ denote the associated interval group. Assume that the following hold:
\bit
\item $(U,\leq_L)$ and $(U,\leq_R)$ are meet-semilattices.
\item For any $a,u,v,w \in U$ such that $a \cdot u,a \cdot v \in U$ and $w$ is the join for $\leq_L$ of $u$ and $v$, then $a \cdot w \in U$.
\item For any $a,u,v,w \in U$ such that $u \cdot a,v \cdot a \in U$ and $w$ is the join for $\leq_R$ of $u$ and $v$, then $w \cdot a \in U$.
\item For any $a,b,u,v \in U$ such that $a \cdot u,a \cdot v,b \cdot u,b \cdot v \in U$, either $a,b$ have a join for $\leq_R$, or $u,v$ have a join for $\leq_L$.
\eit
Then the group $G_U \times \Z$ is a Garside group, with Garside element $(e,1)$.
\ethm

\brk
We may remark that there are very simple situations where we can apply Theorem~\ref{thm:simple_criterion_U_garside}. For instance, let us consider the free group $F$ over a finite set $S$, and let $U=S \cup \{e\}$. Then $U$ satisfies the assumptions of Proposition~\ref{pro:I_lattice}, and in particular the group $F \times \Z$ is Garside. This particular case can also be deduced from \cite{picantin2022cyclic} via different methods.
We will however see, in the rest of the article, more interesting applications of this result.
\erk

\subsection{Some examples where Theorem~\ref{thm:simple_criterion_U_garside} applies}

Recall that for a group $G$ with a given presentation, there is a \emph{presentation complex}, which is a 2-dimensional cell complex with the 0-skeleton being a single vertex, the 1-skeleton being a wedge of oriented circles in 1-1 correspondence with elements in the generating set (each circle is labeled by a generator), and for each relator we attached a 2-cell such that the boundary of this 2-cell goes to the loop in the 1-skeleton tracing out the word of this relator. The presentation complex has fundamental group being $G$.

For a (not necessarily simplicial) graph, its \emph{girth} is defined to the minimal possible number of edges in an embedded $\mathbb S^1$ in the graph. 

\bthm \label{thm:presentation_girth_Garside}
Let us consider a group $G$ given by a finite presentation $\<S \st r_1=r'_1,\dots,r_n=r'_n\>$. Assume that the following hold:
\bit
\item For each $1 \leq i \leq n$, the words $r_i,r'_i$ are positive words in $S$, without common prefix or suffix.
\item For each $1 \leq i \leq n$, the word $r_i$ (and $r'_i$) does not appear as a subword of some of the other $2n-1$ words.
\item For each distinct $s,t \in S$, there exist at most one $1 \leq i \leq n$ such the the first letters of $\{r_i,r'_i\}$ are $\{s,t\}$.
\item For each distinct $s,t \in S$, there exist at most one $1 \leq i \leq n$ such the the last letters of $\{r_i,r'_i\}$ are $\{s,t\}$.
\item The presentation is $T(5)$, i.e. the link of the vertex in the presentation complex has girth at least $5$. 
\eit
Then $G \times \Z$ is Garside.
\ethm

\bp
Let $\F^+(S)$ be the positive free monoid on the generating set $S$.
Let $U$ denote the quotient of the set of subwords of the words $R=\{r_1,r'_1,\dots,r_n,r'_n\}$ in $\F^+(S)$, under the equivalence relation defined by $r_i \sim r'_i$, for each $1 \leq i \leq n$. Given two positive words $u_1,u_2$ of $\F^+(S)$, we will write $u_1=u_2$ if they are the same word in $\F^+(S)$, and $u_1\equiv u_2$ if they give the same element in $U$. Note that if $u$ is a proper non-trivial subword of $r_i$ and $u'$ is a proper non-trivial subword of $r'_i$, then $u\not\equiv u'$ - this follows from the second bullet point of the assumption of the theorem.

Let us endow $U$ with the partial multiplication induced by the free monoid $\F^+(S)$ on $S$. 
For any $u \in U$, let us define $u \cdot e\equiv e \cdot u\equiv u$. For any $u,v \in U \bs \{e\}$, then $u \cdot v$ exists and is equal to $uv \in \F^+(S)$ if and only if there exists $r \in R$ such that $uv$ is a subword of $r$.

We will show that this defines a positive partial multiplication on $U$.

\mk

Assume that $u,v,w \in U$ are such that $u \cdot v \in U$ and $(u \cdot v) \cdot w \in U$. We will consider $u,v,w$ as representatives inside $\F^+(S)$. Then there exists $r_1\in R$ such that $uv$ is a subword of $r_1$. Moreover, $uv,w\notin R$ and $uvw$ is a subword of $r\in R$. Hence $vw$ is a subword of $r$. Note that $v\notin R$, otherwise we will contradict the second assumption of the theorem. Thus by the definition of product, $v\cdot w\equiv vw\in U$.
So $U$ is left associative, and similarly we can show it is right associative.

\mk

The identity element $e \in \F^+(S)$ is an identity element for $(U,\cdot)$. Now we check positivity. If $u\cdot v\in U$ and $u\cdot v\equiv e$, then $u,v\notin R$ and $uv$ is a subword of $r\in R$. This forces $u=v=e$ as $u$ and $v$ are positive words.
For the cancellability, if $u\cdot v\equiv u\cdot w$, then either $uv=uw$ in $\F^+(S)$ and $v=w$ follows from the cancellability in $\F^+(S)$, or there exists $i$ such that $uv=r_i$ and $uw=r'_i$, which implies $u=e$ as we assume that for each $1 \leq i \leq n$, the words $r_i,r'_i$ have no common prefix or suffix.

\mk

So $U$ satisfies Definition~\ref{def:partial mul}. In particular, we may consider the left and right orders on $U$.

\mk

Let us prove that $(U,\leq_L)$ is weakly boundedly graded. It is clear that $e \in U$ is the minimum of $U$. For each $u \in U$, let $r(u) \in \N$ denote the maximal length of a representative for $u$ in $\F^+(S)$ (note that the second bullet point in the statement of the theorem implies that we can not keep replacing subwords in $u$ by longer subwords, so there is a bound on the length of all representatives). The map $r:U \ra \N$ is a rank map with respect to $\leq_L$ and $\leq_R$. Since $R$ is finite, $r$ has finite image. So $(U,\leq_L)$ and $(U,\leq_R)$ are weakly boundedly graded.

\mk

Let us prove that $(U,\leq_L)$ is a meet-semilattice. According to Proposition~\ref{pro:bowtie}, it is sufficient to prove that $(U,\leq_L)$ does not contain a bowtie. By contradiction, assume that $u,v \leq_L x,y$ is a bowtie in $U$, with $r(x)+r(u)+r(v)+r(y)$ minimal, where $r:U \ra \N$ is a weak grading. Let $s,t \in S$ denote the first letters of $u,v$ respectively. If $s=t$, then $s$ is a common lower bound (with respect to $\leq_L$) of $\{u,v,x,y\}$, which enables us to find a bowtie with a smaller value of $r(x)+r(u)+r(v)+r(y)$. Thus $s \neq t$. We have $s,t \leq_L x,y$. Thus $x$ can be represented by two different positive words $x_s$ and $x_t$ in the free monoid, with $x_s$ starting with $s$ and $x_t$ starting with $t$. If $x_s\notin R$, then $x_s$ is a proper non-trivial subword of some element $r$ in $R$ (by the second bullet point of the theorem). However, we can transform $x_s$ to $x_t$ by using relators finitely many times, this means a subword of $x_s$ (hence a subword of $r$) belongs to $R$, which contradicts the second bullet point of the theorem. Hence $x_s \in R$, and similarly $y$ represents a word in $R$. According to the third bullet point of the theorem, we deduce that there exists $1 \leq i \leq n$ such that $x,y$ both represent $r_i \equiv r'_i$. Hence $x=y$ in $U$, so $u,v \leq_L x,y$ is not a bowtie: this is a contradiction. Hence $(U,\leq_L)$ is a meet-semilattice.

\mk

Let $a,u,v \in U$ such that $a \cdot u,a \cdot v \in U$ and $u \vee_L v=w \in U$. We want to prove that $a \cdot w \in U$. We will actually prove that $u \leq_L v$ or $v \leq_L u$: if not, this means that there exists $1 \leq i \leq n$ such that $u,v$ are prefixes of $r_i,r'_i$ respectively (up to switching $u$ and $v$). Then the words $au,av,u^{-1}v$ give rise to a triangle in the link of the vertex in the presentation complex, which contradicts the $T(5)$ assumption. So $u \leq_L v$ or $v \leq_L u$, and hence $a \cdot w \in \{a \cdot u,a \cdot v\} \subset U$.

\mk

Let us assume that $a,b,u,v \in U$ are pairwise distinct such that $a \cdot u, a \cdot v, b \cdot u, b \cdot v \in U$. We will prove that either $a,b$ are comparable for $\leq_R$, or $u,v$ are comparable for $\leq_L$. If not, then the words $au,av,bv,bu$ give rise to a $4$-cycle in the link of the vertex in the presentation complex, which contradicts the $T(5)$ assumption. So for instance $a,b$ are comparable for $\leq_R$, in which case $a$ and $b$ have a join for $\leq_R$.

\mk

According to Proposition~\ref{pro:I_lattice}, we deduce that $G \times \Z$ is a Garside group.
\ep

\bcor \label{cor:surface_group_garside}
For any surface $S$ of finite type (possibly non-orientable), except the projective plane, $\pi_1(S) \times \Z$ is a Garside group.
\ecor

\bp
If $S$ is a surface with boundary, its fundamental group is a free group.

\mk

If $S$ is the torus, then $\pi_1(S) \simeq \Z^2$, which is a Garside group, so $\pi_1(S) \times \Z \simeq \Z^3$ is a Garside group.

\mk

If $S$ is the closed orientable surface with genus $g \geq 2$, consider the standard presentation, see e.g. \cite[Section 1.2]{hatcher2002algebraic},
$$G_g =\<a_1,b_1,\dots,a_g,b_g \st [a_1,b_1] \dots [a_g,b_g]=1\>.$$
This presentation is not positive, so we will modify it as follows, where $\bar a_1$ and $\bar b_1$ play the role of $a^{-1}_1$ and $b^{-1}_1$:
\beq G_g&=&\<\bar a_1,\bar b_1,a_2,b_2,\ldots,a_g,b_g, h_2,h_3,\dots,h_{g-1}\st \bar a_1\bar b_1h_2h_3 \dots h_{g-1}a_gb_g =\bar b_1\bar a_1b_ga_g,\\
&& a_2b_2=h_2b_2a_2, a_3b_3=h_3b_3a_3,\dots, a_{g-1}b_{g-1}=h_{g-1}b_{g-1}a_{g-1}\>.\eeq
Note that these two presentations give the same group - $\bar a_1$ and $\bar b_1$ go to $a^{-1}_1$ and $b^{-1}_1$, and $h_i$ goes to the commutator of $a_i$ and $b_i$ for $2\le i\le g-1$.
Then this new presentation satisfies the assumptions of Theorem~\ref{thm:presentation_girth_Garside}.

\mk

If $S$ is the projective plane, then $\pi_1(S) \simeq \Z/2\Z$, so $\pi_1(S) \times \Z$ has torsion, hence it is not a Garside group.

\mk

If $S$ is the closed non-orientable surface with genus $2$, i.e. the Klein bottle, then its fundamental group has the following  presentation
$$\pi_1(S) =\<a,b \st a^2=b^2\>,$$
which is a Garside presentation with Garside element $\Delta=a^2=b^2$. Hence $\pi_1(S) \times \Z$ is a Garside group.

\mk

If $S$ is the closed non-orientable surface with genus $g \geq 3$, consider the (almost) standard presentation (see e.g. \cite[Section 1.2]{hatcher2002algebraic})
$$G =\<a_1,\dots,a_g \st a_1^2 \dots a_{g-1}^2=a_g^2\>,$$
then one readily checks that it satisfies the assumptions of Theorem~\ref{thm:presentation_girth_Garside}. Here we explain the verification of the $T(5)$ condition. Note that the link of the unique vertex in the presentation complex has $2g$ vertices - for each generator $a_i$, it gives two vertices in the link, denoted by $a^+_i$ and $a^-_i$, corresponding to the outgoing and incoming direction of the oriented loop corresponding to $a_i$. The each $i$, there is an edge joining $a^+_i$ and $a^-_i$. For $1\le i\le g-2$, there is an edge joining $a^-_i$ and $a^+_{i+1}$. There is an edge joining $a^-_{g-1}$ and $a^-_{g}$, and an edge joining $a^+_g$ and $a^+_1$. These are all the edges. Thus the link has girth $2g$. Since $g\ge 3$, we know the presentation is $T(5)$.
\ep

Another easy class of groups for which we can apply Theorem~\ref{thm:I_lattice_implies_Garside} is the following class of groups with a systolic presentation. They have been defined and studied by Soergel in~\cite{soergel_systolic}.

Recall that  a flag simplicial complex is \emph{systolic} if it is simply connected and all of its vertex links are 6-large, i.e.
all cycles of length 4 or 5 have diagonals, see \cite{januszkiewicz2006simplicial}.

\bdf[Soergel~\cite{soergel_systolic}]
\label{def:systolic}
A finite presentation $\<S \st R\>$ of a group $G$ is called a \emph{systolic restricted presentation} if the following hold:
\bit
\item Each relation $r \in R$ is of the form $r=abc^{-1} \in \F(S)$, where $a,b,c \in S$ and $\F(S)$ is the free group on $S$.
\item The flag completion of the Cayley graph of $G$ with respect to $S$ is simplicial and systolic.
\eit
\edf

Note that asking that the Cayley graph of $G$ with respect to $S$ is simplicial is equivalent to asking that any $s \in S$ has image in $G$ different from $e$, and also for any distinct $s,t \in S$, their image in $G$ are neither equal nor inverse. Soergel gives a complete characterization of such systolic restricted presentations in~\cite[Theorem~1]{soergel_systolic}.

\mk

Among Garside presentations, Soergel gives a characterization of those which are systolic, see~\cite[Theorem~2]{soergel_systolic}. There are essentially amalgams of the following Garside groups $G_{n,m}$, for $n,m \geq 1$, defined by the following systolic restricted presentation:
$$G_{n,m} = \<x_1,\dots,x_n \st x_1x_2 \dots x_m=x_2x_3 \dots x_{m+1}= x_nx_1x_2 \dots x_{m-1}\>.$$
While the presentation of $G_{n,m}$ is not systolic restricted on the nose, \cite[Theorem~2]{soergel_systolic} proves that these groups, as well as some amalgamations of them, actually have systolic restricted presentations.

\mk

Among 2-dimensional Artin groups, Soergel gives a sufficient criterion in terms of orientations of the edges of the Coxeter presentation graph, see~\cite[Theorem~3]{soergel_systolic}. As a very restricted example, if the Coxeter presentation graph $\Gamma$ has no triangles and no squares, then $A_\Gamma$ admits a systolic restricted presentation. 

\bthm
\label{cor_systolic_garside}
Let $G$ denote a group with a systolic restricted presentation. Then $G \times \Z$ is a Garside group.
\ethm

\bp
Let us denote by $U=S \cup \{e\}$ in $G$, and let us consider the induced partial multiplication from $G$. Since the Cayley graph of $G$ with respect to $S$ is simplicial, we deduce that $U$ embeds in $G$. We will denote the flag completion of the Cayley graph of $G$ by $X$.

\mk

The only non-trivial assumption to check for this partial multiplication is the positivity: if there exist $s,t \in S$ such that $st=e$ in $G$, this contradicts the fact that the Cayley graph of $G$ with respect to $S$ is simplicial.

\mk

Since $U$ is finite, it is weakly boundedly graded, and it has minimum $e$. Let us show that $(U,\leq_L)$ is a meet-semilattice by contradiction: let us assume that we have a bowtie $a,b <_L u,v$, with $a,b,u,v \in S$. Then this corresponds to a loop of length $4$ in the link of the vertex $e$ in $X$. By systolicity, we deduce that there exists a diagonal: either $a,b$ are comparable, or $u,v$ are comparable. Hence $a,b <_L u,v$ is not a bowtie. Similarly, $(U,\leq_R)$ is a meet-semilattice.

\mk

Let us now consider $a,u,v \in U$ such that $au,av \in U$ and $u,v$ have a join $w \in U$ for $\leq_L$. Then $a^{-1},u,w,v$ form a loop of length $4$ in the link of the vertex $e$ in $X$. By systolicity, we deduce that there exists a diagonal: either $u,v$ are comparable, in which case $w \in \{u,v\}$ and $aw \in \{au,av\} \subset U$, or there is an edge between $a^{-1}$ and $w$, in which case $aw \in U$.

\mk

Let us now consider $a,b,u,v \in U$ such that $au,av,bu,bv \in U$. Then $a^{-1},u,b^{-1},v$ form a loop of length $4$ in the link of the vertex $e$ in $X$. By systolicity, we deduce that there exists a diagonal. If there is an edge between $a^{-1}$ and $b^{-1}$, this means that $a$ and $b$ are $\leq_R$-comparable, so they have a right join. If there is an edge between $u$ and $v$, this means that $u$ and $v$ are $\leq_L$-comparable, so they have a left join.

\mk

According to Theorem~\ref{thm:I_lattice_implies_Garside}, we conclude that $G \times \Z$ is a Garside group.
\ep

\begin{defi}
	\label{def:square presentation}
A finite presentation $\langle S\mid R\rangle$ is a \emph{positive square presentation} if each relator $r$ is of form $ab=cd$ where $a,b,c,d$ are (not necessarily distinct) elements in $S$.
\end{defi}

Some natural examples of groups of square presentation include right-angled Artin groups, mock right-angled Artin groups in the sense of \cite{scott_mock} and groups arising from word labeled oriented graphs in the sense of \cite{harlander2015aspherical}. We give a criterion showing some of these groups are Garside groups after taking the product with $\Z$.

Given a finite square presentation, let $X$ be the associated presentation complex. Each edge loop of $X$ is oriented and labeled by an element in $S$. Let $\Theta$ be the link of the unique vertex of $X$. A vertex of $\Theta$ is of type $o$ or $i$ if it corresponds to outgoing or incoming edge at the base vertex. 

\bthm
\label{thm:square garside}
Let $G$ denote a group with a positive square presentation such that the link $\Theta$ of its presentation complex satisfies the following conditions:
\begin{enumerate}
	\item there does not exist embedded 2-cycles in $\Theta$ of type $(o,o)$ (means a 2-cycle with two vertices of type $o$) or type $(i,i)$;
	\item there does not exist embedded 3-cycles in $\Theta$ of type $(o,o,i)$ or $(i,i,o)$;
	\item there does not exist embedded 4-cycles in $\Theta$ of type $(o,i,o,i)$.
\end{enumerate}
Then $G\times \Z$ is a Garside group.
\ethm

\begin{proof}
Suppose the collection of relators are of form $\{a_ib_i=a'_ib'_i\}_{i=1}^k$ where $a_i,b_i,a'_i,b'_i\in S$ for each $1\le i\le k$. Let $U$ be the set of equivalence classes of words in $\{e\}\cup S\cup \{a_ib_i\}_{i=1}^k\cup \{a'_ib'_i\}_{i=1}^k$, under the equivalence relation generated by $a_ib_i\sim a'_ib'_i$ for $1\le i\le k$. We endow $U$ with the partial multiplication as in the proof of Theorem~\ref{thm:presentation_girth_Garside}.
Now we verify the assumptions of Theorem~\ref{thm:simple_criterion_U_garside}: $(U,\le_L)$ is a meet-semilattice follows from the lack of 2-cycles of type $(o,o)$ in $\Theta$. Indeed, a bowtie in $(U,\le L)$ consists of two elements $a,b\in S$, such that they have two different upper bounds in $U$. Thus in a small neighborhood of the base vertex, we see a small outgoing $a$ segment, and a small outgoing $b$ segment, span the corners of two different squares. This gives rise a 2-cycle of type $(o,o)$ in the link of the base vertex, which contradicts assumption (1) of the theorem. Similarly, we know that $(U,\le_R)$ is a meet-semilattice follows from the lack of 2-cycles of type $(i,i)$ in $\Theta$. Now take $a,u,v\in S$ with $u\neq v$ such that $au,av\in U$ and $u$ and $v$ have a left join, then this gives a 3-cycle in $\Theta$ made of vertices of type incoming $a$, outgoing $u$, outgoing $v$, which is excluded by the lack of $3$-cycles of type $(o,o,i)$. Similarly, the third item of Theorem~\ref{thm:simple_criterion_U_garside} follows from the lack of $3$-cycles of $(i,i,o)$. For the last item of Theorem~\ref{thm:simple_criterion_U_garside}, let $a,b,u,v\in S$ with $u\neq v$ and $a\neq b$. If $au,av,bu,bv\in U$, this gives a 4-cycle in $\Theta$ with consecutive vertices of type $a$ incoming, $u$ outgoing, $b$ incoming and $v$ outgoing, which is ruled out by the lack of 4-cycles of type $(o,i,o,i)$. Thus we are done by Theorem~\ref{thm:simple_criterion_U_garside}.
\end{proof}
\section{Application to Artin groups}
\label{sec:Artin}

We will now explain how we can apply Theorem~\ref{thm:simple_criterion_U_garside} for some Artin groups.

Let $(W,S)$ denote a Coxeter group, and let $A$ denote the associated Artin group. Let $R$ denote the set of all reflections (i.e. conjugates of elements in $S$) of $W$. A \emph{reflection decomposition} of an element $w\in W$ is a way of written $w$ as a product of elements in $S$.  A reflection decomposition of $w$ is \emph{minimal} if the length of the product is as small as possible. Let $|w|$ denote the reflection length of $w$, i.e. the length of a minimal reflection decomposition of $w$. We define two relations $\le_L$ and $\le_R$ in $W$ as follows. We say $u\le_L v$ if $u$ appears as a prefix of a minimal reflection decomposition of $v$; and $u\le_R v$ if $u$ appears as a suffix of a minimal reflection decomposition of $v$. One readily sees that $\le_R$ and $\le_L$ are partial orders on $W$.

We will define the set $U$ inside $W$, and we will look for conditions on $U$ and $W$ ensuring that the assumptions from Theorem~\ref{thm:simple_criterion_U_garside} are satisfied.

\mk

Given a finite subset $U \subset W$, we consider the following partial multiplication $\cdot$ on $U$: if $u,v \in U$ are such that their product $uv$ in the Coxeter group $W$ lies in $U$ and furthermore $|u v| = |u|+|v|$, we define $u \cdot v=uv \in U$. Let $R_U = R \cap U$, i.e. the collection of reflections that are in $U$. Similarly, if $u_1,u_2,\ldots,u_n\in U$, we will write $w=u_1\cdot u_2\cdot\ldots\cdot u_n\in U$ if $u_1u_2\cdots u_n\in U$ and $|w|=|u_1|+|u_2|+\ldots+|u_n|$.

\blem
\label{lem:U partial order}
Suppose $U$ satisfies the following conditions:
\begin{enumerate}
\item For every $u \in U$, there exist $r_1,\dots,r_n \in R_U$ such that $u=r_1\cdot r_2 \cdot \ldots \cdot r_n$.
\item For every $r_1,\ldots,r_n \in R_U$ such that $r_1\cdot r_2 \cdot \ldots \cdot r_n \in U$, we have $r_1\cdot r_2 \cdot \ldots \cdot r_{n-1} \in U$ and $r_2\cdot r_3 \cdot \ldots \cdot r_n \in U$.
\end{enumerate}
Then the set $(U,\cdot)$ satisfies Definition~\ref{def:partial mul}. 
\elem

\bp
It suffices to verify that
for $u,v,w \in U$ such that $u\cdot v\cdot w \in U$, we have $u\cdot v\in U$ and $v\cdot w\in U$. Indeed, by Assumption 1, let us write reflection factorizations in $R_U$: $u=r_1\cdot r_2 \cdot \ldots \cdot r_n$, $v=r'_1\cdot r'_2 \cdot \ldots \cdot r'_{n'}$ and $w=r''_1\cdot r''_2 \cdot \ldots \cdot r''_{n''}$. We then have $u\cdot v\cdot w \in U$, so by assumption 2, we have both $u\cdot v \in U$ and $v\cdot w \in U$.
\ep

Note that the relations $\le_L$ and $\le_R$ on $U$ defined right after Definition~\ref{def:partial mul} is the same as the restriction of the relations $(W,\le_L)$ and $(W,\le_R)$ (as defined in the beginning of this section) to $U$. So there is no ambiguity in our notation.

\mk
Let $E$ be the poset constructed from $U$ as in Section~\ref{subsec:general}.
We have a criterion for $E$ to be a lattice.

\bpro \label{pro:E_lattice}
Suppose $U\subset W$ is finite. Assume that all the following requirements are satisfied.
\ben
\item \label{p:factorization} For every $u \in U$, there exist $r_1,\dots,r_n \in R_U$ such that $u=r_1\cdot r_2 \ldots  r_n$.
\item \label{p:prefix_suffix} For every $r_1,\ldots,r_n \in R_U$ such that $r_1\cdot r_2 \cdot \ldots \cdot r_n \in U$, we have $r_1\cdot r_2 \cdot \ldots \cdot r_{n-1} \in U$ and $r_2\cdot r_3 \cdot \ldots \cdot r_n \in U$.
\item Now we consider the partial orders $\le_L$ and $\le_R$ restricted on $U$. We require that these partial order satisfy that, for every $r_1\in R_U$ and $r_2\in R_U$ with a common left upper bound in $(U,\le_L)$, they have a left join in $(U,\le_L)$; similarly, if $r_1$ and $r_2$ have a common right upper bound in $(U,\le_R)$, then they have a right join in $(U,\le_R)$.
\item For every $a \in U$, for any $u,v \in R_U$ such that $u,v$ have a left join $w\in (U,\le_L)$ and $a\cdot u,a\cdot v \in U$, we have $a\cdot w \in U$.
\item For every $a \in U$, for any $u,v \in R_U$ such that $u,v$ have a right join $w\in (U,\le_R)$ and $u\cdot a,v\cdot a \in U$, we have $w\cdot a \in U$.
\item For every $a,b,u,v \in R_U$ and any $x \in U$ such that $a\cdot x\cdot u,a\cdot x\cdot v,b\cdot x\cdot u,b\cdot x\cdot v \in U$, we have that either $u,v$ have a left join, or $a,b$ have a right join. 
\een
Then $U$ satisfies all the conditions in Proposition~\ref{pro:I_lattice}. In particular $E$ is a lattice.
\epro

\bp
By Lemma~\ref{lem:order_on_U}, it remains to verify that $U$ satisfies the assumptions of Proposition~\ref{pro:I_lattice}.

\mk

Assumption 1 of Proposition~\ref{pro:I_lattice} follows by considering the reflection length on $U$ (the bounded assumption follows from finiteness of $U$).

\mk

We will now prove that $(U,\leq_L)$ is a meet-semilattice. We artificially add a largest element $\hat 1$ to $U$, so $P=(U\cup \{\hat 1\},\le_L)$ is a bounded poset of finite length (i.e. there is a finite upper
bound on the lengths of its chains). Recall that an element $p_1\in P$ \emph{covers} $p_2\in P$ if $p_1>p_2$ and there does not exist $p\in P$ with $p_1>p>p_2$.
We claim that if $u_1,u_2$ are two distinct elements in $P$ that covers $v$, then $u_1$ and $u_2$ has a join. By assumption 2, we can write $u_1=v\cdot r_1$ and $u_2=v\cdot r_2$ with $r_1,r_2\in R_U$. If $\hat 1$ is the only common left upper bound of $u_1$ and $u_2$, then clearly they have a join. If $u_1$ and $u_2$ have a common left upper bound $u'$ other than $\hat 1$, then we can write $u'=u_i\cdot w_i$ with $w_i\in U$ for $i=1,2$. Thus $r_1$ and $r_2$ has a common left upper bound, which is $w'=r_1\cdot w_1=r_2\cdot w_2\in U$. Let $r$ be the left join of $r_1$ and $r_2$. Then $r\le_L w'$, which implies that $v\cdot r\le_L v\cdot w'=u'$. Thus $v\cdot r$ is the join for $u_1$ and $u_2$. Now it follows from \cite[Lemma 2.1]{bjorner1990hyperplane} that $P$ is a lattice. Thus $(U,\le_L)$ is a meet-semilattice. Similarly we can prove $(U,\le_R)$ is a meet-semilattice.

\mk

We now prove Assumption 3 of Proposition~\ref{pro:I_lattice}, i.e. for all $a,u,v \in U$ such that $u,v$ have a join $w \in U$ for $\leq_L$, and $a\cdot u,a\cdot v \in U$, then $a\cdot w \in U$.

We will prove it by decreasing induction on $|a|$, and for a fixed value of $|a|$ by increasing induction on $|u|+|v|$. Since $U$ is finite, if $|a|$ is maximal, then $u,v =e$, so the property is true. Now consider $a,u,v \in U$, and assume that the property holds for any larger value of $|a|$. If $u,v \in R_U$, the property holds by assumption. So let us assume that the property holds for smaller values of $|u|+|v|$. Let us assume that $u \not\in R_U$, and write $u=u_1\cdot r$, with $u_1 \in U$, $r \in R_U$ and $|u_1|=|u|-1$. According to Properties~\ref{p:factorization} and \ref{p:prefix_suffix}, we know that $a\cdot u_1 \in U$. Since $u_1,v$ have an upper bound $w$, and since $U$ is a meet-semilattice, they have a left join $w_1=u_1\cdot w'$, with $w' \in U$. Since $|u_1|+|v| < |u|+|v|$, we deduce by induction that $a\cdot w_1 = a\cdot u_1\cdot w' \in U$. Now $u_1\cdot w'$ and $u=u_1\cdot r$ have a left upper bound $w \in U$, so we deduce by Lemma~\ref{lem:U partial order}, Properties~\ref{p:factorization} and \ref{p:prefix_suffix} that $w'$ and $r$ have a left upper bound in $U$, hence they also have a left join: let us write $w' \vee_L r=w'' \in U$. Since $|au_1| > |a|$, we deduce by induction that $au_1\cdot w'' \in U$, in particular $|au_1w''|=|au_1|+|w''|=|a|+|u_1|+|w''|$. Note that $u=u_1r\le_L u_1\cdot w''$ and $v\le_L w_1=u_1w'\le_L u_1w''$, we know $w \leq_L u_1\cdot w''$.
On the other hand, $u^{-1}_1w$ is a left common upper bound for $r$ and $w'$. Hence $w''\le_L u^{-1}_1 w$ and $u_1\cdot w''\le_L w$. Then $w=u_1\cdot w''$. In particular $|aw|=|au_1w''|=|a|+|u_1|+|w''|=|a|+|w|$. As $aw=wu_1w''\in U$, we know $a\cdot w\in U$.

Assumption 4 of Proposition~\ref{pro:I_lattice} can be proved in a similar way.

\mk

We will now prove Assumption 5 of Proposition~\ref{pro:I_lattice}, i.e. for every $a,b,u,v,x \in U$ such that $a\cdot x\cdot u,a\cdot x\cdot v,b\cdot x\cdot u,b\cdot x\cdot v \in U$, we have that either $u,v$ have a left join, or $a,b$ have a right join. 

We will prove it by decreasing induction on $|x|$, and for a fixed value of $|x|$ by increasing induction on $|a|+|b|+|u|+|v|$. Since $U$ is finite, if $|x|$ is maximal, then $a,b,u,v$ are all equal to $e$, so the property is true. Now consider $a,b,u,v,x \in U$, and assume that the property holds for any larger value of $|x|$. If $a,b,u,v \in R_U$, then the property holds by assumption. Without loss of generality, assume that $a \in U \bs R_U$, and write $a=r\cdot a'$, for some $r \in R_U$ and $a' \in U \bs \{e\}$ so that $|a'|=|a|-1$. So $a'\cdot x\cdot u,a'\cdot x\cdot v,b\cdot x\cdot u,b\cdot x\cdot v \in U$. Since $|a'|<|a|$, we deduce by induction that either $a',b$ have a right join or $u,v$ have a left join, and in the latter case we have the desired conclusion. Let us then assume that $a',b$ have a right join $c \in U$. Let us write $c=c'\cdot a'$, where $c' \in U$.

Since $a'\cdot xu,b\cdot xu \in U$ and $a',b$ have a right join $c \in U$, according to Property 5, we deduce that $c\cdot xu \in U$, and similarly $c\cdot xv \in U$. We now consider the four elements $$c'\cdot a'x\cdot u=cxu,c'\cdot a'x\cdot v=cxv,r\cdot a'x\cdot u,r\cdot a'x\cdot v$$ in $U$. Since $a'\cdot x\in U$ and $|a'x| > |x|$, we deduce by induction that either $c',r$ have a right join or $u,v$ have a left join, and in the latter case we have the desired conclusion. Let us then assume that $c',r$ have a right join $d \in U$. Since $c'\cdot a'=c,r\cdot a'=a \in U$ and $c',r$ have a right join $d \in U$, according to Property 5, we deduce that $d\cdot a' \in U$. Now remark that $a=ra' \leq_R da'$ and $b \leq_R c=c'a' \leq_R da'$, so $a,b$ have a common right upper bound for $\leq_R$. Since $(U,\leq_R)$ is a meet-semilattice, we conclude that $a,b$ have a right join (which is the meet of all common right upper bounds of $\{a,b\}$).
\ep

Let us denote by $K_E$ the interval complex of the poset $E$ as in Definition~\ref{def:interval_group} (the interval labeling on $E$ is defined after Lemma~\ref{lem:E partial order}), and let $G_E$ denote the corresponding interval group. We will find a simple criterion ensuring that the interval group $G_E$ is isomorphic to $A \times \Z$, where $A$ is the Artin group associated to $W$.

\bthm \label{thm:homotopy_type_salvetti}
\label{thm:homotopy}
Assume that, for each spherical $T \subset S$, there is a choice of Coxeter element $w_T \in W_T$ such that, for every spherical $T' \subset T$, we know that $w_{T'}\le_L w_T$ in $(W_S,\le_L)$. Assume that
$$U = \bigcup_{T \subset S \mbox{ spherical}} [e,w_T].$$
Here $[e,w_T]$ denotes the interval between $e$ and $w_T$ in $(W_S,\le_L)$.
Then $K_E$ has the homotopy type of the Salvetti complex of the Artin group $A \times \Z$. In particular, the interval group $G_E$ is isomorphic to $A \times \Z$.
\ethm

\bp
For each spherical $T \subset S$, let us denote $U_T=[e,w_T] \subset U$. Consider the subposet $E_T=\left(U_T \times \{0\}\right) \sqcup \left(\ov{U_T} \times \{1\}\right) \subset E$, and denote by $K_{E_T} \subset K_E$ the subcomplex corresponding to the quotient of the geometric realization of $E_T$.

We claim that $K_{E_T}$ has the homotopy type of the Salvetti complex $X_T$ of the Artin group $A_T \times \Z$. Indeed, let us denote by $s_0 \in A_T \times \Z$ a generator of $\Z$, so that the Artin group $A_T \times \Z$ has standard generating system $T'=T \cup \{s_0\}$. Now $w'_T=w_Ts_0$ is a Coxeter element for the spherical Coxeter group $W'_T=W_T \times \Z/2\Z$, and $K_{E_T}$ coincides with the dual Salvetti complex for $w'_T$ as described in~\cite[Section~5]{paolini_salvetti}. According to~\cite[Remark~5.4]{paolini_salvetti}, we deduce that $K_{E_T}$ has the same homotopy type as the standard Salvetti complex $X_T$ for the spherical Artin group $A_T \times \Z$.

By assumption on $U$, it is clear that $K$ is equal to the union of all $K_{E_T}$, for $T \subset S$ spherical. Also remark that the standard Salvetti complex $X$ for the Artin group $A \times \Z$ is equal to the union of all $X_T$, for $T \subset S$ spherical.

According to the proof of~\cite[Theorem~5.5]{paolini_salvetti}, we deduce that $K_E$ has the homotopy type of $X$.
In particular, the interval group $G_E$ of $E$, which is the fundamental group of $K_E$, is naturally isomorphic to the Artin group $A \times \Z$. 
\ep

\bcor \label{cor:Artin_nice_U_garside}
Assume that $W$ is a Coxeter group, with a subset $U \subset W$ satisfying the conditions of Proposition~\ref{pro:E_lattice} and of Theorem~\ref{thm:homotopy_type_salvetti}. Let $A$ be the Artin group associated with $W$.
Then $A \times \Z$ is Garside, with Garside element $(e,1)$. Moreover, the $K(\pi,1)$ conjecture holds for $A$.
\ecor

\bp
This follows from a combination of Proposition~\ref{pro:E_lattice}, Proposition~\ref{pro:I_lattice}, Theorem~\ref{thm:I_lattice_implies_Garside} and Theorem~\ref{thm:homotopy}. Moreover, the standard Salvetti complex $X$ of $A \times \Z$ is aspherical, so in particular the standard Salvetti complex of $A$ itself is aspherical: we deduce that the $K(\pi,1)$ conjecture holds for $A$.
\ep

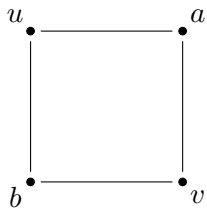
\begin{figure}[H]
\centering
\begin{tikzpicture}
\def \p {0.05}
\def \r {1}
\def \op {1}
\def \gris {black!10}

\draw[fill] (1,1) circle (\p) node(s1) {};
\draw[fill] (-1,1) circle (\p) node(s2) {};
\draw[fill] (-1,-1) circle (\p) node(s3) {};
\draw[fill] (1,-1) circle (\p) node(s4) {};

\draw [-] (s1) edge (s2) (s2) edge (s3) (s3) edge (s4) (s4) edge (s1);

\node (label) at (1.2,1.2) {$a$};
\node (label) at (-1.2,1.2) {$u$};
\node (label) at (-1.2,-1.2) {$b$};
\node (label) at (1.2,-1.2) {$v$};
\end{tikzpicture}
\caption{The right-angled Artin group $A$ over a square.}
\label{fig:raag_square}
\end{figure}

\brk
There are some Artin groups for which it is not possible to find a subset $U \subset A$ satisfying the conditions of Proposition~\ref{pro:E_lattice} and of Theorem~\ref{thm:homotopy_type_salvetti}. Here are two simple examples.
\ben
\item Consider the right-angled Artin group $A \simeq \F_2 \times \F_2$ with defining graph a square with vertices $a,u,b,v$ in this cyclic order (see Figure~\ref{fig:raag_square}), and assume that the conditions of Theorem~\ref{thm:homotopy_type_salvetti} hold. Then we have $au,av,bu,bv \in U$, but neither $a,b$ nor $u,v$ have a join for $\leq_L$. Then the conditions of Proposition~\ref{pro:E_lattice} do not hold.

\begin{figure}[H]
\centering
\begin{tikzpicture}
\def \p {0.05}
\def \r {1}
\def \op {1}
\def \gris {black!10}

\draw[fill] (-3,0) circle (\p) node(s1) {};
\draw[fill] (-2,0) circle (\p) node(s2) {};
\draw[fill] (-1,0) circle (\p) node(s3) {};
\draw[fill] (0,0) circle (\p) node(s4) {};
\draw[fill] (1,0) circle (\p) node(s5) {};
\draw[fill] (2,0) circle (\p) node(s6) {};
\draw[fill] (3,0) circle (\p) node(s7) {};

\draw [-] (s1) edge (s7);

\node (label) at (-2.5,0.2) {$4$};
\node (label) at (-1.5,0.2) {$4$};
\node (label) at (-0.5,0.2) {$4$};
\node (label) at (0.5,0.2) {$4$};
\node (label) at (1.5,0.2) {$4$};
\node (label) at (2.5,0.2) {$4$};

\node (label) at (-3,-0.3) {$s_1$};
\node (label) at (-2,-0.3) {$s_2$};
\node (label) at (-1,-0.3) {$s_3$};
\node (label) at (0,-0.3) {$s_4$};
\node (label) at (1,-0.3) {$s_5$};
\node (label) at (2,-0.3) {$s_6$};
\node (label) at (3,-0.3) {$s_7$};

\node (label) at (-2.5,-0.8) {$a$};
\node (label) at (-1.5,-0.8) {$b$};
\node (label) at (1.5,-0.8) {$u$};
\node (label) at (2.5,-0.8) {$v$};
\end{tikzpicture}
\caption{The Dynkin diagram of an Artin group for which Corollary~\ref{cor:Artin_nice_U_garside} does not apply.}
\label{fig:problem_artin}
\end{figure}
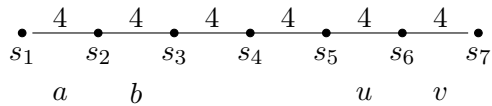

\item Consider the $A$ with defining graph a complete graph over $7$ vertices, whose Dynkin diagram is a line with vertices $s_1,s_2,\dots,s_7$, with all edge labels equal to $4$ (see Figure~\ref{fig:problem_artin}). Assume that the conditions of Theorem~\ref{thm:homotopy_type_salvetti} hold. Consider the four elements of $U$: $a=s_1s_2$ or $s_2s_1$ (depending on the ordering on $S$), $b=s_2s_3$ or $s_3s_2$, $u=s_5s_6$ or $s_6s_5$ and $v=s_6s_7$ or $s_7s_6$. Then $au,av,bu,bv \in U$, but neither $a,b$ nor $u,v$ have a join for $\leq_L$. Then the conditions of Proposition~\ref{pro:E_lattice} do not hold.

\een
\erk

\bcor \label{cor:complex}
Assume that, for each spherical $T \subset S$, there is a choice of Coxeter element $w_T \in W_T$ such that, for every spherical $T' \subset T$, we have $w_{T'} \leq_L w_T$. Assume that
$$U = \bigcup_{T \subset S \mbox{ spherical}} [e,w_T].$$
Let $\widehat U$ be the lift of $U$ from $W_S$ to $A_S$ via the (compatible) isomorphism between the dual Artin group associated with $A_T$ for each $T\subset S$ spherical and $A_T$ (cf. Theorem~\ref{thm:bessis}, and more precisely \cite[Theorem 2.2.5]{bessis}).

Assume that $U$ satisfying the conditions of Proposition~\ref{pro:E_lattice}. Let $X_S$ be the flag complex of the Cayley graph of $A_S$ with generating set $\widehat U$. Then $X_S$ admits an $A_S$-equivariant CUB metric such that each simplex of $X_S$ is equipped with a polyhedral norm as in  \cite{haettel_simplicial_npc}.
\ecor

\bp
By Corollary~\ref{cor:Artin_nice_U_garside}, $A_S\times \Z$ is a Garside group with the choice of fundamental interval $E$ as in Section~\ref{sec:general_construction}. Note that Bestvina complex (cf. Section~\ref{subsec:garside_complex}) for the Garside group $A_S\times \Z$ is isomorphic to flag complex of the Cayley graph of $A_S$ with generating set $\widehat U$. Thus we are done by Theorem~\ref{thm:garside_cub}.
\ep

\section{Cyclic-type Artin groups}
\label{sec:cyclic_type}
In this section and the next, we will describe families of Artin groups for which we can find a set $U$ satisfying the conditions of Proposition~\ref{pro:E_lattice} and of Theorem~\ref{thm:homotopy_type_salvetti}. This section will focus on cyclic type Artin groups.
The main goal of this section will be Proposition~\ref{prop_cyclic}.

\subsection{Spherical Artin group with linear Dynkin diagram}
The goal of this subsection is to prove Corollary~\ref{cor:left_right}, which is a main ingredient towards Proposition~\ref{prop_cyclic}.  We will first establish a special case of Corollary~\ref{cor:left_right}, which is Lemma~\ref{lem:left_right}, then use Lemma~\ref{lem:left_right} to deduce the more general Corollary~\ref{cor:left_right}.

\blem
\label{lem:trivial_word}
Let $W_S$ be an arbitrary Coxeter group (not necessarily spherical). Let $w$ denote a word in $S$ representing the trivial element of $W_S$. Then each letter of $w$ appears at least twice.
\elem

\bp
By contradiction, assume that we can write $w=usv$, where $s \in S$ and $u,v$ are words in $S \bs \{s\}$. Then, in the Coxeter groups $W_S$, the words $s$ and $u^{-1}v^{-1}$ represent the same element. Since $s$ in the support of $s$ and not of $u^{-1}v^{-1}$, this is a contradiction.
\ep

\blem
\label{lem:commutation_letter}
Let $W_S$ be an arbitrary Coxeter group (not necessarily spherical). Let $s \in S$, and let $w$ denote a reduced word in $S \bs \{s\}$ representing an element commuting with $s$. Then every letter of $w$ commutes with $s$.
\elem

\bp
Since $w$ and $s$ have disjoint supports, the words $sw$ and $ws$ are reduced. We can pass from the reduced word $sw$ to the reduced word $ws$ by applying only standard relations (see for instance~\cite[Theorem~3.4.2]{davis_coxeter}). This implies that $s$ commutes with every letter of $w$.
\ep

\blem
\label{lem:reducible_word}
Let $W_S$ be an arbitrary Coxeter group (not necessarily spherical). Let $\{s_1,s_2,\ldots, s_n\}\subset S$ such that, for each $1 \leq i \leq n-1$, there exists $i < j \leq n$ such that $s_i$ and $s_j$ do not commute. Then the word $$s_1\cdots s_{n-1}s_ns_{n-1}\cdots s_1$$
is reduced.
\elem

\bp
We induct on $n$. Then case $n=1$ is trivial. For the general case, by contradiction, assume that the word $w=s_1\cdots s_{n-1}s_ns_{n-1}\cdots s_1$ is not reduced. According to the deletion condition (see for instance~\cite[Corollary~5.8]{humphreys1990reflection}), $w$ can also be represented by a word $w'$ obtained from $w$ by deleting two letters.

\mk

Since $w$ represents a reflection of $W_S$, $w'$ also represents a reflection. According to the strong exchange condition (see~\cite[Theorem~5.8]{humphreys1990reflection}), if we remove one letter from $w'$ we may obtain the trivial element. According to Lemma~\ref{lem:trivial_word}, there are only two possibilities:
\begin{enumerate}
	\item there exists $1 \leq i \leq n-1$ such that $w'$ is obtained from $w$ by removing the two occurences of $s_i$;
	\item there exists $1\leq i\leq n-1$ such that $w'$ is obtained from $w$ by removing $s_n$ and one of the $s_i$.
\end{enumerate}

\mk

First we rule out possibility (1). In this case $$w'=s_1s_2 \dots s_{i-1}s_{i+1} \dots s_n \dots s_{i+1} s_{i-1} \dots s_1.$$ By conjugating by $s_1s_2 \dots s_{i-1}$, we deduce that the words $$s_i\cdots s_{n-1}s_ns_{n-1}\cdots s_i=s_{i+1}\cdots s_{n-1}s_ns_{n-1}\cdots s_{i+1}.$$  In particular, the element $u=s_{i+1}\cdots s_{n-1}s_ns_{n-1}\cdots s_{i+1}$ commutes with $s_i$. As $u=s_{i+1}\cdots s_{n-1}s_ns_{n-1}\cdots s_{i+1}$ is reduced by induction assumption, according to Lemma~\ref{lem:commutation_letter}, we deduce that $s_i$ commutes with every letter $s_{i+1},\dots,s_n$, which contradicts the assumption.

\mk

Now we rule out possibility (2). In this case, $w'=s_1\ldots s_{i-1}s_is_{i-1}\ldots s_1$. By conjugation by $s_1s_2\ldots s_{i-1}$, we deduce that $s_i=s_i\ldots s_{n-1}s_ns_{n-1}\ldots s_i$. Thus $s_{i+1}\ldots s_ns_{n-1}\ldots s_{i}$ represents the trivial element, contradicting Lemma~\ref{lem:trivial_word}.
\ep

\blem \label{lem:left_right}
Suppose $S$ is spherical with linear Dynkin diagram. We label elements of $S$ as $\{s_1,\ldots,s_n\}$ using one of the two linear orders induced by the Dynkin diagram. 
 Let $U$ be the dual Garside interval with respect to the dual Garside element $\delta=s_1s_2\cdots s_n$ and let $R_U$ be the set of reflections in $U$.
Assume that $u,v \in R_U$ are such that $u\cdot v \in U$. Let $I=\supp(u)$ and $J=\supp(v)$. Define $\min (I)$ be the smallest element in $I$ with respect to the chosen linear order.

If $\min(I)-1\in J$, then $I\subset J$.
\elem

Here we slightly abuse notation, and will use $s_i-1$ to denote $s_{i-1}$, and $s_i+1$ to denote $s_{i+1}$.

\bp
Suppose we fix an order of $\{s_1,s_2,\ldots,s_n\}$ from left to right of the Dynkin diagram. Then we have the following cases to consider.
\ben
\item Type $I_m$ with $m\ge 3$, and $n=2$.
\item Type $A_n$.
\item Type $B_n$, with Dynkin diagram $3-3-\cdots-3-4$.
\item Type $B_n$, with Dynkin diagram $4-3-\cdots-3-3$.
\item Type $H_3$, with Dynkin diagram $3-5$, and $n=3$.
\item Type $H_3$, with Dynkin diagram $5-3$, and $n=3$.
\item Type $H_4$, with Dynkin diagram $3-3-5$, and $n=4$.
\item Type $H_4$, with Dynkin diagram $5-3-3$, and $n=4$.
\item Type $F_4$, with Dynkin diagram $3-4-3$, and $n=4$.
\een

A priori, two different linear orders on the Dynkin diagrams gives two different cases. However, we claim that if the lemma holds for a given linear order on the Dynkin diagram, then it also holds with respect to the reserved linear order. 
Indeed, suppose $\{s_1,\ldots,s_n\}$ are named according a given linear order on the Dynkin diagram. Let $u,v,I,J$ and $\delta$ be as in the lemma. For each element $g\in W_S$ with its reduced expression $g=s_{i_1}s_{i_2}\cdots s_{i_k}$, we define $\bar g=s_{i_k}\cdots s_{i_2}s_{i_1}$. Note that $\bar g\in W_S$ does not depend on the choice of reduced expression of $g$. Then $uv\leq_L \delta$ if and only if $\bar v\bar u\leq_R \bar\delta$ if and only if $\bar v\bar u\leq_L \bar \delta$ (by Theorem~\ref{thm:bessis} and Definition~\ref{def:garside} (4)). As $\min(\supp(u))-1\in \supp(v)$ with respect to the given order, and $I$ and $J$ are irreducible (Lemma~\ref{lem:irreducible}), we deduce that either $\min(\supp(\bar v))-1\in \supp(\bar u)$ with respect to the reserved order, or $I\subset J$. In the latter case we are already done. In the former case, if we know the lemma holds with respect to the original order, then we deduce that $\supp(\bar v)\subset\supp(\bar u)$, hence $J\subset I$, which contradicts that $\min(I)-1\in J$ with respect to the original order, so the former case cannot happen. This proves the claim.

\mk

Assume first that $|S|=2$, i.e. $W_S$ is of type $I_m$ with $m\ge 3$. Assume that $I,J \neq S=\{s_1,s_2\}$, so that $u,v \in \{s_1,s_2\}$. We only have to prove that $s_2s_1 \not\leq_L \delta$.

By contradiction, assume that $s_2s_1 \leq_L \delta$. Since $\delta=s_1s_2$ has reflection length $2$, we deduce that $s_2s_1=\delta$, so $s_1s_2=s_2s_1$. This contradicts $m \geq 3$.

\mk

Suppose $|S|>2$. We argue by contradiction and assume $I$ is not a subset of $J$. As $J$ contains $\min(I)-1$, and $S=I\cup J$, by Lemma~\ref{lem:irreducible}, we assume that $I=\{s_i,\dots,s_{n}\}$ and $J=\{s_1,\dots,s_j\}$. Note that $j \leq n-1$, otherwise we already have $I\subset J$. 

%In the following we will show that $i \geq j+2$, which will contradicts the assumption that $s_{\min(I)-1}\in J$.

\mk

%Note that for each element $g\in W_S$ with its reduced expression $g=s_{i_1}s_{i_2}\cdots s_{i_k}$, we define $\bar g=s_{i_k}\cdots s_{i_2}s_{i_1}$. Note that $\bar g\in W_S$ does not depend on the choice of reduced expression of $g$. Then $uv\leq_L \delta$ if and only if $\bar v\bar u\leq_R \bar\delta$ if and only if $\bar v\bar u\leq_L \bar \delta$. This allows us to reduce the $5-3$ case to the $3-5$ case, the $5-3-3$ case to the $3-3-5$ case, and the $4-3-\cdots-3-3$ case to the $3-3-\cdots-3-4$ case. 

Assume that $W_S$ is not of type $F_4$ with $i=2$ and $j=3$. Then either $J$ is of type $A_j$ or $I$ is of type $A_{n-i+1}$. By the order-reversing argument in the above claim, we can assume without loss of generality that $J$ is of type $A_j$.
Then we know that $v=s_1s_2\cdots s_{j-1}s_js_{j-1}\cdots s_2s_1$. Indeed, this is the unique reflection in $A_j$ with support being $J$, and we can see $v$ must be of this form by treating $A_j$ as the permutation group of $j+1$ elements $\{\alpha_1,\ldots,\alpha_{j+1}\}$, with $s_i$ being the transposition permuting $\alpha_i$ and $\alpha_{i+1}$. Then $v$ is the transposition permuting $\alpha_1$ and $\alpha_{j+1}$. 
We also have
\begin{align*}
	&\delta = s_1\cdots s_js_{j+1} \dots s_{n} = vs_1\cdots s_{j-1}s_{j+1} \dots s_{n} \\
	&= (vs_1v^{-1})\cdots (vs_{j-1}v^{-1})(vs_{j+1} \dots s_{n}v^{-1})v.
\end{align*}
By the Garside property, there exists $w \in [1,\delta]$ with reflection length $n-1$ such that $w\cdot v=\delta$. Thus 
$$
(vs_1v^{-1})\cdots (vs_{j-1}v^{-1})(vs_{j+1}v^{-1}) \cdots (vs_{n}v^{-1})
$$
is a minimal reflection decomposition of $w$. By \cite[Lemma 3.7]{mccammond_failure_lattice_property}, $vs_kv^{-1}\leq_L w$ for $k\neq j$. Since $u\cdot v \leq_L \delta$, we know that $u\cdot v\leq_R \delta$ by Theorem~\ref{thm:bessis}. Thus $u \leq_R w$.

According to \cite[Lemma 1.4.3]{bessis}, the element $w$ is a Garside element for the parabolic subgroup $P_w$ of $P$ generated by the reflections which are $\leq_L$-smaller than $w$. Since $w$ has reflection length $n-1$, this subgroup $P_w$ equals
$$P_w=\<vs_1v^{-1},\cdots vs_{j-1}v^{-1},vs_{j+1}v^{-1},\dots,vs_{n}v^{-1}\> =v\left(\<s_1,\dots,s_{j-1}\> \times \<s_{j+1},\dots,s_{n}\>\right)v^{-1}.$$
As $u \leq_R w$, we know that $u \in P_w$ by \cite[Theorem 1.4]{baumeister2014note}. Hence $v^{-1}uv\in W_{S\setminus\{s_j\}}$. By Lemma~\ref{lem:irreducible},  $\supp(v^{-1}uv)\subset \{s_1,\ldots,s_{j-1}\}$ or $\supp(v^{-1}uv)\subset \{s_{j+1},\ldots,s_n\}$. Thus $u\in v\left(\<s_1,\ldots,s_{j-1}\>\right)v^{-1}$ or $u\in  v\left(\<s_{j+1},\ldots,s_n\>\right)v^{-1}$. Also since $$v\left(\<s_1,\ldots,s_{j-1}\>\right)v^{-1} \subset W_{\{s_1,\ldots,s_j\}},$$ we rule out that $u\in v\left(\<s_1,\ldots,s_{j-1}\>\right)v^{-1}$, hence $u\in  v\left(\<s_{j+1},\ldots,s_n\>\right)v^{-1}$. In particular,
$$
u\in  v\left(\<s_{j+1},\ldots,s_n\>\right)v^{-1}\cap W_I=\<vs_{j+1}v^{-1},s_{j+2},\ldots,s_n\>\cap \<s_i,\ldots,s_n\>.
$$

As $s_{\min(I)-1}\in J$, we know $i \leq j+1$. Let $P=\<vs_{j+1}v^{-1},s_{j+2},\ldots,s_n\>\cap \<s_i,\ldots,s_n\>$. By \cite{qi2007note}, $P$ is a parabolic subgroup of $W_S$. Note that $P\supset W_{\{s_{j+2},\ldots,s_n\}}$. On the other hand, $vs_{j+1}v^{-1}\notin W_{\{s_i,\ldots,s_n\}}$ as $vs_{j+1}v^{-1}=s_1s_2\cdots s_{j}s_{j+1}s_j\cdots s_1$ and the word $s_1s_2\cdots s_{j}s_{j+1}s_j\cdots s_1$ is reduced by Lemma~\ref{lem:reducible_word}. Hence $u\in P=W_{\{s_{j+2},\ldots,s_n\}}$. As $j+2\ge i+1$, this contradicts that $\supp(u)=\{s_i,\ldots,s_n\}$. We deduce this contradiction as we were assuming $I$ is not a subset of $J$ in the beginning. This proves that $I\subset J$.

%It remains to prove the claim. Note that $vs_{j+1}v^{-1}=s_1s_2\cdots s_{j}s_{j+1}s_j\cdots s_1$. Note that this word is already in reduced form with respect to $S$, indeed, if $m_{s_j,s_{j+1}}=4$, then there is no way to transfer to another word by applying a relation, thus the word is reduced by Tits's solution to the word problem of Coxeter groups; if $m_{s_j,s_{j+1}}=3$, then $W_{s_1,\ldots,s_{j+1}}$ is the symmetry group on $j+2$ elements, and the word is a reduced word in this symmetry group. As this word is reduced and contains $s_1$, it is not contained in $W_{\{s_i,\ldots,s_n\}}$ as $i>1$. Then the claim is proved.

\mk

The remaining case is in type $F_4$ with $i=2$ and $j=3$. Then $v$ is a reflection inside the Coxeter group of $B_3$ generated by $s_1,s_2$ and $s_3$, which has Dynkin diagram $3-4$. Consider the canonical representation of Coxeter group of type $B_3$ acting on $\R^3$. Then $s_1$ acts by the orthogonal reflection along $x_1=x_2$, $s_2$ acts by the orthogonal reflection along $x_2=x_3$, and $s_3$ acts by the orthogonal reflection along $x_3=0$. Note that there are nine reflection in $W_{s_1,s_2,s_3}$, whose reflection hyperplanes are $x_i=\pm x_j$ for $1\le i\neq j\le 3$ and $x_i=0$ for $1\le i\le 3$. Note that reflections along $x_i=x_j$ for $1\le i\neq j\le 3$ are supported on $W_{s_1,s_2}$; reflections along $x_2=\pm x_3$ or $x_i=0$ for $i=2,3$ are supported on $W_{s_2,s_3}$. This gives 6 reflections in total. The remaining three reflections in $W_{s_1,s_2,s_3}$ give all the possibilities of $v$. More precisely, reflection along $x_1=0$ gives $v=s_1s_2s_3s_2s_1$, reflection along $x_1+x_3=0$ gives $v=s_1s_3s_2s_3s_1$, and reflection along $x_1+x_2=0$ gives $v=s_2s_3s_2s_1s_2s_3s_2$.

The case $v=s_1s_2s_3s_2s_1$ is identical to before. Now we assume $v=s_1s_3s_2s_3s_1=s_3s_1s_2s_1s_3=s_3s_2s_1s_2s_3$. Then
\begin{align*}
	&\delta = s_1s_2s_3s_4= s_1s_2s_3s_4s_3s_2s_1s_2s_3v \\
	&= (s_1s_2s_3s_4s_3s_2s_1)(s_2)(s_3)v.
\end{align*}
By the same argument as before, we know $u\in \<s_1s_2s_3s_4s_3s_2s_1,s_2,s_3\>$. Let
\begin{align*}
	&P=\<s_1s_2s_3s_4s_3s_2s_1,s_2,s_3\>=s_1s_2s_3\<s_4,s_1,s_2s_3s_2\>s_3s_2s_1=s_1s_2s_3s_2\<s_4,s_2s_1s_2,s_3\>s_2s_3s_2s_1\\
	&=s_1s_2s_3s_2\<s_4,s_1s_2s_1,s_3\>s_2s_3s_2s_1=s_1s_2s_3s_2s_1\<s_4,s_2,s_3\>s_1s_2s_3s_2s_1.
\end{align*} 
In particular, $P$ is a parabolic subgroup. Note that $u\in P\cap W_{s_2,s_3,s_4}$. 

By \cite{qi2007note}, $P\cap W_{s_2,s_3,s_4}$ is a parabolic subgroup of $W_S$. Note that $\<s_2,s_3\>\subset P\cap W_{s_2,s_3,s_4}$. Moreover, $s_1s_2s_3s_4s_3s_2s_1\in P\setminus W_{s_2,s_3,s_4}$ as $s_1s_2s_3s_4s_3s_2s_1$ is a reduced word by Lemma~\ref{lem:reducible_word}. Thus $P\cap W_{s_2,s_3,s_4}=\<s_2,s_3\>$. Thus $s_4\notin \supp(u)$, which is a contradiction.

It remains to look at the case $v=s_2s_3s_2s_1s_2s_3s_2$. Then 
\begin{align*}
	&\delta = s_1s_2s_3s_4= s_1s_2s_3s_4s_2s_3s_2s_1s_2s_3s_2v \\
	&= (s_1s_2s_3s_4s_3s_2s_1)(s_1s_2s_3s_2s_3s_2s_1)(s_2s_3s_2)v.
\end{align*}
Note that $s_2s_3s_2=s_1s_2s_3(s_3s_2s_1s_2s_3s_2s_1s_2s_3)s_3s_2s_1=s_1s_2s_3(s_1s_2s_3s_2s_1)s_3s_2s_1$. Thus by repeating the previous discussion, we know 
$$
u\in s_1s_2s_3\<s_4,s_2,s_1s_2s_3s_2s_1\>s_3s_2s_1=s_1s_2s_3s_1s_2\<s_4,s_1,s_3\>s_2s_1s_3s_2s_1.
$$
As $u$ belongs to a parabolic subgroup which splits as a product and $s_4\in \supp(u)$, we argue as before to deduce that $$u\in s_1s_2s_3s_1s_2\<s_4,s_3\>s_2s_1s_3s_2s_1:=P.$$
Thus $u\in P\cap W_{s_2,s_3,s_4}$. Note that 
\begin{align*}
	&	s_1s_2s_3s_1s_2 (s_3)s_2s_1s_3s_2s_1=s_1s_2s_1(s_3s_2s_3s_2s_3)s_1s_2s_1 \\
	&=s_2s_1s_2 (s_2s_3s_2)s_2s_1s_2=s_2s_1s_3s_1s_2=s_2s_3s_2\in P\cap W_{s_2,s_3,s_4}.
\end{align*}
Thus $P\cap W_{s_2,s_3,s_4}$ is a parabolic subgroup of rank $\ge 1$. On the other hand,
$$
s_1s_2s_3s_1s_2 (s_4)s_2s_1s_3s_2s_1=s_1s_2s_3s_4s_3s_2s_1.
$$
As $s_1s_2s_3s_4s_3s_2s_1$ is a reduced word in $W_S$ by Lemma~\ref{lem:reducible_word}, it can not be contained in $W_{s_2,s_3,s_4}$. Thus $P\cap W_{s_2,s_3,s_4}\subsetneq P$. It follows that $P\cap W_{s_2,s_3,s_4}=\<s_2s_3s_2\>$, hence $s_4\notin \supp(u)$, which is a contradiction.
\ep

\bcor \label{cor:left_right}
Suppose $S$ is spherical with linear Dynkin diagram. We label elements of $S$ as $\{s_1,\ldots,s_n\}$ using one of the two linear orders induced by the Dynkin diagram. 
Let $U$ be the dual Garside interval with respect to the dual Garside element $\delta=s_1s_2\cdots s_n$ and let $R_U$ be the set of reflections in $U$.
Assume that $u,v \in U$ are such that $\supp(u)$ and $\supp(v)$ are irreducible, and $u\cdot v \in U$. Let $I=\supp(u)$ and $J=\supp(v)$. 

If $\min(I)-1\in J$, then $I\subset J$.
\ecor

Recall that we slightly abuse notation, use $s_i-1$ to denote $s_{i-1}$, and $s_i+1$ to denote $s_{i+1}$.
\bp
	Let $v=r_1r_2\cdots r_k$ with $r_i\in R_U$ be a minimal reflection decomposition of $v$. As $\min(I)-1\in J$, there exists $1\le i\le k$ such that $\min(I)-1\in \supp(r_i)$. Suppose $I\subset J$ is not true. Then $\max(\supp(r_i))<\max(I)$, where $\max(I)$ is defined to be the biggest element in $I$. Hence $\max(\supp(r_i))+1\in I$. We write a minimal reflection decomposition of $u$ as $u=t_1t_2\cdots t_m$. Then there exists $1\le j\le m$ such that $\max(\supp(r_i))+1\in\supp(t_j)$. As $u\cdot v\in \delta$, we know $t_1\cdots t_mr_1\cdots r_k$ is a minimal reflection decomposition of $uv$. In particular $t_1\cdots t_mr_1\cdots r_k\le_L \delta$. By \cite[Lemma 3.7]{mccammond_failure_lattice_property}, $t_j\cdot r_i\le_L\delta$. By construction we have $\min(\supp(t_j))-1\in \supp(r_i)$, and $\supp(r_i)$ does not contain $\supp(t_j)$ (as $\max(\supp(r_i))+1\in\supp(t_j)$), which contradicts Lemma~\ref{lem:left_right}. Thus the corollary is proved.
\ep

\subsection{Cyclic-type Artin groups}
Let $W_S$ be a cyclic type Coxeter group (cf.  Table~\ref{tab:cyclic}).
We take a cyclic order on $S$ coming from its Dynkin diagram, and denote elements in $S$ as elements in $\Z/n\Z$.
For each $i \in \Z/n\Z$, consider the dual Garside interval $U_i$ in $P_{S \bs i}$ with respect to the dual Garside element $\delta_i=s_{i+1}s_{i+2} \dots s_n s_1 \dots s_{i-1}$. Let $U=\cup_{i \in \Z/n\Z} U_i$. It is clear that this set $U$ satisfies the assumptions of Theorem~\ref{thm:homotopy_type_salvetti}.

\mk

For each $i \in \Z/n\Z$, consider the set $R_{U_i} \subset U_i$ of reflections in the spherical parabolic subgroup $P_{S \bs i}$, and let $R_U=\cup_{i \in \Z/n\Z} R_{U_i} \subset U$.

\bpro
\label{prop_cyclic}
The sets $R$ and $U$ satisfy all assumptions from Proposition~\ref{pro:E_lattice}. In particular, if $A_\Gamma$ is of cyclic type, then $A_\Gamma\times \Z$ is a Garside group.
\epro

\bp
We verify each assumption of Proposition~\ref{pro:E_lattice} as follows.
\ben
\item Any $u \in U_i$ can be written as a product of elements in $R_{U_i}$ which is a minimal length reflection factorization.
\item Let $r_1,\dots,r_m \in U_R$ such that $u =r_1\cdot r_2 \cdot \ldots \cdot r_m \in U$. Let $i \in \Z/n\Z$ such that $u \in U_i$. Then $r_1,\dots,r_m \in P_{S \bs i}$ by \cite[Theorem 1.4]{baumeister2014note},  and $r_1\cdot r_2\cdot \ldots \cdot r_m \leq_L \delta_i$. So both $r_1 \cdot r_2\cdot \ldots\cdot r_{m-1}$ and $r_2 \cdot r_3\cdot\ldots \cdot r_m$ belong to $P_{S \bs i}$, and also they are both prefixes of $\delta_i$. Hence $r_1 \cdot r_2\cdot \ldots\cdot r_{m-1}$ and $r_2 \cdot r_3\cdot\ldots \cdot r_m$ belong to $U_i$.
\item Let $r,r' \in U_R$ which admit a common left upper bound $u \in U$. Let $i \in \Z/n\Z$ such that $u \in U_i$. By \cite[Theorem 1.4]{baumeister2014note}, $r,r'\in U_i$ and $u$ is a common left upper bound for $r,r'$ in $(U_i,\le_L)$. In particular, $r,r' \leq_L \delta_i$, so $r$ and $r'$ admit a unique left join $u_i$ in $(U_i,\le_L)$. Now we show $u_i$ is also the join of $r$ and $r'$ in $(U,\le_L)$. Indeed, take an arbitrary left upper bound $u'$ of $r,r'$ in $U$. Suppose $u'\in U_j$. Then as before we know $r,r'\in U_j$ and $u'$ is a common left upper bound of $r,r'$ in $(U_j,\le_L)$ by Theorem~\ref{thm:bessis}. Let $u_j$ (resp. $u_{ij}$) be the left join of $r,r'$ in $(U_j,\le_L)$ (resp. $(U_i\cap U_j,\le_L)$). One readily verifies that $u_{ij}=u_i$, $u_{ij}=u_j$ and $u_j\le_L u'$. Thus $u_i\le_L u'$, implying $u_i$ is the left join of $r$ and $r'$ in $U$. The case of common right upper bound is similar.
\item Let $a,u,v,w$ be as in Proposition~\ref{pro:E_lattice} (4).
Let $I,J,K \subset \Z/n\Z$ denote $\supp(a)$, $\supp(u)$ and $\supp(v)$ respectively. We first prove $I\cup J\cup K\subsetneq \Z/n\Z$. Suppose by contradiction that $I\cup J\cup K=\Z/n\Z$. Let $\{I_i\}_{i=1}^k$ be the irreducible components of $I$. By Lemma~\ref{lem:irreducible}, $a=a_1\cdot a_2\cdot\ldots\cdot a_k$ such that $\supp(a_i)=I_i$. By \cite[Lemma 3.7]{mccammond_failure_lattice_property}, $a_i\cdot u\leq_L a\cdot u$, thus $a_i\cdot u\in U$. Similarly, $a_i\cdot v\in U$ for $1\le i\le k$. As $I\cup J\cup K=\Z/n\Z$, for each $a_i$, we know either $\min(I_i)-1\in J$ or $\min(I_i)-1\in K$ (here $I_i$ inherits a linear order from the cyclic order on $\Z/n\Z$, hence it make senses to talk about $\min(I_i)$ and $\min(I_i)-1$). We first look at the case that $\min(I_j)-1\in J$. Consider $I_i\cup J$, which is irreducible. As $a_i\cdot u\in U$, there exists $i_0$ such that $a_i\cdot u\in U_{i_0}$. Then $I_i\cup J\subset S\setminus \{i_0\}$ by \cite[Theorem 1.4]{baumeister2014note}. We endow $S\setminus\{i_0\}$ with the linear order induced from the cyclic order on $S$, then Corollary~\ref{cor:left_right} implies $I_i\subset J$. This shows that $I_i\subset J\cup K$ for each $i$. Thus $I\cup J\cup K=J\cup K$. However, as $u$ and $v$ has a left join $w$ in $U$, there exists $i'_0$ such that $w\in U_{i'_0}$, hence $u,v\in U_{i'_0}$ by  \cite[Theorem 1.4]{baumeister2014note}. Thus $J\cup K\subsetneq S$, which is a contradiction. The case that $\min(I_j)-1\in K$ can be treated similarly.

Let $i\in S$ such that $I\cup J\cup K\subset S\setminus\{i\}$.
As $w$ is a left join of $u$ and $v$ in $U$, by the discussion in item 3, $\supp(w)\in J\cup K$ and $w$ is the left join of $u$ and $v$ in $(U_i,\le_L)$. Thus $\supp(a)\cup \supp(w)\subset U_i$. Clearly $aw\in U_i\subset U$. It remains to show $|aw|=|a|+|w|$.
As $(U_i,\le_L)$ is a lattice, $a\cdot u$ and $a\cdot v$ has a left join in $U_i$, denoted by $a'$. As $a\le_L a'$, we know $a'=a\cdot w'$ for some $w'\in U_i$. By cancellation property in $U_i$, we know $u\le_L w'$ and $v\le_L w'$. Thus $w\le_L w'$. Then $w'=w\cdot w_0$ for $w_0\in U_i$. Then $|a'|=|a|+|w'|=|a|+|w|+|w_0|$. As $a'=aww_0$, $|a'|\le |aw|+|w_0|$. Thus $|aw|=|a|+|w|$.

\item This is similar to the previous item.

\item Let $a,b,u,v,x$ be as in Proposition~\ref{pro:E_lattice} (6). Let $I_a=\supp(a)$. Similarly we define $I_b,I_u$ and $I_v$.
We claim either $I_a\cup I_b\subsetneq S$ or $I_u\cup I_v\subsetneq S$. Assume by contradiction that $I_a\cup I_b=S$ and $I_u\cup I_v=S$. As $a\cdot x\cdot u\in U$, there there exists $i\in S$ such that $a\cdot x\cdot u\le_L \delta_i$. By \cite[Lemma 3.7]{mccammond_failure_lattice_property}, $a\cdot u\le_L a\cdot x\cdot u\le_L\delta_i$, thus $a\cdot u\in U_i\subset U$. Similarly, $a\cdot v,b\cdot u,b\cdot v\in U$. As $I_a\cup I_b=S$, either $\min(I_u)-1\subset I_a$ or $\min(I_u)-1\subset I_b$. If $\min(I_u)-1\subset I_a$, as $a\cdot u\in U_i$, we know from Lemma~\ref{lem:left_right} that $I_u\subset I_a$. As $I_u\cup I_v=S$, $I_a\cup I_v=S$, which contradicts $a\cdot v\in U$. The case of $\min(I_u)-1\subset I_b$ is similar.

If $I_a\cup I_b\subsetneq S$, then there is $i\in S$ such that $a,b\in U_i$. As $(U_i,\le_R)$ is lattice by Theorem~\ref{thm:bessis}, $a$ and $b$ have a right join in $(U_i,\le_R)$, which is also a right join of $a$ and $b$ in $U$ by the argument in item 3. The case $I_u\cup I_v\subsetneq S$ is similar.
\een
\ep

Now the following is a consequence of Proposition~\ref{prop_cyclic}, Corollary~\ref{cor:Artin_nice_U_garside} and the main result of \cite{jankiewicz_schreve_center_kpi1}.

\bcor
\label{cor:K(pi,1)1}
Assume that $A_\Gamma$ is of cyclic type. Then $A_\Gamma$ satisfies the $K(\pi,1)$ conjecture and has trivial center.
\ecor
\section{Combination of cyclic type and spherical type Artin groups}
\label{sec:combination}

In this section, we will describe another family of Artin groups for which we can find a set $U$ satisfying the conditions of Proposition~\ref{pro:E_lattice} and of Theorem~\ref{thm:homotopy_type_salvetti}. The main goal of this section is to prove Theorem~\ref{thm_main}.

An edge of a Coxeter presentation graph is \emph{large} if it its label is $\ge 3$. For each induced subgraph $\Lambda\subset\Gamma$, we define $\Lambda^{\perp}$ to be the induced subgraph spanned by vertices of $\Gamma\setminus \Lambda$ which commute with every vertex in $\Lambda$.

We will be considering orientation of each large edges of $\Gamma$. For the moment suppose $\Gamma$ is spherical and we orient each large edge. We say a Coxeter element of $\Gamma$ is compatible with such orientation if whenever there is an oriented edge from $s_1\in S$ to $s_2\in S$, then $s_1$ appears before $s_2$ in the expression of the Coxeter element.

\begin{lem}
	\label{lem:coxeter}
	Given a spherical Coxeter presentation graph $\Gamma$ with an orientation on each of its large edges, any two Coxeter elements that are compatible with the orientation are equal.
\end{lem} 

\bp
We will prove it by induction on the rank of $\Gamma$.

Let us assume that $w=s_1 \dots s_n$ and $w'=s'_1 \dots s'_n$ are two Coxeter elements that are compatible with the orientation. Let $i \in \{1,\dots,n\}$ such that $s'_i=s_1$. Assume that, among all possible reduced expressions of $w'$ that are compatible with the orientation, the position $i$ of $s_1$ is minimal. We will prove that $i=1$.

Assume by contradiction that $i>1$. Since $i$ is minimal, we deduce that the edge between $s'_{i-1}$ and $s'_i$ has label $\ge 3$. As $w'$ is compatible with the orientation, we deduce that the edge between $s'_{i-1}$ and $s'_i$ is oriented from $s'_{i-1}$ to $s'_i$. As $w$ is also compatible with the orientation and $s'_i=s_1$, we deduce that this edge is oriented from $s'_{i}$ to $s'_{i-1}$. This is a contradiction.

So $s'_1=s_1$. By induction, we deduce that $s_2 \dots s_n=s'_2 \dots s'_n$, hence $w=w'$.
\ep

\begin{lem}
	\label{lem:order}
	Given a spherical Coxeter presentation graph $\Gamma$ with orientation on its large edges. Let $\delta$ be the Coxeter element which is compatible with the orientation. Let $[1,\delta]$ be the collection of elements in $W_\Gamma$ that are prefixes of $\delta$ with respect to the reflection length on $W_\Gamma$.
	
	Give two reflections $r_1,r_2\in [1,\delta]$ such that $r_1r_2\in [1,\delta]$. Suppose there exist $s_1$ and $s_2$ such that $s_1\in \supp(r_1)\setminus \supp(r_2)$ and $s_2\in \supp(r_2)\setminus\supp(r_1)$.
	Then either $s_1$ and $s_2$ commute, or there is an oriented edge from $s_1$ to $s_2$.
\end{lem}
\bp
We argue by contradiction and assume there is an oriented edge from $s_2$ to $s_1$. Let $\Lambda$ be the Dynkin diagram, which is a tree. Then we cut $\Lambda$ along the midpoint of the edge $\overline{s_2s_1}$ into two subtrees with $s_i\in \Lambda_i$ for $i=1,2$. By Lemma~\ref{lem:irreducible}, $\supp(r_i)\subset\Lambda_i$ for $i=1,2$. The edge orientation on $\Lambda$ induces edge orientation on $\Lambda_i$ for $i=1,2$. Let $\delta_i$ be the Coxeter element in $A_{\Lambda_i}$ which is compactible with the edge orientation on $\Lambda_i$. As vertices in $\Lambda_1\setminus \{s_1\}$ commute with vertices in $\Lambda_2\setminus\{s_2\}$,  Lemma~\ref{lem:coxeter} implies that $\delta=\delta_2\delta_1$. As $r_i$ is a reflection in $A_{\Lambda_i}$, we know $r_1\leq_L \delta_1$ and $r_2\leq_R \delta_2$ by \cite[Lemma 1.3.3]{bessis}. In particular $\delta_1$ has a minimal reflection decomposition of form $\delta_1=r_1\cdot r'_1\cdot r'_2\cdots r'_k$, and $\delta_2$ has a minimal reflection decomposition of form $\delta_2= r''_1\cdots r''_m\cdot r_2$. Thus 
$$
r''_1\cdots r''_m\cdot r_2\cdot r_1\cdot r'_1\cdot r'_2\cdots r'_k=\delta.
$$
By \cite[Lemma 3.7]{mccammond_failure_lattice_property}, $r_2\cdot r_1\leq_L \delta$. Thus $r_2\cdot r_1, r_1\cdot r_2\in [1,\delta]$, and these two elements are both common upper bound for $r_1$ and $r_2$ with respect to $\le_L$. Then $r_2r_1=r_1r_2$ as $([1,\delta],\le_L)$ is a lattice. We write $r_i$ as an reducible word $w_i$ in $W_S$. Then $w_i$ only uses from letters from $\Lambda_i$, and $w_1w_2,w_2w_1$ are reduced words.
Then by Tits's solution to the word problem of Coxeter group, we know that it is possible to apply the relators finitely many times to transform $w_1w_2$ into $w_2w_1$. However, as $s_2$ is on the right side of $s_1$ in $w_1w_2$ and $m(s_1,s_2)\ge 3$, and the property of having at least one $s_2$ on the right side of $s_1$ is preserved under applying the relations, this leads to a contradiction.
\ep

Given a 4-cycle $\omega\subset \Gamma$ with consecutive vertices $\{x_i\}_{i=1}^4$, a pair of antipodal vertices in $\omega$ means either the pair $\{x_1,x_3\}$, or the pair $\{x_2,x_4\}$. A 4-cycle in $\Gamma$ has a \emph{diagonal} means it has a pair of antipodal vertices of $\omega$ which are connected by an edge in $\Gamma$.  We say the Coxeter presentation graph $\Gamma$ is a \emph{join} of two Coxeter presentation subgraphs $\Gamma_1$ and $\Gamma_2$ if $\Gamma$ is a join of $\Gamma_1$ and $\Gamma_2$ as graphs (i.e. $\Gamma$ is obtained from a disjoint union of $\Gamma_1$ and $\Gamma_2$ by adding a single edge between each vertex of $\Gamma_1$ and each vertex of $\Gamma_2$), and each vertex of $\Gamma_1$ commute with every vertex of $\Gamma_2$. We say a Coxeter presentation graph $\Gamma$ is spherical (resp. of cyclic type) if $A_\Gamma$ is a spherical Artin group (resp. an Artin of cyclic type). 

\bthm
\label{thm_main}
Let $\Gamma$ be a Coxeter presentation graph such that
\bit
\item each complete subgraph of $\Gamma$ is a join of a cyclic type graph and a spherical type graph (we allow one of the join factors to be empty);
\item for any cyclic type induced subgraph $\Lambda\subset \Gamma$, $\Lambda^{\perp}$ is spherical.
\eit
We assume in addition that there exists an orientation of all large edges of $\Gamma$ such that
\bit
\item the orientation restricted to each cyclic type subgraph of $\Gamma$ gives a consistent orientation on the associated circle;
\item if $\omega$ is a 4-cycle in $\Gamma$ with a pair of antipodal points $x_1$ and $x_2$ such that each edge of $\omega$ containing $x_i\in \{x_1,x_2\}$ is either not large or oriented towards $x_i$, then the cycle has a diagonal.
\eit Then $A_\Gamma\times \Z$ is a Garside group. 
\ethm

\bp
Let $S$ be the vertex set of $\Gamma$.
Let $I\subset S$ be a spherical subset. We define $\delta_I$ be a product of all elements in $I$ in an order which is compatible with the orientation of $\Gamma$ in sense explained before Lemma~\ref{lem:coxeter}. Then $\delta_I$ is well-defined by Lemma~\ref{lem:coxeter}.

We also view $\delta_I$ as an element in the Coxeter group $W_\Gamma$. Let $\mathcal S$ be the collection of all spherical subset of $S$. Define $U=\cup_{I\in \mathcal S}[1,\delta_I]$, where $[1,\delta_I]$ denotes the interval in $(W_\Gamma,\le_L)$ with $\le_L$ being the partial order induced considering prefix with respect to minimal reflection decompositions. It is clear that $U$ satisfies the assumptions of Theorem~\ref{thm:homotopy_type_salvetti}.

\mk

We now verify that $U$ satisfies all the requirements in Proposition~\ref{pro:E_lattice}. 

By \cite[Theorem 1.4]{baumeister2014note}, any minimal length reflection decomposition of an element $a\in [1,\delta_I]$ only involves reflections in $W_I$. On the other hand, by \cite[Lemma 1.3.3]{bessis}, for any reflection $r\in W_I$, there exists a minimal length reflection decomposition of $\delta_I$ starting with $r$, thus $r\in U$. Now Assumptions 1 and 2 of Proposition~\ref{pro:E_lattice} follow.

Recall that $R_U$ denotes the collection of reflections in $U$.
For Assumption 3 of Proposition~\ref{pro:E_lattice}, if $r_1,r_2\in R_U$ has a left common upper bound $a\in U$, then there exists a spherical subset $I\in S$ such that $a\in [1,\delta_I]$. By \cite[Theorem 1.4]{baumeister2014note}, $r_1,r_2\in W_I$. As in the previous paragraph, we know $r_1,r_2\in [1,\delta_I]$. As $([1,\delta_I],\leq_L)$ is a lattice (Theorem~\ref{thm:bessis} and \cite[Theorem 1.4]{baumeister2014note}), we know $r_1$ and $r_2$ has left join $a$ in $([1,\delta_I],\leq_L)$. By the same argument as in the verification of Assumption 3 in Proposition~\ref{pro:E_lattice}, we know $a$ is also the left join of $r_1$ and $r_2$ in $(U,\leq_L)$. Similarly, we know if $r_1,r_2\in R_U$ have a common right upper bound in $(U,\leq_R)$, then they have a right join in $(U,\leq_R)$.

Now we verify Assumption 4 of Proposition~\ref{pro:E_lattice}. 

For any $w\in W_\Gamma$, let $I_w=\supp(w)$. We claim that if $a,b\in U$ and $a\cdot b\in U$ (recall that $a\cdot b$ means $|ab|=|a|+|b|$ with $|\cdot|$ denotes the reflection length), then $I_{ab}=I_a\cup I_b$. 
Note that $I_{ab}\subset I_a\cup I_b$ is clear. Now let $a=r_1r_2\cdots r_n$ and $b=r'_1r'_2\cdots r'_m$ be minimal length reflection decomposition of $a$ and $b$. By \cite[Theorem 1.4]{baumeister2014note}, $\supp(r_i)\subset I_a$ for each $i$, thus $I_a=\cup_{i=1}^n\supp(r_i)$. Similarly $I_b=\cup_{i=1}^m\supp(r'_i)$.
As $a\cdot b\in U$, $r_1\cdots r_nr'_1\cdots r'_m$ is a minimal length reflection decomposition of $a$ and $b$. As $ab\in A_{I_{ab}}$, we know from \cite[Theorem 1.4]{baumeister2014note} that $r_i,r'_i\in A_{I_{ab}}$. 
Similarly as before $I_{ab}=(\cup_{i=1}^n\supp(r_i))\cup (\cup_{i=1}^m\supp(r'_i))$. Thus $I_a\subset I_{ab}$ and $I_b\subset I_{ab}$. Now the claim follows.

Let $a,u,v,w$ be as in Assumption 4 of Proposition~\ref{pro:E_lattice}. Then $I_a\cup I_u=I_{au}$, which spans a complete subgraph of $\Gamma$. Similarly, $I_a\cup I_v$ spans a complete subgraph of $\Gamma$. By the previous paragraph, if $u\le_L w$ and $v\le_L w$, then $I_u\subset I_w$ and $I_v\subset I_w$. Hence $I_v\cup I_u$ spans a complete subgraph of $\Gamma$. Thus $I=I_v\cup I_u\cup I_a$ spans a complete subgraph of $\Gamma$. Then by the first bullet point of the theorem, $I=I_1\cup I_2$ where $I_1$ is a cyclic type irreducible component of $I$ and $I_2$ is the union of all irreducible spherical components of $I$. Note that $a=a_1\cdot a_2$ for $a_i\in W_{I_i}\cap U$ for $i=1,2$ - to see this, we consider a minimal reflection decomposition of $a$, and by Lemma~\ref{lem:irreducible} and \cite[Theorem 1.4]{baumeister2014note}, each reflection in this decomposition has support contained in either $I_1$ and $I_2$. As elements supported in $I_1$ commute with elements supported in $I_2$, we can group reflections supported in $I_i$ to form $a_i$ for $i=1,2$, and have $a=a_1\cdot a_2$, as well as  $a_i\in U$ for $i=1,2$. Lemma~\ref{lem:irreducible} also implies that 
$u$ belongs to either $W_{I_1}$ or $W_{I_2}$, and $v$ belongs to either $W_{I_1}$ or $W_{I_2}$. We consider the following cases.
\begin{enumerate}
	\item If $u,v\in W_{I_1}$, then $a_1u,a_1v\in W_{I_1}$. As $a\cdot u=a_1\cdot a_2\cdot u\in U$, by \cite[Lemma 3.7]{mccammond_failure_lattice_property}, $a_1\cdot u\le_L a\cdot u$. Thus $a_1\cdot u\in U$. Similarly, $a_1\cdot v\in U$. Thus $a_1\cdot v,a_2\cdot u\in U\cap W_{I_1}$. Proposition~\ref{prop_cyclic} implies that $a_1\cdot w\in U\cap W_{I_1}$. Note that $\supp(a_2)$ and $\supp(a_1\cdot w)$ commute, so their union is spherical. As $a_2\in U$ and $a_1\cdot w\in U$, it follows from the definition of $U$ that $a\cdot w=a_2\cdot (a_1\cdot w)$ also belongs to $U$.
	\item If exactly one of $\{u,v\}$, say $u$, is in $W_{I_1}$, then $w=u\cdot v$ as $\supp(u)$ and $\supp(v)$ commute. As $\supp(a\cdot u)$ and $\supp(v)$ commute, their union is spherical. As $a\cdot u\in U$ and $v\in U$, it follows from the definition of $U$ that $a\cdot w=(a\cdot u)\cdot v$ is contained in $U$.
	\item Suppose each of $u,v$ is in $W_{I_2}$. As $a\cdot u\in U$ and $a\cdot u=a_1\cdot (a_2\cdot u)$ with $\supp(a_1)$ and $\supp(a_2\cdot u)$ commuting, by a similar argument as before using Lemma~\ref{lem:irreducible} and \cite[Theorem 1.4]{baumeister2014note}, we know $a_2\cdot u\in U\cap W_{I_2}$. Thus $a_2\cdot u\in [1,\delta_{I_2}]$. Similarly, $a_2\cdot v\in [1,\delta_{I_2}]$. As $([1,\delta_{I_2}],\leq_L)$ is a lattice (Theorem~\ref{thm:bessis} and \cite[Theorem 1.4]{baumeister2014note}), $a_2\cdot v$ and $a_2\cdot u$ have a join in $([1,\delta_{I_2}],\leq_L)$, which must be $a_2\cdot w$ as $w=u\vee_L v$. In particular, $a_2\cdot w\in U$. As $\supp(a_2\cdot w)$ and $\supp(a_1)$ commute, we know 
	$a\cdot w=a_1\cdot (a_2\cdot w)\in U$.
\end{enumerate}

Assumption 5 of Proposition~\ref{pro:E_lattice} can be verified similarly.

Now we verify Assumption 6 of Proposition~\ref{pro:E_lattice}. Let $a,b,u,v,x$ be as in Assumption 6.
As $a\cdot x\cdot u\in U$, by previous discussion we know that $I_a\cup I_x\cup I_u=I_{axu}$. Thus $I_a\cup I_u$ spans a complete subgraph of $\Gamma$. Similarly, $I_a\cup I_v$, $I_b\cup I_u$, $I_b\cup I_v$ span complete graphs of $\Gamma$.

First we consider the case when $I_a\cup I_b$ spans a complete subgraph. By Lemma~\ref{lem:irreducible}, $I_a$ and $I_b$ are irreducible. If $I_a$ commute with $I_b$, then $I_a\cup I_b$ is spherical and $a\cdot b$ is the right join of $a$ and $b$ in $(U,\leq_R)$. It remains to consider the case when $I_a\cup I_b$ is irreducible. By the first bullet point of the theorem, either $I_a\cup I_b$ is spherical, or $I_a\cup I_b$ is contained in an induced subgraph of $\Gamma$ of cyclic type.
If $I_a\cup I_b$ is spherical, then $a$ and $b$ have a right join in $([1,\delta_{I_a\cup I_b}],\leq_L)$ by Theorem~\ref{thm:bessis}, hence in $(U,\leq_L)$. Now suppose $I_a\cup I_b$ is contained in a cyclic type subgraph $\Lambda$ of $\Gamma$. It suffices to consider the case that $I_a\cup I_b=\Lambda$, otherwise $I_a\cup I_b$ is spherical and we finish as before. Note that $I_a\cup I_b\cup I_u$ spans a complete subgraph of $\Gamma$. As $I_u$ is irreducible by Lemma~\ref{lem:irreducible}, thus either $I_u\subset I_a\cup I_b$, or $I_u\subset (I_a\cup I_b)^\perp$ by our assumption on complete subgraphs of $\Gamma$.
Similarly,  either $I_v\subset I_a\cup I_b$, or $I_v\subset (I_a\cup I_b)^\perp$. If both $I_u\subset I_a\cup I_b$ and $I_v\subset I_a\cup I_b$ hold, then $a\cdot u,b\cdot u,a\cdot v,b\cdot v\in U\cap W_{I_a\cup I_b}$ by \cite[Theorem 1.4]{baumeister2014note} and we are reduced to Proposition~\ref{prop_cyclic}. If at least one of the two statements $I_u\subset I_a\cup I_b$ and $I_v\subset I_a\cup I_b$ is false, then $I_u\cup I_v$ is spherical (this uses the second bullet point in the assumption of the theorem), which implies that $u$ and $v$ have a left join in $U$. 

The case when $I_u\cup I_v$ spans a complete subgraph is similar. It remains to consider the case that $I_u\cup I_v$ does not span a complete subgraph of $\Gamma$, and $I_a\cup I_b$ does not span a complete subgraph of $\Gamma$. Now we will show this remaining case actually can not happen, hence finishes the proof.

Suppose $I_u\cup I_v$ is not complete. Take $s_u\in I_u$ and $s_v\in I_v$ such that they are not adjacent in $\Gamma$. We hope to show $I_a\cup I_b$ spans a  complete subgraph of $\Gamma$. Take $s\in I_a$ and $t\in I_b$. If $s\in I_u$, then $s$ and $t$ are adjacent in $\Gamma$ as $I_u\cup I_b$ spans a complete subgraph. Now we assume $s\notin I_u$. Note that $s_u\notin I_a$, otherwise $s_u$ and $s_v$ are adjacent in $\Gamma$. 
As $a\cdot x\cdot u\in U$, we know $a\cdot u\in U$ by \cite[Lemma 3.7]{mccammond_failure_lattice_property}. Now by Lemma~\ref{lem:order}, either $s$ and $s_u$ commute, or there is an oriented edge from $s$ to $s_u$. Similarly, we know this sentence is still true if we replace the ordered pair $(s,s_u)$ in the statement by $(t,s_u)$, $(s,s_v)$ and $(t,s_v)$. Thus by our assumption, the 4-cycle $s,s_u,t,s_v$ in $\Gamma$ must have a diagonal. The diagonal must connect $s$ and $t$, as $s_u$ and $s_v$ are not adjacent. Thus $I_a\cup I_b$ spans a complete subgraph of $\Gamma$.
\ep

The following is a consequence of Theorem~\ref{thm_main} and Corollary~\ref{cor:complex}.
\begin{cor}
	\label{cor:complex1}
Let $\Gamma$ be a Coxeter presentation graph with vertex set $S$ satisfying all the assumptions in Theorem~\ref{thm_main}. For each spherical $T \subset S$, we choose a Coxeter element $w_T \in W_T$ compatible with the orientation of $\Gamma$. Let
$$U = \bigcup_{T \subset S \mbox{ spherical}} [e,w_T].$$
Let $\widehat U$ be the lift of $U$ from $W_S$ to $A_S$ via the  isomorphism between the dual Artin group associated with $A_T$ for each $T\subset S$ spherical and $A_T$ (\cite[Theorem 2.2.5]{bessis}).

Let $X_S$ be the flag complex of the Cayley graph of $A_S$ with generating set $\widehat U$. Then $X_S$ admits an $A_S$-equivariant CUB metric such that each simplex of $X_S$ is equipped with a polyhedral norm as in  \cite{haettel_simplicial_npc}.
\end{cor}

\bibliographystyle{alpha}
\bibliography{bibli}

\newcommand{\etalchar}[1]{$^{#1}$}
\def\polhk#1{\setbox0=\hbox{#1}{\ooalign{\hidewidth
  \lower1.5ex\hbox{`}\hidewidth\crcr\unhbox0}}}
\begin{thebibliography}{BDSW14}

\bibitem[BD19]{behrstock2019combinatorial}
Jason Behrstock and Cornelia Dru{\c{t}}u.
\newblock Combinatorial higher dimensional isoperimetry and divergence.
\newblock {\em Journal of Topology and Analysis}, 11(03):499--534, 2019.

\bibitem[BDSW14]{baumeister2014note}
Barbara Baumeister, Matthew Dyer, Christian Stump, and Patrick Wegener.
\newblock A note on the transitive hurwitz action on decompositions of
  parabolic coxeter elements.
\newblock {\em Proceedings of the American Mathematical Society, Series B},
  1(13):149--154, 2014.

\bibitem[Bes99]{bestvina_artin}
Mladen Bestvina.
\newblock Non-positively curved aspects of {A}rtin groups of finite type.
\newblock {\em Geom. Topol.}, 3:269--302, 1999.

\bibitem[Bes03]{bessis}
David Bessis.
\newblock The dual braid monoid.
\newblock {\em Ann. Sci. \'Ecole Norm. Sup. (4)}, 36(5):647--683, 2003.

\bibitem[Bes06]{bessis_free}
David Bessis.
\newblock A dual braid monoid for the free group.
\newblock {\em J. Algebra}, 302(1):55--69, 2006.

\bibitem[Bes15]{bessis_finite_complex_Kpi1}
David Bessis.
\newblock {Finite complex reflection arrangements are K ($\pi$, 1)}.
\newblock {\em Annals of mathematics}, pages 809--904, 2015.

\bibitem[BEZ90]{bjorner1990hyperplane}
Anders Bj{\"o}rner, Paul~H Edelman, and G{\"u}nter~M Ziegler.
\newblock Hyperplane arrangements with a lattice of regions.
\newblock {\em Discrete \& computational geometry}, 5(3):263--288, 1990.

\bibitem[BH99]{bridson_haefliger}
Martin~R. Bridson and Andr\'e Haefliger.
\newblock {\em {Metric \,spaces \,of \,non-positive \,curva\-ture}}, volume
  {319\!} of {\em {Grund.~math.~Wiss.}}
\newblock {Springer}, 1999.

\bibitem[BKL98]{birman_ko_lee}
Joan Birman, Ki~Hyoung Ko, and Sang~Jin Lee.
\newblock A new approach to the word and conjugacy problems in the braid
  groups.
\newblock {\em Adv. Math.}, 139(2):322--353, 1998.

\bibitem[BLR08]{bartels_luck_holger_farrell_jones}
Arthur Bartels, Wolfgang L\"{u}ck, and Holger Reich.
\newblock On the {F}arrell-{J}ones conjecture and its applications.
\newblock {\em J. Topol.}, 1(1):57--86, 2008.

\bibitem[BM00]{brady2000three}
Thomas Brady and Jonathan~P McCammond.
\newblock Three-generator {A}rtin groups of large type are biautomatic.
\newblock {\em Journal of Pure and Applied Algebra}, 151(1):1--9, 2000.

\bibitem[BM10]{brady_mccammond}
Tom Brady and Jon McCammond.
\newblock Braids, posets and orthoschemes.
\newblock {\em Algebr. Geom. Topol.}, 10(4):2277--2314, 2010.

\bibitem[Bou02]{bourbaki_lie_456}
Nicolas Bourbaki.
\newblock {\em Lie groups and {L}ie algebras. {C}hapters 4--6}.
\newblock Elements of Mathematics (Berlin). Springer-Verlag, Berlin, 2002.
\newblock Translated from the 1968 French original by Andrew Pressley.

\bibitem[BS72]{brieskorn_saito}
Egbert Brieskorn and Kyoji Saito.
\newblock Artin-{G}ruppen und {C}oxeter-{G}ruppen.
\newblock {\em Invent. Math.}, 17:245--271, 1972.

\bibitem[CCG{\etalchar{+}}20]{helly_groups}
J{\'e}r{\'e}mie Chalopin, Victor Chepoi, Anthony Genevois, Hiroshi Hirai, and
  Damian Osajda.
\newblock Helly groups.
\newblock {\em arXiv preprint arXiv:2002.06895}, 2020.

\bibitem[CCHO21]{chalopin_chepoi_hirai_osajda}
J{\'e}r{\'e}mie Chalopin, Victor Chepoi, Hiroshi Hirai, and Damian Osajda.
\newblock Weakly modular graphs and nonpositive curvature.
\newblock {\em Mem. Amer. Math. Soc.}, 2021.

\bibitem[CD95]{charney_davis_kpi1}
Ruth Charney and Michael~W. Davis.
\newblock The {$K(\pi,1)$}-problem for hyperplane complements associated to
  infinite reflection groups.
\newblock {\em J. Amer. Math. Soc.}, 8(3):597--627, 1995.

\bibitem[Cha]{charney_problems}
Ruth Charney.
\newblock {Problems related to Artin groups}.
\newblock {American Institute of Mathematics, http://people.brandeis.edu/\tild
  charney/papers/\verb|_|probs.pdf}.

\bibitem[Cha92a]{charney1992artin}
Ruth Charney.
\newblock Artin groups of finite type are biautomatic.
\newblock {\em Mathematische Annalen}, 292(1):671--683, 1992.

\bibitem[Cha92b]{charney_biautomatic}
Ruth Charney.
\newblock Artin groups of finite type are biautomatic.
\newblock {\em Math. Ann.}, 292(4):671--683, 1992.

\bibitem[Cho10]{chouraqui2010garside}
Fabienne Chouraqui.
\newblock Garside groups and {Y}ang--{B}axter equation.
\newblock {\em Communications in Algebra}, 38(12):4441--4460, 2010.

\bibitem[CLL15]{corran2015braid}
Ruth Corran, Eon-Kyung Lee, and Sang-Jin Lee.
\newblock Braid groups of imprimitive complex reflection groups.
\newblock {\em Journal of Algebra}, 427:387--425, 2015.

\bibitem[CMW04]{charney2004bestvina}
Ruth Charney, John Meier, and Kim Whittlesey.
\newblock Bestvina's normal form complex and the homology of {G}arside groups.
\newblock {\em Geometriae Dedicata}, 105(1):171--188, 2004.

\bibitem[CP05]{crisp2005representations}
John Crisp and Luis Paris.
\newblock Representations of the braid group by automorphisms of groups,
  invariants of links, and {G}arside groups.
\newblock {\em Pacific journal of mathematics}, 221(1):1--27, 2005.

\bibitem[CP11]{corran2011new}
Ruth Corran and Matthieu Picantin.
\newblock A new {G}arside structure for the braid groups of type (e, e, r).
\newblock {\em Journal of the London Mathematical Society}, 84(3):689--711,
  2011.

\bibitem[Dav15]{davis_coxeter}
Michael~W. Davis.
\newblock The geometry and topology of {C}oxeter groups.
\newblock In {\em Introduction to modern mathematics}, volume~33 of {\em Adv.
  Lect. Math. (ALM)}, pages 129--142. Int. Press, Somerville, MA, 2015.

\bibitem[Deh02]{dehornoy2002groupes}
Patrick Dehornoy.
\newblock Groupes de {G}arside.
\newblock In {\em Annales scientifiques de l'Ecole normale superieure},
  volume~35, pages 267--306. Elsevier, 2002.

\bibitem[Deh15]{garside}
Patrick Dehornoy.
\newblock {\em Foundations of {G}arside theory}, volume~22 of {\em EMS Tracts
  in Mathematics}.
\newblock European Mathematical Society (EMS), Z\"urich, 2015.
\newblock With Fran\c{c}ois Digne, Eddy Godelle, Daan Krammer and Jean Michel,
  Contributor name on title page: Daan Kramer.

\bibitem[Del72]{deligne}
Pierre Deligne.
\newblock Les immeubles des groupes de tresses g\'en\'eralis\'es.
\newblock {\em Invent. Math.}, 17:273--302, 1972.

\bibitem[Dig06]{digne_garside_An}
F.~Digne.
\newblock Pr{\'e}sentations duales des groupes de tresses de type affine {A}.
\newblock {\em Commentarii Mathematici Helvetici}, 81(1):23--47, 2006.

\bibitem[Dig12]{digne_garside_Cn}
F.~Digne.
\newblock A {G}arside presentation for {A}rtin-{T}its groups of type
  {$\tilde{C_n}$}.
\newblock {\em Ann. Inst. Fourier (Grenoble)}, 62(2):641--666, 2012.

\bibitem[DL03]{dehornoy2003homology}
Patrick Dehornoy and Yves Lafont.
\newblock Homology of gaussian groups.
\newblock In {\em Annales de l'institut Fourier}, volume~53, pages 489--540,
  2003.

\bibitem[DL15]{descombes_lang_hyperbolicity}
Dominic Descombes and Urs Lang.
\newblock Convex geodesic bicombings and hyperbolicity.
\newblock {\em Geom. Dedicata}, 177:367--384, 2015.

\bibitem[DL16]{descombes_lang_flats}
Dominic Descombes and Urs Lang.
\newblock Flats in spaces with convex geodesic bicombings.
\newblock {\em Anal. Geom. Metr. Spaces}, 4(1):68--84, 2016.

\bibitem[DP99]{dehornoy_paris_gaussian}
Patrick Dehornoy and Luis Paris.
\newblock Gaussian groups and {G}arside groups, two generalisations of {A}rtin
  groups{A}rtin groups.
\newblock {\em Proc. London Math. Soc. (3)}, 79(3):569--604, 1999.

\bibitem[DPS22]{delucchi_paolini_salvetti_rank_3}
Emanuele Delucchi, Giovanni Paolini, and Mario Salvetti.
\newblock {Dual structures on Coxeter and Artin groups of rank three}.
\newblock {\em arXiv preprint arXiv:2206.14518}, 2022.

\bibitem[FO20]{fukaya_oguni}
Tomohiro Fukaya and Shin-ichi Oguni.
\newblock A coarse {C}artan-{H}adamard theorem with application to the coarse
  {B}aum-{C}onnes conjecture.
\newblock {\em J. Topol. Anal.}, 12(3):857--895, 2020.

\bibitem[Gar69]{garside1969braid}
Frank~A Garside.
\newblock The braid group and other groups.
\newblock {\em The Quarterly Journal of Mathematics}, 20(1):235--254, 1969.

\bibitem[GP12]{godelle_paris}
Eddy Godelle and Luis Paris.
\newblock Basic questions on {A}rtin-{T}its groups.
\newblock In {\em Configuration spaces}, volume~14 of {\em CRM Series}, pages
  299--311. Ed. Norm., Pisa, 2012.

\bibitem[GT16]{gebhardt2016zappa}
Volker Gebhardt and Stephen Tawn.
\newblock Zappa--sz{\'e}p products of {G}arside monoids.
\newblock {\em Mathematische Zeitschrift}, 282(1-2):341--369, 2016.

\bibitem[Hae22]{haettel_simplicial_npc}
Thomas Haettel.
\newblock {A link condition for simplicial complexes, and CUB spaces}.
\newblock {\em arXiv preprint arXiv:2211.07857}, 2022.

\bibitem[Hae24]{haettel_helly_kpi1}
Thomas Haettel.
\newblock {Lattices, injective metrics and the $K(\pi,1)$ conjecture}.
\newblock {\em Algebraic \& Geometric Topology}, 2024.
\newblock To appear.

\bibitem[Hat02]{hatcher2002algebraic}
Allen Hatcher.
\newblock {\em {Algebraic Topology}}.
\newblock Cambridge University Press, 2002.

\bibitem[HH24]{haettel_huang_weakly_modular}
Thomas Haettel and Jingyin Huang.
\newblock {Lattices, Garside structures and weakly modular graphs}.
\newblock {\em Journal of Algebra}, 656:226--4258, 2024.

\bibitem[HO21a]{haettel_osajda_locally_elliptic}
Thomas Haettel and Damian Osajda.
\newblock {Locally elliptic actions, torsion groups, and nonpositively curved
  complexes}.
\newblock {\em arXiv:2110.12431}, 2021.

\bibitem[HO21b]{huang_osajda_helly}
Jingyin Huang and Damian Osajda.
\newblock Helly meets {G}arside and {A}rtin.
\newblock {\em Invent. Math.}, 225(2):395--426, 2021.

\bibitem[HR15]{harlander2015aspherical}
Jens Harlander and Stephan Rosebrock.
\newblock Aspherical word labeled oriented graphs and cyclically presented
  groups.
\newblock {\em Journal of Knot Theory and Its Ramifications}, 24(05):1550025,
  2015.

\bibitem[Hua23]{huangfourcycle}
Jingyin Huang.
\newblock Labled four cycles and the ${K}(\pi,1)$ problem for {A}rtin groups.
\newblock {\em arXiv preprint arXiv:2305.16847}, 2023.

\bibitem[Hum90]{humphreys1990reflection}
James~E Humphreys.
\newblock {\em Reflection groups and Coxeter groups}.
\newblock Number~29. Cambridge university press, 1990.

\bibitem[JS06]{januszkiewicz2006simplicial}
Tadeusz Januszkiewicz and Jacek \'{S}wi\c{a}tkowski.
\newblock Simplicial nonpositive curvature.
\newblock {\em Publications Math{\'e}matiques de l'Institut des Hautes
  {\'E}tudes Scientifiques}, 104(1):1--85, 2006.

\bibitem[JS23]{jankiewicz_schreve_center_kpi1}
Kasia Jankiewicz and Kevin Schreve.
\newblock {The $K(\pi,1)$-conjecture implies the center conjecture for Artin
  groups}.
\newblock {\em J. Algebra}, 615(1):455--463, 2023.

\bibitem[KIP02]{kent2002geometric}
Richard~P Kent~IV and David Peifer.
\newblock A geometric and algebraic description of annular braid groups.
\newblock {\em International Journal of Algebra and Computation},
  12(01n02):85--97, 2002.

\bibitem[Kno19]{knopf2019acylindrical}
Svenja Knopf.
\newblock Acylindrical actions on trees and the farrell--jones conjecture.
\newblock {\em Groups, Geometry, and Dynamics}, 13(2):633--676, 2019.

\bibitem[KR17]{kasprowski_rueping}
Daniel Kasprowski and Henrik R\"uping.
\newblock The {F}arrell-{J}ones conjecture for hyperbolic and {CAT}(0)-groups
  revisited.
\newblock {\em J. Topol. Anal.}, 9(4):551--569, 2017.

\bibitem[Lan50]{lanner1950complexes}
Folke Lann{\'e}r.
\newblock On complexes with transitive groups of automorphisms.
\newblock {\em (No Title)}, 1950.

\bibitem[Lan13]{lang}
Urs Lang.
\newblock Injective hulls of certain discrete metric spaces and groups.
\newblock {\em J. Topol. Anal.}, 5(3):297--331, 2013.

\bibitem[LL07]{lee_lee_garside_translation}
Eon-Kyung Lee and Sang~Jin Lee.
\newblock Translation numbers in a {G}arside group are rational with uniformly
  bounded denominators.
\newblock {\em J. Pure Appl. Algebra}, 211(3):732--743, 2007.

\bibitem[McC05]{mccammond_intro_garside}
Jon McCammond.
\newblock {An introduction to Garside structures}.
\newblock {\em {https://web.math.ucsb.edu/\tild
  jon.mccammond/papers/intro-garside.pdf}}, 2005.

\bibitem[McC15]{mccammond_failure_lattice_property}
Jon McCammond.
\newblock Dual euclidean {A}rtin groups and the failure of the lattice
  property.
\newblock {\em J. Algebra}, 437:308--343, 2015.

\bibitem[McC17]{mccammond_mysterious}
Jon McCammond.
\newblock {The mysterious geometry of Artin groups}.
\newblock {\em {http://web.math.ucsb.edu/\tild
  jon.mccammond/papers/mysterious-geometry.pdf}}, {2017}.

\bibitem[Mos97]{mosher_biautomatic}
Lee Mosher.
\newblock Central quotients of biautomatic groups.
\newblock {\em Comment. Math. Helv.}, 72(1):16--29, 1997.

\bibitem[MS17]{mccammond_sulway}
Jon McCammond and Robert Sulway.
\newblock Artin groups of {E}uclidean type.
\newblock {\em Invent. Math.}, 210(1):231--282, 2017.

\bibitem[Pao21]{paolini_dual_approach}
Giovanni Paolini.
\newblock {The dual approach to the $K(\pi,1)$ conjecture}.
\newblock {\em arXiv preprint arXiv:2112.05255}, 2021.

\bibitem[Par14]{paris_kpi1}
Luis Paris.
\newblock {$K(\pi,1)$} conjecture for {A}rtin groups.
\newblock {\em Ann. Fac. Sci. Toulouse Math. (6)}, 23(2):361--415, 2014.

\bibitem[Pic22]{picantin2022cyclic}
Matthieu Picantin.
\newblock Cyclic amalgams, {HNN} extensions, and {G}arside one-relator groups.
\newblock {\em Journal of Algebra}, 607:437--465, 2022.

\bibitem[PS21]{paolini_salvetti}
Giovanni Paolini and Mario Salvetti.
\newblock Proof of the {$K(\pi,1)$} conjecture for affine {A}rtin groups.
\newblock {\em Invent. Math.}, 224(2):487--572, 2021.

\bibitem[Qi07]{qi2007note}
Dongwen Qi.
\newblock A note on parabolic subgroups of a coxeter group.
\newblock {\em Expositiones Mathematicae}, 25(1):77--81, 2007.

\bibitem[Sco08]{scott_mock}
Richard Scott.
\newblock Right-angled mock reflection and mock artin groups.
\newblock {\em Transactions of the American Mathematical Society},
  360(8):4189--4210, 2008.

\bibitem[Soe21]{soergel_systolic}
Mireille Soergel.
\newblock Systolic complexes and group presentations.
\newblock {\em arXiv preprint arXiv:2105.01345}, 2021.

\bibitem[Tit66]{tits_artin}
J.~Tits.
\newblock Normalisateurs de tores. {I}. {G}roupes de {C}oxeter \'{e}tendus.
\newblock {\em J. Algebra}, 4:96--116, 1966.

\bibitem[VdL83]{vanderlek}
Harm Van~der Lek.
\newblock {The homotopy type of complex hyperplane complements}.
\newblock {1983}.
\newblock {PhD thesis, Radboud University of Nijmegen}.

\bibitem[Wen03]{wenger2003isoperimetric}
Stefan Wenger.
\newblock Isoperimetric inequalities of {E}uclidean type in metric spaces.
\newblock {\em arXiv preprint math/0306089}, 2003.

\end{thebibliography}

\end{document}